\documentclass[12pt]{iopart}
\expandafter\let\csname equation*\endcsname=\relax
\expandafter\let\csname endequation*\endcsname=\relax
\usepackage{amsmath,amssymb,amsthm}
\usepackage{algorithm}
\usepackage{algorithmicx}
\usepackage[noend]{algpseudocode}
\usepackage{algpseudocode}
\usepackage{graphicx}
\usepackage{epstopdf}
\usepackage{float} 
\usepackage{booktabs} 

\usepackage{subcaption}
\usepackage{amsmath} 

\usepackage[usenames]{color}
\usepackage{amsthm,amsmath,amssymb}
\usepackage{amsgen, amsfonts, amsbsy, amssymb}
\usepackage{iopams}
\usepackage{mathrsfs}
\usepackage{color}
\usepackage{tikz}
\usepackage{lineno}

\usetikzlibrary{shapes,positioning}

\tikzstyle{inout}=[trapezium, trapezium left angle=60, trapezium right angle=120, draw] 
\tikzstyle{end}=[rectangle, rounded corners, draw]   
\tikzstyle{endn}=[rounded rectangle, draw]   
\tikzstyle{exec}=[rectangle, draw]   
\tikzstyle{decide}=[kite, kite vertex angles=120, draw]

\usepackage{hyperref}
\hypersetup{colorlinks=true,
	linkcolor=blue,
	anchorcolor=blue,
	citecolor=blue}
\floatname{algorithm}{Algorithm} 
\newtheorem{theorem}{Theorem}[section]

\newtheorem{assumption}{Assumption}

\newtheorem{condition}{Condition}

\newtheorem{definition}{Definition}

\newtheorem{lemma}{Lemma}[section]

\newtheorem{proposition}{Proposition}[section]

\usepackage{titlesec}
\usepackage{titletoc}
\usepackage[title]{appendix}
\numberwithin{equation}{section}

\begin{document}

\title[Nonconvex weighted variational MAR via convergent primal-dual  algorithms ]{Nonconvex weighted variational metal artifacts removal via convergent primal-dual  algorithms}

\author{Lianfang Wang$^1$, Zhangling Chen$^1$, Zhifang Liu$^1$, Yutong Li$^1$, Yunsong Zhao$^2$, Hongwei Li$^2$, and
 Huibin Chang$^1$\footnote{ Corresponding author.}
 }

\address{$^1$ School of Mathematical Sciences, Tianjin Normal University, Tianjin,  300387, China}

\address{$^2$ School of Mathematical Sciences, Capital Normal University, Beijing, 100048, China}
\ead{changhuibin@gmail.com}

\begin{abstract}
Direct reconstruction through filtered back projection engenders metal artifacts in polychromatic computed tomography images, attributed to highly attenuating implants, which further poses great challenges for subsequent image analysis. Inpainting the metal trace directly  in the Radon domain for the extant  variational method leads to strong edge diffusion and potential inherent artifacts. With normalization based on pre-segmentation, the inpainted outcome can be notably ameliorated. However, its reconstructive fidelity is heavily contingent on the precision of the pre-segmentation, and highly accurate segmentation of images with metal artifacts is non-trivial in actuality. In this paper, we propose a nonconvex weighted variational approach for metal artifact reduction. Specifically, in lieu of employing a binary function with zeros in the metal trace, an adaptive weight function is designed in the Radon domain, with zeros in the overlapping regions of multiple disjoint metals as well as areas of highly attenuated projections, and the inverse square root of the measured projection in other regions. A nonconvex $L^1-\alpha L^2$ regularization term is incorporated to further enhance edge contrast, alongside a box-constraint in the image domain. Efficient first-order primal-dual algorithms, proven to be globally convergent and of low computational cost owing to the closed-form solution of all subproblems, are devised to resolve such a constrained nonconvex model.  Both simulated and real experiments are conducted with comparisons to other variational algorithms, validating the superiority of the presented method. Especially in comparison to Reweighted JSR, our proposed algorithm can curtail the total computational cost to at most one-third, and for the case of inaccurate pre-segmentation, the recovery outcomes by the proposed algorithms are notably enhanced.
\end{abstract}
\maketitle

\section{Introduction}
X-ray computed tomography (CT) is one of the most common  means of medical diagnosis.
However, when metallic implants present, the reconstructed CT images by conventional reconstruction algorithms like the filtered back projection (FBP) \cite{katsevich2002theoretically} and simultaneous algebraic reconstruction technique \cite{andersen1984simultaneous} may suffer from serious metal artifacts, thus potentially engendering misdiagnoses. Consequently, conceiving efficacious methodologies for metallic artifact reduction (MAR) is imperative.

One effective approach is to modify the model by incorporating  imaging physics \cite{park2015metal, lee2019direct, hur2021metal, hegazy2023metal, zhu2023physics}. 
This methodology can mitigate artifacts whilst preserving boundaries and specifics. However, under conditions of elevated noise levels, a divergence transpires betwixt the aforementioned modelling and the authentic model, inevitably engendering unsatisfactory outcomes. 
Another well-known way is to identify the metal-affected projection as missing data and use interpolation to recover it. {\color{blue}Some are to directly inpaint the domain with different interpolation algorithms \cite{kalender1987reduction, bruyant2000streak, roeske2003reduction, gu2006method, veldkamp2010development}, leading to inaccurate boundaries and the introduction of new artifacts.} In order to enhance the performance \cite{meyer2010normalized, axente2015clinical}, Meyer et al. \cite{meyer2010normalized} later  proposed a normalized metal artifact reduction (NMAR), which utilized  prior information to normalize the measurement and then interpolated the sinogram of the metal region. Interpolation alone cannot precisely restore the true missing projection, and thus still introduce new artifacts. With the remarkable success of deep learning in medical image processing, recent studies have implemented deep neural networks (DNNs) to tackle the problem of reducing metal artifacts. {\color{blue}The existing researches consist of the image-to-image learning \cite{huang2018metal, liao2019adn, gjesteby2019dual, nielsen2019magnetic, li2023marganvac}, the sinogram domain network \cite{gjesteby2017deep, zhang2018convolutional, ghani2019fast, yu2020deep, trapp2022empirical} and dual domain (both the image and sinogram domains) \cite{lyu2020encoding, ikuta2022deep, zhou2022dudodrdual, hyun2022deep, wang2023indudonet+dual} that utilizes either residual learning or adversarial learning techniques.} However, DNNs commonly require large, representative training datasets and extensive computational resources, imposing practical limitations on their widespread adoption.

As a critical mathematical technique for MAR, variational regularization methods offer interpretability, stability, and computational efficiency with modest local resource requirements. Since the correction process involves solving a mathematically ill-posed problem, regularization-based methods play a crucial role. {\color{blue}Successful examples include the total variation (TV) \cite{sidky2008image, ritschl2011improved, chen2013limited, zhang2016iterative} and  the wavelet frame–based approach \cite{dong2013x, zhan2016ct, zhang2018reweighted, choi2018pet}. Furthermore,  different variants of TV regularization methods \cite{liu2012adaptiveTV, lou2015weighted, faggiano2016metal, gong2019adaptive, deng2021customized} have been presented  in order to achieve more accurate reconstruction.} Due to the non-uniformity of metal artifacts, modelling solely based on uncontaminated projection information and regularization may not be sufficient for effectively removing metal artifacts. Therefore, many researchers proposed incorporating preprocessed images into the reconstruction process. For instance, Zhang, Dong and Liu \cite{zhang2018reweighted} proposed to normalize the original data in the sinogram domain by projections of pre-segmentation. Then an effective iterative reconstruction was obtained relying on the regularization on dual domains. 

The existing methods  directly inpainted the metal trace (the regions related to the metal in the sinogram domain) \cite{mehranian2013x, zhang2016computed, zhang2016iterative}, which were acknowledged to produce strong diffusion around the metal (See Fig. \ref{metalmask1}), and even new artifacts (See Fig. \ref{metalmask}). 
Other methods further normalized the Radon domain \cite{meyer2010normalized, zhang2018reweighted} based on the projection of a segmentation, which required relatively accurate pre-segmentation. However, it is also challenging to segment the CT images contaminated by metal artifacts. Furthermore, the convex regularization models, either with TV or tight frames, were prone to a reduction of the edge contrast \cite{wu2018general}. In this paper, 
in lieu of employing a binary function with zeros in the metal trace, an adaptive weight function in the Radon domain is designed, with zeros in the overlapping region of multiple separated metals as well as the area of high attenuated projections, and the inverse square root of the measured projection in other regions.   To further enhance the edge contrast, a nonconvex $L^1-\alpha L^2$ regularization term  is considered, such that we develop a non-convex weighted variational method for MAR, together with a box-constraint in the image domain. 
To resolve such a constrained non-convex model, we reformulate the proposed model into saddle-point problems based on the predual forms of TV. We then design efficient first-order primal-dual algorithms  that are proven to be globally convergent under mild conditions.
In both the simulated and real experiments, compared with other existing  methods, such as beam-hardening corrector (BCMAR) \cite{park2015metal}, NMAR \cite{meyer2010normalized}, {\color{blue}inpainting with compound prior modelling both sinogram and image sparsity (TV-TV inpainting) \cite{zhang2016computed} and the reweighted joint spatial-Radon domain (Reweighted JSR )} \cite{zhang2018reweighted}, the proposed algorithms produce higher accuracy reconstruction with low computational cost.

The main contributions of this paper are listed as follows:
\begin{itemize}
	\item We present a novel weighted nonconvex variational model to correct metal artifacts by combining the nonconvex regularization with an adaptive weighted norm. Namely,  the weight function is specially constructed based on the hybrid scheme utilizing the measured projection other than simply setting it to a binary matrix. 
	
	\item We develop a first-order preconditioned primal-dual hybrid gradient algorithm to solve the reformulated penalized saddle-point problem for the proposed model, proving its convergence given proper parameters. Then by introducing an auxiliary variable in the Radon domain, an even faster fully-splitting primal-dual hybrid gradient algorithm is further proposed, now with guaranteed convergence under the additional assumption that the auxiliary variable is bounded. Both proposed algorithms can be efficiently implemented, as each subproblem has a closed-form solution.
	\item 
	We conduct numerous experiments to evaluate the various aspects including convergence, parameter impact and performance for the proposed algorithms. Numerically, the proposed methods can produce competitive results in terms of peak signal-to-noise ratio (PSNR) and structural similarity index measure (SSIM), among all the compared variational methods. Especially, for less precise pre-segmentation, the proposed algorithms can reconstruct significantly superior result in comparison to Reweighted JSR. Moreover, the proposed fully splitting variant algorithm with much faster convergence, can reduce the total computational cost to at most one-third than Reweighted JSR.
\end{itemize}

The remainder of the paper is organized as follows. In section \ref{Pre-knowledge}, we briefly review some of the basic concepts involved. The nonconvex weighted MAR model with a box-constraint is presented, and the efficient splitting algorithms are given in section \ref{proposed model weighted}, with the convergence guarantee in section \ref{convergence guarantee}. In section \ref{numerical experiment}, the effectiveness of the proposed algorithm is validated by numerous experiments. Section \ref{conclusion} summarizes this work.

\section{Preliminaries}\label{Pre-knowledge}
\subsection{Polychromatic X-ray CT and related variational reconstruction methods}
For a multichromatic energy X-ray, the adoption of the monochromatic energy assumption allows the imaging model to be rewritten as 
\begin{eqnarray}
\label{linear app}
    \mathcal P u \approx Y,
\end{eqnarray}
where $Y$ is the measured projection data, $\mathcal P: X \to \mathbb{R}^{m_1 \times m_2}$ is the Radon transform representing the discrete line integrals at $m_1$ different projection angles and along a total of $m_2$ different beams and $u \in X=\mathbb{R}^{n \times n}$ denotes the target image at a specific but unknown energy level. Components with high attenuation have greater energy dependence, such as metal implants, and therefore, solving \eqref{linear app} leads to severe artifacts that greatly reduce the image quality.

By introducing TV regularization in the projection domain, one can  treat the metal trace as missing data and try to recover it  via the following model \cite{zhang2016computed}
\begin{eqnarray}
\label{mask proj model}
   \mathop {\min }\limits_f \frac{1}{2}{\big\| {{B_\Omega^c} \odot \big( {f - Y} \big)} \big\|^2} + \lambda {\big\| \nabla f \big\|_{1}},
\end{eqnarray}
where for all $1\leq i\leq m_1, 1\leq j\leq m_2$, ${(B_\Omega^c)_{i,j}}:=1-(B_\Omega)_{i,j}$ with $B_\Omega$ being a binary matrix related to metal trace $\Omega$ defined as follows
\begin{eqnarray*}
(B_{\Omega})_{i,j}:= \left\{
\begin{aligned}
& 0 \qquad {\rm{if}} \, (i,j) \in \Omega,\\
& 1  \qquad {\rm{otherwise,}}
\end{aligned}
\right .
\end{eqnarray*}   
the notations $\odot$, $ \| \cdot \|$ and $\nabla$ denote the Hadamard product, the standard Frobenius  norm in  $\mathbb R^{m_1 \times m_2}$ and the gradient operator, respectively, and the norm $\|\cdot\|_1$ is defined as:
\begin{eqnarray*}
    {\big\| p \big\|_1} := \sum\limits_{1\le i,j \le n} {{\big| {{{\big( {{p_x}} \big)}_{i,j}}} \big| + \big| {{{\big( {{p_y}} \big)}_{i,j}}} \big|}} ~~~\forall  p=(p_x,p_y) \in X \times X.
\end{eqnarray*}
After arriving at the projection data by solving the above model, the final reconstruction is immediately derived by FBP. 
 Similarly, one can directly consider the sparse restoration in the image domain via the following TV-MAR \cite{wang2023indudonet+dual} model
\begin{eqnarray}
\label{mask model}
 \mathop {\min }\limits_u \frac{1}{2}{\big\| {{B_\Omega^c} \odot \big( {\mathcal P u - Y} \big)} \big\|^2} + \lambda \| \nabla u \big\|_{2,1}.
\end{eqnarray}
where the norm $\| \cdot \|_{2,1}$ is defined as
\begin{eqnarray*}
    {\big\| p \big\|_{2,1}} := \sum\limits_{1\le i,j \le n}{{\sqrt {{{\big| {{{\big( {{p_x}} \big)}_{i,j}}} \big|}^2} + {{\big| {{{\big( {{p_y}} \big)}_{i,j}}} \big|}^2}} } }~~ ~ \forall p \in X \times X. 
\end{eqnarray*}
They have good performance in suppressing artifacts and noises, but the textures in the image may be erased.

Subsequently, in order to balance detail preservation and artifact reduction, Zhang, Dong and Liu \cite{zhang2018reweighted} considered a dual domain model
\begin{eqnarray*}
  \begin{aligned}
      \mathop {\min }\limits_{u,f} \frac{1}{2}{\big\| {\mathcal P u - {Y_s}f} \big\|^2} &+ {\big\| {{\boldsymbol{\lambda}_1} \cdot {W_1}u} \big\|_{1,2}} + {\big\| {{\boldsymbol{\lambda} _2} \cdot {W_2}f} \big\|_{1,2}}\\
      &+ \frac{1}{2}{\big\| {{B_\Omega^c} \odot \big( {f - {Y}/{{{Y_s}}}} \big)} \big\|^2},
  \end{aligned}  
\end{eqnarray*}
where $Y_s$ is the Radon transform of a pre-segmentation based on the level set method (the level set method \cite{li2011level} is used to obtain segmentation including low-density components (soft tissues) and high-density components such as bones and metals), 
\begin{eqnarray*}
{\big\| {{\boldsymbol{\lambda} _i} \cdot {W_i}u} \big\|_{1,2}} = {\Big\| {\sum\limits_s {{{\big( {\sum\limits_l {{\boldsymbol{\lambda} _i}_{( {l,s} )}{{\big| {{{\big( {{W_i}u} \big)}_{l,s}}} \big|}^2}} } \big)}^{\frac{1}{2}}}} } \Big\|_1} \qquad\forall i=1,2   
\end{eqnarray*}
with ${\boldsymbol{\lambda}}_i> 0$ and $W_i,\,i=1,2$ being tight wavelet frames, and ``$/$'' denotes dot division of two matrices. Although the MAR results  are greatly improved, such a method highly relies on the accuracy of the pre-segmentation.

\subsection{Anisotropic-isotropic regularization}
Variational methods with TV  \cite{Pragliola2023on} can preserve edge information for  piecewise-constant images. In order to enhance the sparsity, one may consider the $L^0$ pseudo-norm of the image gradient ${\| {\nabla u} \|_0}$. However, directly solving the minimization using such a norm  is  NP-hard, and its nonconvex approximation is a better alternative. A typical example is  the weighted anisotropic and isotropic TV (AITV) 
\begin{eqnarray*}
\mathrm{AITV}(u):={\big\| {\nabla u} \big\|_1} - \alpha {\big\| {\nabla u} \big\|_{2,1}}.    
\end{eqnarray*}
It has been widely used in image processing \cite{lou2015weighted, park2016weighted, wang2021limited}, which produces recovery results with better  contrast and fruitful details compared to traditional TV. 

\section{ The proposed  variational method for MAR}\label{proposed model weighted}
In this section, we propose a new nonconvex weighted model with a box-constraint to reduce metal artifacts produced by polychromatic X-ray imaging systems and develop effective convergent algorithms. Since the pixel value of the digital image is finite, the box constraint is reasonable. Especially for CT images, the pixel value represents the attenuation coefficient in a particular case \cite{kak2001principles} in practice. In fact, the box-constraint has also been widely used in image denoising \cite{chan2012multiplicative},  image classification \cite{kan2021pnkh} and CT reconstruction \cite{wang2021limited}.
\subsection{ A nonconvex weighted MAR with a box-constraint}\label{subsection1}
A non-convex weighted model with a box-constraint can be given below 
\begin{eqnarray}\label{equ1}
	\mathop {\min }\limits_u \frac{1}{{2\lambda }}\big\| {W \odot \big( {\mathcal P u - Y} \big)} \big\|^2 + {\big\| {\nabla u} \big\|_1} - \alpha {\big\| {\nabla u} \big\|_{2,1}}\qquad {\rm{s.t.}}\,\,u \in \left[ {0,c} \right ],
\end{eqnarray}
where the constant $c$ is the upper bound of the reconstructed image, $\lambda$ and $\alpha$ are two positive parameters, 
\begin{eqnarray}
    \label{weight W}
   W := \frac{{\boldmath 1}_{m_1,m_2}}{\max\{Y^{\frac{1}{2}}, \varepsilon {\boldmath 1}_{m_1,m_2}\}} \odot B_{\Omega_t}, 
\end{eqnarray}
where ${\boldmath 1}_{m_1, m_2}\in\mathbb R^{m_1\times m_2}$ denotes the all-ones matrix, $\frac{~\cdot~}{~\cdot~}$, and $\max\{\cdot, \cdot\}$ denote the elementwise division and maximum of two matrices, respectively.
As expressed in (\ref{weight W}), the parameter $\varepsilon$  serves to preclude division by zero.

\begin{figure}[htbp]
     \centering
	\begin{subfigure}{0.45\linewidth}
		\centering		\includegraphics[width=0.9\linewidth]{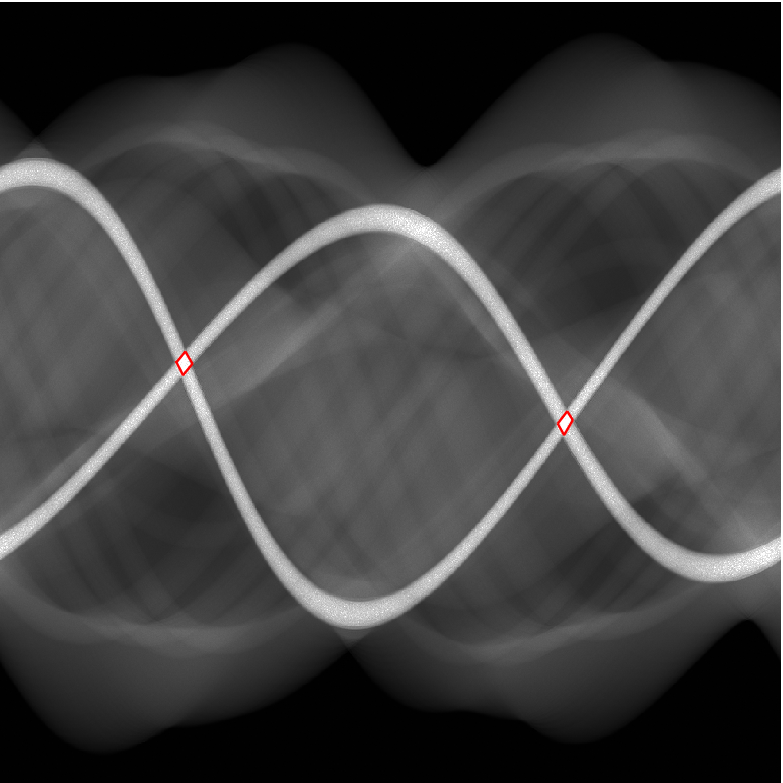}
		\caption{ }
		\label{Y}   
	\end{subfigure}
	\centering
	\begin{subfigure}{0.45\linewidth}
		\centering
		\includegraphics[width=0.9\linewidth]{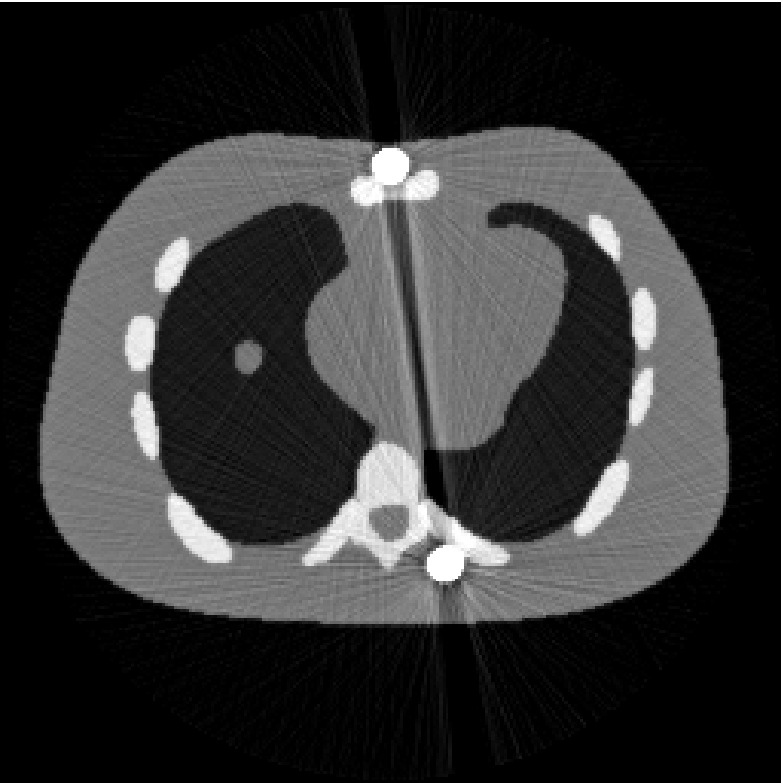}
		\caption{ }
		\label{maCT}   
	\end{subfigure}
 
 	\begin{subfigure}{0.45\linewidth}
		\centering		\includegraphics[width=0.9\linewidth]{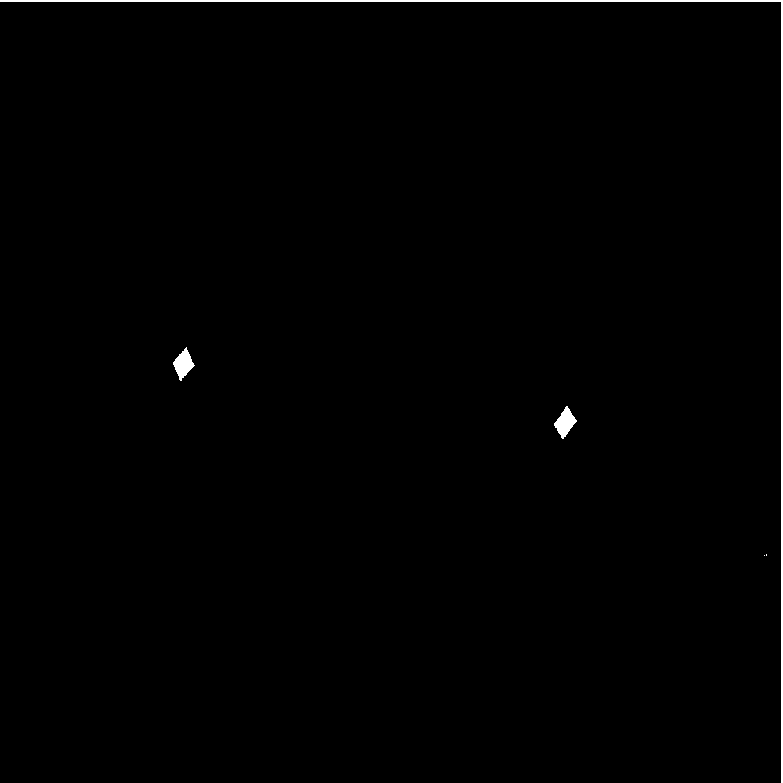}
		\caption{ }
		\label{mask}   
	\end{subfigure}
	\centering
	\begin{subfigure}{0.45\linewidth}
		\centering
		\includegraphics[width=0.9\linewidth]{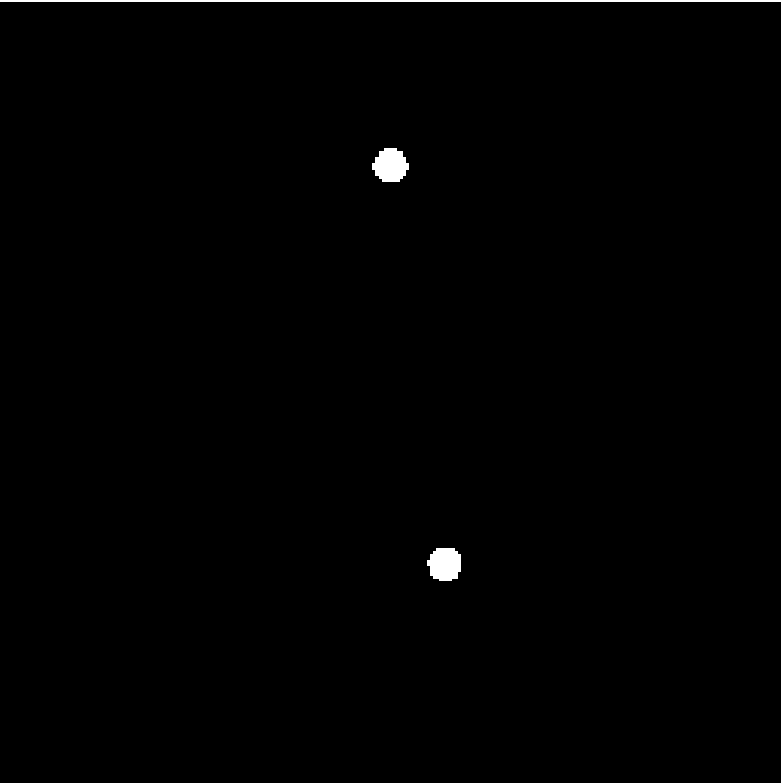}
		\caption{ }
		\label{metal}   
	\end{subfigure}
 \label{ncat ncat}
 \caption{(a) The measured projection data. (b) The reconstructed image $u_a$ by the analysis model (\ref{CG reconstruct}) (4600 Hounsfield Units (HU) window, 1300 HU level). (c) The overlap area $\Omega_t$. (d) The metal image.}
\end{figure}

Here
\begin{eqnarray*}
 \Omega_t:= O_m\cup O_t   
\end{eqnarray*}
is a subset of the metal trace $\Omega$, where 
$O_m$ represents the projection region jointly through each two separated metals (as shown in the red box area in Figure \ref{Y}), and the region for  highly attenuated projections 
\begin{eqnarray*}
O_t:=\left\{ {(i,j)\in\Omega~|~{\mkern 1mu} {Y_{i,j}} \ge t \times \big\| Y\big\|_{\rm{max}} ~~\forall 1\le i\le m_1,1\le j\le m_2} \right\}    
\end{eqnarray*}
with thresholding level $t$ and $\| Y \|_{\rm{max}} := \mathop {\max }\limits_{(\hat{i},\hat{j}) \in \{1,2,\cdots, m_1\}\times \{1,2,\cdots, m_2\}}  |Y_{\hat i,\hat j}|$.

 In (\ref{equ1}), the weighted fidelity term $W$ makes the recovery mechanism more efficient. If only the non-metallic information of the projection data is considered (i.e. use $B_\Omega^c$ as the weighted matrix instead of $W$) as \eqref{mask proj model}, the boundary in the image domain may be diffused, and  new artifacts also appear (as shown in Figure \ref{metalmask}). That may be caused by discarding too many projections in the metal trace. Hence, we here propose to only discard the most severely contaminated projections, which are either the intersection of two metals as $O_m$ or the region $O_t$ related to other relatively high attenuation materials (e.g. bones). Moreover, an inverse square root of the measured projections is combined as \eqref{weight W} in order to adaptively balance the error distribution, i.e. bigger errors in the Radon domain are allowed for the higher attenuated projections.

In addition, we show how to determine the region $O_m$. A rough reconstruction is obtained by  FBP or using the conjugate gradient (CG) method by solving the following least-square problem 
\begin{eqnarray}
    \label{CG reconstruct}
    \mathop {\min }\limits_u \frac{1}{2}{\big\| {\mathcal P u - Y} \big\|^2}.
\end{eqnarray}
The reconstructed image, denoted by $u_a$, is shown in Figure \ref{maCT}. And then the position of the metal in the image domain can be obtained by simple threshold processing. Further the different metals are projected separately to get the overlap area $O_m$. 

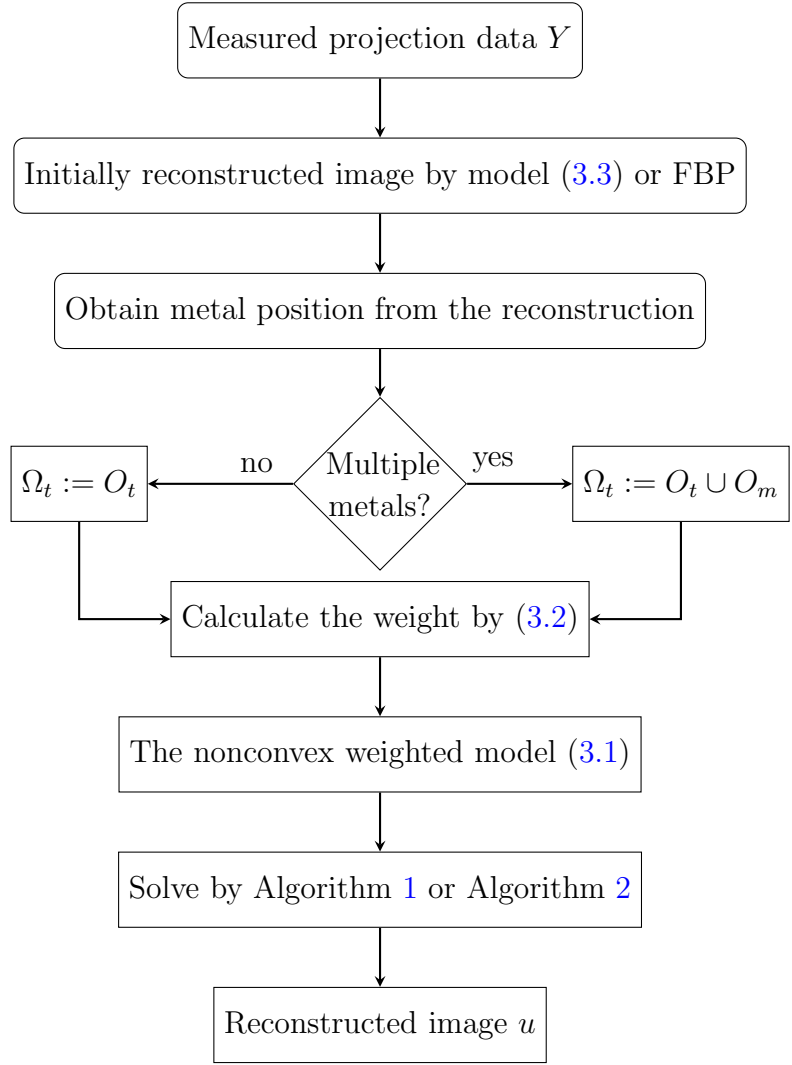
\begin{figure}
\begin{center}
\tikzstyle{startstop} = [rectangle, rounded corners, minimum width=1cm, minimum height=1cm,text centered, draw=black]
\tikzstyle{decision} = [diamond, draw, text width=3.5em, text badly centered, inner sep=0pt]
\tikzstyle{process} = [rectangle, minimum width=1cm, minimum height=1cm, draw=black]
\tikzstyle{arrow} = [thick,->,>=stealth]
\begin{tikzpicture}[node distance = 1.8cm]
\node (start) [startstop]at (-9,5) {Measured projection data $Y$};
\node (inout) [startstop, below of = start] {Initially reconstructed image by model (\ref{CG reconstruct}) or FBP};
\node (metal) [startstop, below of = inout] {Obtain metal position from the reconstruction};
\node (decP) [decision, below of = metal, yshift=-0.5cm] {Multiple metals?};
\node (one) [process, below of=metal, xshift = -4cm,yshift=-0.5cm] {$\Omega_t:= O_t$};
\node (multi) [process, below of = metal, xshift = 4cm,yshift=-0.5cm]{$\Omega_t:= O_t \cup O_m$};

\node (weight) [process, below of=decP ] { 
  Calculate the weight by \eqref{weight W} };
\node (compute) [process, below of=weight] {The nonconvex weighted model (\ref{equ1})};
\node (solve) [process, below of=compute] {Solve by Algorithm \ref{Pre Pdhg} or Algorithm \ref{Pdhg}};
\node (result) [process, below of=solve] {Reconstructed image $u$};
\draw [arrow] (start) -- (inout);
\draw [arrow] (inout) -- (metal);
\draw [arrow] (metal) -- (decP);
\draw [arrow] (decP) -- node[near start, above] {yes} (multi);
\draw [arrow] (decP) -- node[near start, above] {no} (one);

\draw [arrow](multi) |- (weight);
\draw [arrow](one) |- (weight);
\draw [arrow] (weight) -- (compute);
\draw [arrow] (compute) -- (solve);
\draw [arrow] (solve) -- (result);
\end{tikzpicture}
\end{center}
\caption{Flowchart of the proposed CT image reconstruction with reduced metal artifacts.}
\label{Flowchart}
\end{figure}

To better demonstrate the proposed CT image reconstruction method for reducing metal artifacts, we summarize the entire procedure in the flowchart shown in Figure \ref{Flowchart}.

\subsection{ Primal-dual optimization for the penalized model}\label{subsection2}
Due to the nonconvex regularization, (\ref{equ1}) is a nonconvex optimization problem. A classical approach for solving such non-convex model is to use difference-of-convex algorithm (DCA) \cite{lou2015weighted, park2016weighted, lou2018fast}. Writing it as the difference of two convex models, and linearizing the convex one containing ${\| \cdot \|_{2,1}}$ term, one can solve it using the split Bregman \cite{goldstein2009split} or primal-dual \cite{zhu2008efficient, esser2010general, chambolle2011first} technique.  Here we consider  a primal-dual splitting algorithm based on the predual form of TV. 

Based on the predual form of  the ${\|  \cdot \|_{2,1}}$, the corresponding saddle-point problem  can be written as
\begin{eqnarray*}
\mathop {\min }\limits_{u } \frac{1}{{2\lambda }}\left\| {W \odot \left( {\mathcal P u - Y} \right)} \right\|^2 + {\mathbb{I} _U}\left( u \right)+ {\left\| {\nabla u} \right\|_1} - \alpha\mathop {\max }\limits_{\left\| q \right\|_{2,\infty} \le 1}-  \left\langle {\nabla u,q} \right\rangle, 
\end{eqnarray*}
where ${\mathbb{I}_U}( u )$ is the indicator function regarding a closed set $U$ :
\begin{eqnarray*}
    {\mathbb{I} _U}\left( u \right) :=
    \begin{cases}
     0 \qquad  &\qquad {\rm{if}} \,u \in U,\\
	+ \infty  &\qquad {\rm{otherwise}}   
    \end{cases}
\end{eqnarray*}
with $U := \left\{ {u \in X:0 \le {u_{i,j}} \le c ~~\forall i,j} \right\}$, ${\big\| q \big\|_{2,\infty }}: = \mathop {\max }\limits_{1 \le i,j \le n} \sqrt {{{\| {{{( {{q_x}} )}_{i,j}}} \|}^2} + {{\| {{{( {{q_y}} )}_{i,j}}} \|}^2}}$ and $\langle { \cdot , \cdot } \rangle$ denotes the inner product. It can be rewritten equivalently as the following optimization problem
\begin{eqnarray*}
\mathop {\min }\limits_{ u ,q} \frac{1}{{2\lambda }}\left\| {W \odot \left( {\mathcal P u - Y} \right)} \right\|^2 + {\mathbb{I} _U}\left( u \right)+ {\left\| {\nabla u} \right\|_1} + \alpha \left\langle {\nabla u,q} \right\rangle  + {\mathbb{I} _Q}\left( q \right)
\end{eqnarray*}
with $Q:=\left\{ {q \in X \times X: \sqrt {{{ {({q_x})_{i,j}^2} }} + {{ {({q_y})_{i,j}^2} }}}  \le 1} ~~\forall 
 i,j\right\}$. Further based on the predual form of ${\left\|  \cdot  \right\|_1}$, it can be rewritten as
\begin{equation}
    \label{old saddle point problem}
\mathop {\min }\limits_{ u,q} \mathop {\max }\limits_p \frac{1}{{2\lambda }}\left\| {W \odot \left( {\mathcal P u - Y} \right)} \right\|^2 + {\mathbb{I} _U}\left( u \right)+ \left\langle {\nabla u,p + \alpha q} \right\rangle  -  {\mathbb{I} _S}\left( p \right) + {\mathbb{I} _Q}\left( q \right),
\end{equation}
where $S:=\left\{ {p \in {X \times 
 X},\big| {({p_x})_{i,j}} \big| \le 1,\big| {({p_y})_{i,j}} \big| \le 1} ~~\forall {i,j}\right\}$.
 
 In order to guarantee the  convergence of the splitting  algorithm for such nonconvex optimization model, an additional quadratic  term  is added  to   (\ref{old saddle point problem}), such that one immediately gets the penalized model below
\begin{eqnarray}
\label{non multiplier }
\begin{aligned}
 \mathop {\min }\limits_{u,q} \mathop {\max }\limits_p \frac{1}{{2\lambda }}\left\| {W \odot \left( {\mathcal P u - Y} \right)} \right\|^2+ {\mathbb{I} _U}\left( u \right) + \left\langle {\nabla u,p + \alpha q} \right\rangle& \\
 -  {\mathbb{I} _S}\left( p \right) + {\mathbb{I} _Q}\left( q \right)-\frac{\eta}{2}{\left\| p \right\|^2}&.
\end{aligned}
\end{eqnarray}
The above optimization problem w.r.t. the variable $p$ is strongly concave. 
Note that it can return to the original model (\ref{old saddle point problem}) only if $\eta=0$. We remark that such penalization is to enhance the smoothness of the normalized gradient of the image. 

Letting ${G }( u) := \frac{1}{{2\lambda }}\left\| {W \odot ( {\mathcal P u - Y} )} \right\|^2+{\mathbb{I} _U}( u)  $, $ f(  \cdot  ) := {\| {\cdot} \|_1} $ and $ g( \cdot ) := {\| {\cdot} \|_{2,1}}$, $f_\eta ^*( p ) :={f^*}( p )+\frac{\eta }{2}{\| p \|^2}$ with ${f^*}$ being itself the convex conjugate, the saddle point problem \eqref{non multiplier } is rewritten as follows
\begin{eqnarray}
\label{new saddle-point problem}
 \mathop {\min }\limits_{u ,q} \mathop {\max }\limits_p L_{PD}\left( {u,q,p} \right)
\end{eqnarray}
 with 
 \begin{eqnarray*}
   L_{PD}( {u,q,p} ):={G}( u) + \langle {\nabla u,p+\alpha q} \rangle  - {f_\eta^*}( p )  + g ^*(q ).  
 \end{eqnarray*}
 
\subsubsection{Preconditioned primal-dual hybrid gradient} 
A natural scheme to solve the above saddle point problem is to split them, which consists of following preconditional version \cite{pock2011diagonal} of the $u$ subproblem and three-step iterations for the generalized  primal-dual hybrid gradient (PDHG):
\begin{eqnarray}
\label{Pre iter}
\begin{cases}
  \texttt{Step 1:}\  & {u^{k + 1}} = \arg \mathop {\min }\limits_{u} {L_{PD}}( {u,{q^k},{p^k}} ) + \frac{\| {u - {u^k}} \|_{{\mathcal M_\lambda}}^2}{2},\\ 
 \texttt{Step 2:} \ & {{\bar u}^{k + 1}} = 2{u^{k + 1}} - {u^k},\\
\texttt{Step 3:} \  & {q^{k + 1}} = \arg \mathop {\min }\limits_q {L_{PD}}( {{{\bar u}^{k + 1}},q,{p^k}} ) + \frac{{\| {q - {q^k}} \|^2}}{{2\tau }},\\
\texttt{Step 4:} \ &{p^{k + 1}} = \arg \mathop {\max }\limits_p {L_{PD}}( {{{\bar u}^{k + 1}},{q^{k + 1}},p} ) - \frac{{\| {p - {p^k}} \|^2}}{{2\beta }}.  
\end{cases}
\end{eqnarray}
with positive definite operator ${\mathcal M_\lambda}: X \to X$ (we call it preconditioning operator and give the specific definition later), where $\left\| u \right\|_{{\mathcal M_\lambda}}^2 := \left\langle {{\mathcal M_\lambda}u,u} \right\rangle $. Here we remark that the operator $\mathcal M_\lambda$ is selected such that  not only the $u-$subproblem in Step 1 has a closed form solution, but also the convergence of the Algorithm \ref{Pre Pdhg} can be guaranteed. Besides, since each subproblem is strictly convex, the solutions to the minimization problems are unique and therefore  use ``$=$''.

In order to derive the expression for $\mathcal M_\lambda$, rewrite the subproblem in Step 1 of (\ref{Pre iter}) as 
\begin{eqnarray*}
    \begin{aligned}
{u^{k + 1}} = \arg \mathop {\min }\limits_{u \in U} \left\{{ 
\frac{1}{{2\lambda }}\| {W \odot ( {\mathcal P u - Y} )} \|^2  + \langle {\nabla u,{p^k} + \alpha {q^k}} \rangle } \right.&\\
\left. +\frac{1}{2}\| {u - {u^k}} \|_{{\mathcal M_\lambda }}^2  \right\}&.           \end{aligned}
\end{eqnarray*}
Regardless of the constraint, the derivative of the subproblem is calculated as
\begin{eqnarray*}
    \begin{aligned}
&\frac{1}{\lambda }{{\mathcal P}^T}\big(W \odot W \odot ( {\mathcal P u - Y} )\big) - {\rm{div}} ( {{p^k} + \alpha {q^k}} ) + {\mathcal M_\lambda}( {u - {u^k}} )\\
 =& \left( \frac{1}{\lambda }{{{\mathcal P}^T}( {W \odot W \odot \mathcal P }) + {\mathcal M_\lambda}} \right)u - \frac{1}{\lambda }{{\mathcal P}^T}\big( {W \odot W \odot Y} \big)\\
& - {\rm{div}} \big( {{p^k} + \alpha {q^k}} \big) -  {\mathcal M_\lambda}{u^k},
    \end{aligned}
\end{eqnarray*}
where $\rm{div}$ is the divergence operator. Let  $\mathcal M_\lambda$ satisfies
\begin{eqnarray}\label{gamma condition}
    \frac{1}{\lambda }{{{\mathcal P}^T}\left( {W \odot W \odot \mathcal P} \right)+ {\mathcal M_\lambda}} = \gamma \mathcal I
\end{eqnarray}
with $\gamma >0$ and $\mathcal I$ being identity operator. Immediately for $u$-subproblem, one just needs to solve the following constrained quadratic  problem 
\begin{eqnarray*}
    \begin{aligned}
{u^{k + 1}}  = \arg \mathop{\min }\limits_{u \in U} \left\{ 
\frac{1}{2}\left\|  u - \frac{1}{\gamma}\left( \frac{1}{\lambda }{{\mathcal P}^T}\left( {W \odot W \odot Y} \right) + {\rm{div}}\big( {{p^k} + \alpha {q^k}} \big)\right.\right.\right.&\\ 
\left. \left. \left.+ {\mathcal M_\lambda }{u^k} \right) \right\|^2\right\}&.
\end{aligned}
\end{eqnarray*}
The close-form solution of Step 1 is given as
\begin{eqnarray}
\label{Pre u}
    \begin{aligned}            
        {u^{k + 1}} = \frac{1}{\gamma }{\rm{Proj}}\Big({\frac{1}{\lambda }{{\mathcal P}^T}\left( {W \odot W \odot Y} \right) + {\rm{div}}\big( {{p^k} + \alpha {q^k}} \big) + {\mathcal M_\lambda }{u^k};U} \Big)
    \end{aligned}
\end{eqnarray}
with the projection operator defined as 
\begin{eqnarray*}
{\rm{Proj}}\left( {u;U} \right) := \min \left\{ {\max \left\{ {u,0} \right\},c} \right\}.    
\end{eqnarray*}

Second, consider the subproblem in Step 3 of (\ref{Pre iter}) w.r.t. the variable $q$ as
\begin{eqnarray*}
  {q^{k + 1}} = \arg \mathop {\min }\limits_{q \in Q} \left\{ {\alpha \langle {\nabla {{\bar u}^{k + 1}},q} \rangle  + \frac{1}{{2\tau }}{{\| {q - {q^k}} \|}^2}} \right\}.  
\end{eqnarray*}
 Consequently the close-form solution of Step 3 is given as
\begin{eqnarray}
\label{Pre q}
    \begin{aligned}
    {q^{k + 1}} = {\rm{Proj}}\Big( {{q^k} - {\tau\alpha}\nabla {\bar u^{k+1}}};Q \Big)
    \end{aligned}
\end{eqnarray}
with the projection operator defined as
\begin{eqnarray*}
[{\rm{Proj}}\left( {q;Q} \right)]_{i,j} := \frac{q_{i,j}}{{\max \left( {1,\sqrt {{{ {({q_x})_{i,j}^2} }} + {{ {({q_y})_{i,j}^2} }}}} \right)}}.    
\end{eqnarray*}

Finally, the subproblem w.r.t. the variable $p$ reads
\begin{eqnarray}
\label{Pre p}
    {p^{k + 1}} = {\rm{Proj}}\Big( {\frac{1}{{ {1 + \eta \beta }}}{\big( {{p^k} + \beta \nabla {{\bar u}^{k + 1}}} \big)}};S \Big)
\end{eqnarray}
with the projection operator defined as
\begin{eqnarray*}
[{\rm{Proj}}\left( {p;S} \right)]_{i,j} := \frac{p_{i,j}}{{\max \left( {1,\big| {({p_x})_{i,j}} \big|,\big| {({p_y})_{i,j}} \big|} \right)}}.    
\end{eqnarray*}
Based on the above calculations, the preconditioned  primal-dual hybrid gradient (Pre-PDHG) for Model (\ref{new saddle-point problem}) is summarized in Algorithm \ref{Pre Pdhg}.
\renewcommand{\thealgorithm}{1} 
    \begin{algorithm}
        \caption{ Pre-PDHG for Model (\ref{new saddle-point problem})} 
        \label{Pre Pdhg}
        \begin{algorithmic}[1] 
            \Require Set $q^0=0, \ p^{0}=0$, parameters $\lambda$, $\alpha$, $\eta$, $\gamma$, $\tau$ and $\beta$, maximum iteration number $N$, and tolerate accuracy $tol$.
            \Ensure sample in the spatial domain.
            \For{$k = 1 \to N$} 
                \If{Stopping criteria is not met}
                    \State Compute $u^{k+1}$ by (\ref{Pre u}).
                    \State Update $\bar u^{k+1}$ as Step 2 of (\ref{Pre iter}).
                    \State Compute $q^{k+1}$ by (\ref{Pre q}).
                    \State Compute $p^{k+1}$ by (\ref{Pre p}).        \EndIf
                \EndFor
            \State \Return $u$ 
        \end{algorithmic}
    \end{algorithm}

In Algorithm \ref{Pre Pdhg} the parameter $\gamma$ has to be selected to be big enough in order to satisfy the setting of the precondition, resulting in slow convergence. It is necessary to seek an alternative splitting method with fast convergence. 

\subsubsection{Fully-splitting  primal-dual hybrid gradient}
 A natural idea is to introduce the constraint $v=\mathcal Pu$ and add it to the objective function with the Lagrange multiplier $\Lambda \in \mathbb R^{m_1 \times m_2}$ as below
\begin{eqnarray}
    \begin{aligned}
       L\left( {u,v,q,p,\Lambda } \right) := &\frac{1}{{2\lambda }}\left\| {W \odot \left( {v - Y} \right)} \right\|^2 +{\mathbb{I} _U}\left( u \right)+ \left\langle {\Lambda ,v - \mathcal P u} \right\rangle\\
      & +\left\langle {\nabla u,p + \alpha q} \right\rangle   
       - {\mathbb{I} _S}\left( p \right) + {\mathbb{I} _Q}\left( q \right)  -\frac{\eta}{2}{\left\| p \right\|^2}.   
    \end{aligned}
\end{eqnarray}
Consequently, one needs to solve a saddle point problem as follows:
\begin{eqnarray*}
 \mathop {\min }\limits_{u,q,v} \mathop {\max }\limits_{p,\Lambda } L\left( {u,v,q,p,\Lambda } \right).   
\end{eqnarray*}

The scheme of the fully-splitting PDHG (FS-PDHG) in the $(k+1)_{\rm{th}}$ iteration can be described as:
\begin{eqnarray}
\label{Pdhg iteration}
   \begin{cases}
\texttt{Step 1:} \ & {\Lambda ^{k + 1}} = {\Lambda ^k} + {\rho}\big( {{ v^k} - \mathcal P{ u^k}} \big),\\
\texttt{Step 2:} \ & {u^{k + 1}} = \arg \mathop {\min }\limits_{u } L\big( {u,v^k,q^k,p^k,\Lambda^{k+1} } \big)  + \frac{{\| {u - {u^k}} \|^2}}{{2{\sigma _1}}},\\
\texttt{Step 3:}\ &{\bar u^{k + 1}} = {{2u^{k + 1}} - {u^k}} ,\\
\texttt{Step 4:} \ & {v^{k + 1}} =\arg \mathop {\min }\limits_v L\big( {\bar u^{k + 1},v,q^k,p^k,\Lambda^{k+1} } \big) + \frac{{\| {v - {v^k}} \|^2}}{{2{\sigma _2}}},  \\
\texttt{Step 5:} \  & {q^{k + 1}} = \arg \mathop {\min }\limits_{q} L\big( {\bar u^{k + 1},v^{k+1},q,p^k,\Lambda^{k+1} } \big) + \frac{{\| {q - {q^k}} \|^2}}{{2{\tau}}},\\
\texttt{Step 6:} \ & {p^{k + 1}} =\arg \mathop{\max}\limits_{p}  L\big( {\bar u^{k + 1},v^{k+1},q^{k+1},p,\Lambda^{k+1} } \big)- \frac{{\| {p - {p^k}} \|^2}}{{2\beta }}
   \end{cases} 
\end{eqnarray}
with  positive parameters $\rho$, $\sigma_1$, $\sigma_2$, $\tau$ and $\beta$. Here we remark that these subproblems w.r.t. the variables $u, v, q$, and  $p$ are strictly convex and have closed form solutions.

First, we consider the subproblem w.r.t. the variable $u$:
\begin{eqnarray*}
  {u^{k + 1}} = \arg \mathop {\min }\limits_{u \in U} \big\langle {\nabla u,{p^k} + \alpha {q^k}} \big\rangle  + \big\langle {{\Lambda ^{k + 1}}, - \mathcal P u} \big\rangle  + \frac{1}{{2{\sigma _1}}}{\big\| {u - {u^k}} \big\|^2}.  
\end{eqnarray*}
Similarly to (\ref{Pre u}),  one readily obtains the solution:
\begin{eqnarray}
\label{split u}
    {u^{k + 1}} = {\rm{Proj}}\Big( {{\sigma _1}{\rm{div}}\big( {{p^k} + \alpha {q^k}} \big) + {\sigma _1}{{\mathcal P}^T}{\Lambda ^{k + 1}} + {u^k};U} \Big).
\end{eqnarray}

Next, we consider the subproblem w.r.t. the variable $v$. The closed-form solution can be directly derived as follows:
\begin{eqnarray}
    \label{splitt v}
    {v^{k + 1}} = \left( {\frac{1}{{{\sigma _2}}}{v^k} - {\Lambda ^{k + 1}} + \frac{1}{\lambda }Y \odot W \odot W} \right) / \left( {\frac{1}{{{\sigma _2}}} + \frac{1}{\lambda }  W \odot W} \right).
\end{eqnarray}

Regarding the subproblems w.r.t. $q$ and $p$, their solutions are directly given below following (\ref{Pre q}) and (\ref{Pre p}):
\begin{eqnarray}
\label{split q p}
\begin{aligned}
{q^{k + 1}} &= {\rm{Proj}}\Big( {{q^k} - {\tau\alpha}\nabla {\bar u^{k+1}}};Q \Big),\\
{p^{k + 1}}& = {\rm{Proj}}\Big( {\frac{1}{{{1 + \eta \beta }}}{\big( {{p^k} + \beta \nabla {{\bar u}^{k + 1}}} \big)}};S \Big).
\end{aligned}
\end{eqnarray}

Finally, we summarize the FS-PDHG scheme in Algorithm \ref{Pdhg}.
\renewcommand{\thealgorithm}{2} 
    \begin{algorithm}
        \caption{FS-PDHG for Model (\ref{new saddle-point problem})} 
        \label{Pdhg}
        \begin{algorithmic}[1] 
            \Require Set $u^{0},v^{0},q^0,p^{0}=0$, parameters $\lambda$, $\alpha$, $\eta$, $\rho$, $\sigma_1$, $\sigma_2$, $\tau$ and $\beta$, maximum iteration number $N$, and tolerate accuracy $tol$.
            \Ensure sample in the spatial domain.
            \For{$k = 1 \to N$} 
                \If{Stopping criteria is not met}
                    \State Update the multiplier as Step 4 of (\ref{Pdhg iteration}).
                    \State Compute $u^{k+1}$ by (\ref{split u}). 
                    \State Update $\bar u^{k+1}$ as Step 3 of (\ref{Pdhg iteration}).
                    \State Compute $v^{k+1}$ by (\ref{splitt v}).
                    \State Compute $q^{k+1}$, $p^{k+1}$ in parallel by (\ref{split q p}).
                \EndIf
                \EndFor 
            \State \Return $u$ 
        \end{algorithmic}
    \end{algorithm}

\section{Convergence analysis}\label{convergence guarantee}
First, we give the general framework \cite{attouch2010proximal, lou2018fast, wang2019global} for the convergence analysis of non-convex optimization iterative algorithms. Assuming that $\Phi :{\mathbb{R}^d} \to \left( { - \infty , + \infty } \right]$ is proper and lower semicontinuous, we consider the following minimization problem:
\begin{eqnarray*}
\mathop {\min }\limits_z \Phi \left( z \right).    
\end{eqnarray*}
The iterative sequence $\left\{ {{z^k}} \right\}_{k = 0}^\infty$ is generated by any general algorithm $\mathscr{A}$ for solving the above problem. 
First recall the definition concerning subdifferential calculus.
\begin{definition}(Subdifferentials \cite{rockafellar2009variational}) Let $\Phi :\mathbb{R}{^d} \to \left( { - \infty , + \infty } \right]$ be a proper and lower semicontinuous function.\\
(1) For a given $x \in {\rm{dom}} \Phi$, the Fréchet subdifferential of $\Phi$ at $x$, written as $\hat \partial \Phi \left( x \right)$, is the set of all vectors $u \in \mathbb{R}^d$ which satisfy
\begin{eqnarray*}
 \mathop {\lim }\limits_{y \ne x} \mathop {\inf }\limits_{y \to x} \frac{{\Phi \left( y \right) - \Phi \left( x \right) - \left\langle {u,y - x} \right\rangle }}{{\left\| {y - x} \right\|}} \ge 0   
\end{eqnarray*}
with $ {\rm{dom}} \Phi : = \left\{ {x \in {\mathbb{R}^d}: \Phi \left( x \right) <  + \infty } \right\}$, and the notations $\langle \cdot \rangle$ and $ \| \cdot \|$ being defined as the standard inner product and $\ell_2$   norm in vector space $\mathbb R^d$ respectively. When $x \notin {\rm{dom}}\Phi $, we set $\hat \partial \Phi \left( x \right) = \emptyset$.\\
(2) The limiting-subdifferential, or simply the subdifferential, of $\Phi$ at $x\in \mathbb{R}^d$,
written as $\partial \Phi ( x )$, is defined through the following closure process
\begin{eqnarray*}
\begin{aligned}
  \partial \Phi \left( x \right): = \left\{ u \in {\mathbb{R}^d}:\exists {x^k} \to x, \Phi \big( {{x^k}} \big) \to \Phi \big( x \big) \,{\rm{and}} \right.&\\
 \left. {u^k} \in \hat \partial \Phi \big( {{x^k}} \big) \to u \,{\rm{as}} \, k \to \infty  \right\}&. 
\end{aligned}   
\end{eqnarray*}
\end{definition}

For a general non-convex optimization problem, the expectation  is to prove that the whole sequence generated by the algorithm $\mathscr{A}$ converges to a critical point of $\Phi$. Generally speaking, the subsequence convergence of  $\left\{{z^k}\right\}$ can be obtained easily applying the following analytical framework \cite{attouch2010proximal}. If $\Phi$ satisfies the KL property \cite{attouch2010proximal}, the iterative sequence can be further proved to be a Cauchy sequence such that the global convergence can be reached.
Specifically, the following three conditions are adopted to prove the convergence of the algorithm:
\begin{itemize}
\item Sufficient descent: $\mathscr{A}$ is essentially a descent algorithm, and each step has one lower bound estimation:
\begin{eqnarray}
\label{suff decent}
    \begin{aligned}
\Phi \big( {{z^k}} \big) - \Phi \big( {{z^{k + 1}}} \big) \ge {\rho _1}{\big\| {{z^{k + 1}} - {z^k}} \big\|^2}  \qquad k=0,1,\cdots      
    \end{aligned}
\end{eqnarray}
with a positive constant $\rho_1$.
\item Upper bound of the subgradient:
\begin{eqnarray}
\label{subgradient bound}
\big\| {{\omega ^{k + 1}}} \big\| \le {\rho _2}\big\| {{z^{k + 1}} - {z^k}} \big\|,\,{\omega ^{k + 1}} \in \partial \Phi \big( {{z^{k + 1}}} \big) \qquad   k=0,1,\cdots     
\end{eqnarray}
with a positive constant $\rho_2$.
\item Kurdyka– Lojasiewicz (KL) property \cite{attouch2010proximal}.
\end{itemize} 

Next, we prove the convergence of the proposed algorithms step by step in the following parts. The successive errors of the iterative sequence  are defined as:
\begin{eqnarray*}
E_u^{k + 1}: =  {u^{k + 1}} - {u^k},\,E_q^{k + 1}: = {q^{k + 1}} - {q^k},\,E_p^{k + 1}: = {p^{k + 1}} - {p^k},\\
E_v^{k + 1} : = {v^{k + 1}} - {v^k},\, E_\Lambda ^{k + 1}: = {\Lambda ^{k + 1}} - {\Lambda ^k}.
\end{eqnarray*}

\subsection{Convergence of Algorithm \ref{Pre Pdhg}}
The following analysis relies heavily on the auxiliary sequence defined  below
\begin{eqnarray}
\label{energy function}
{L_\lambda }( {u,q,p,\tilde u}): = {G }( u ) + g ^*( q ) - f_\eta ^*( p ) + \langle {\nabla u,p + \alpha q} \rangle  + \frac{1}{2}\| {u - \tilde u} \|_{{\mathcal M_\lambda }}^2.
\end{eqnarray}
The first-order optimality conditions for subproblems of (\ref{Pre iter}) are given as follows
\begin{eqnarray}
\label{u optimality condition}
{\rm{div}}\big( {{{p}^{k}} + \alpha {{q}^{k}}} \big) + {\mathcal M_\lambda}\big({{u^k} - {u^{k + 1}}}\big)\in \partial G\big( {{u^{k + 1}}} \big), \\ 
 - \alpha \nabla {{\bar u}^{k+1}} + \frac{{{q^{k}} - {q^{k+1 }}}}{\tau } \in \partial {g^*}\big( {{q^{k+1 }}} \big), \\
\nabla {{\bar u}^{k+1}} + \frac{{{p^{k}} - {p^{k+1 }}}}{\beta } \in \partial {f_\eta^*}\big( {{p^{k +1}}} \big).
\end{eqnarray}

\begin{condition}
\label{para condition}
The three positive parameters $\tau, \alpha, \eta$ satisfy
\begin{eqnarray*}
 2-{K\tau\alpha} >0 ,\,{\mathcal M_\lambda }- \left( {\frac{{2{K^2}}}{\eta } + K\alpha } \right) \mathcal I\succ 0,
 \end{eqnarray*}    
\end{condition}
with $\cdot\succ 0$ denoting strictly positive definite and $K$ being a bound on the norm of the linear
operator $\nabla$.

Using formula (\ref{u optimality condition}), $G$, $g^*$ convexity and the $k_{\rm{th}}$ update of $p$, the following Proposition \ref{propo} can be derived readily.
\begin{proposition}
\label{propo}
For all $k \ge 0$,  we have 
\begin{align}
\label{u sub}
{G }\big( {{u^k}} \big) - {G }\big( {{u^{k + 1}}} \big) \ge &\big\langle {\nabla \big( {{u^{k + 1}}-u^k} \big),{p^k} + \alpha {q^k}} \big\rangle   + \big\| {{E_u^{k + 1}} } \big\|_{{\mathcal M_\lambda }}^2, \\        
\label{q sub}
g ^*\big( {{q^k}} \big) - g ^*\big( {{q^{k + 1}}} \big) \ge &\big\langle {\alpha \nabla {{\bar u}^{k + 1}},{q^{k + 1}} - {q^k}} \big\rangle  + \frac{1}{{\tau }}{\big\| {{E_q^{k + 1}} } \big\|^2},  \\ 
\label{y sub}
f_\eta ^*\big( {{p^{k + 1}}} \big) - f_\eta ^*\big( {{p^k}} \big) \ge& \big\langle {\nabla {{\bar u}^k},{p^{k + 1}} - {p^k}} \big\rangle  + \big( {\frac{\eta }{2} + \frac{1}{{2\beta }}} \big){\big\| {{E_p^{k + 1}} } \big\|^2} \\
&+ \frac{{{{\big\| {{E_p^k}} \big\|}^2}}}{{2\beta }} - \frac{{{{\big\| {{p^{k + 1}} - {p^{k - 1}}} \big\|}^2}}}{{2\beta }} .  \nonumber 
\end{align}
\end{proposition}

\begin{lemma}
\label{saddle point problem and original model}
    If $X^* :=\left(u^*, q^*, p^*,\tilde u^*\right)$ is a critical point of the function $L_\lambda$, then $\left(u^*, q^*, p^*\right)$ is a critical point of $L_{PD}$. Meanwhile, when $\eta =0$, $u^*$ is a critical point of the original model (\ref{equ1}).
\end{lemma}
The corresponding proof can be found in Appendix \ref{proof model pre}.

Letting ${X^k}: = \big( {{u^k},{q^k},{p^k},{u^{k - 1}}} \big)$ be generated by Algorithm \ref{Pre Pdhg}, we first estimate a sufficient descent of $L_\lambda$ as shown in (\ref{suff decent}).

\begin{lemma}
\label{sufficient decrease condition}
For all $k > 0$,  under Condition \ref{para condition}, we obtain 
\begin{eqnarray*}
    \begin{aligned}
{L_\lambda }\big( {{X^k}} \big) - {L_\lambda }\big( {{X^{k + 1}}} \big) &\ge C{\big\| {{X^{k + 1}} - {X^k}} \big\|^2},
    \end{aligned}
\end{eqnarray*}
where $C:= \min \left\{ {C_1,~\frac{1}{\tau } - \frac{{K\alpha }}{{{2}}},~{\frac{\eta}{2} - K{\xi _1}}} \right\}$ with $C_1$ being the minimum eigenvalue of the operator ${\mathcal M_\lambda }- \big( {\frac{{2{K^2}}}{\eta } + K\alpha } \big) \mathcal I$.
\end{lemma}
The corresponding proof can be found in Appendix \ref{proof energy}. The following Lemma represents the boundedness of $\left\{{X^k}\right\}$ and $\left\{{L_\lambda(X^k)}\right\}$.
\begin{lemma}
    \label{boundedness}
The  sequence $\left\{ {{X^k}} \right\}$  generated by Algorithm \ref{Pre Pdhg} is bounded. Furthermore, the auxiliary sequence $\left\{ {{L_\lambda(X^k)}} \right\}$ is bounded and nonincreasing.     
\end{lemma}
The proof can be found in Appendix \ref{proof bound}. To prove that the sequence approaches a critical point, similar to (\ref{subgradient bound}), we need the following results.
\begin{lemma}
    \label{Subgradient upper bound}
For all $k \ge 0$, we get 
\begin{eqnarray}
    \begin{aligned}
{\rm{dist}}\big( {0,\partial {L_\lambda }\big( {{X^{k + 1}}} \big)} \big) \le {C_2}\big\| {{X^{k + 1}} - {X^k}} \big\|       \end{aligned}
\end{eqnarray}
with a positive constant $C_2:=2 \max \big\{ {K + K\alpha+{\big\| {{\mathcal M_\lambda }} \big\| },\frac{1}{\tau } + K\alpha ,\frac{1}{\beta } + K} \big\}$.
\end{lemma}
The corresponding proof can be found in Appendix \ref{proof subgradient}. 

Our objective is to prove the convergence of the entire sequence generated by Algorithm \ref{Pre Pdhg} to a critical point of (\ref{energy function}). To accomplish this, we discuss whether the optimization function satisfies the KL property.
It is trivial to see that  ${G}\left( u \right) + \left\langle {\nabla u,p + \alpha q} \right\rangle  + \frac{1}{2}\left\| {u - \tilde u} \right\|_{{\mathcal M_\lambda }}^2 -\frac{\eta}{2}{\left\| p \right\|^2}$ is the semi-algebraic function as well as the indicator functions $g ^*\left( q \right) - f^*\left(p \right)$, such that $L_\lambda$ is semi-algebraic \cite{xu2013block}, which  satisfies the KL property.
\begin{theorem}
		\label{A finite length property}
Under Condition \ref{para condition}, the iterative sequence $\left\{{X^k}\right\}$  generated by Algorithm \ref{Pre Pdhg} converges to a critical point $X^*$ of $L_\lambda$.		
\end{theorem}
The proof of Theorem \ref{A finite length property} follows standard techniques \cite{attouch2010proximal, lou2018fast, wang2019global} and is therefore omitted here. 

\subsection{Convergence of Algorithm \ref{Pdhg}.}
Similar to (\ref{energy function}), according to the iteration scheme (\ref{Pdhg iteration}), the auxiliary sequence is given as follows
\begin{equation}
\label{energy func}
    \begin{aligned}
{L_\sigma }\left( {\Lambda ,u,v,q,p,\tilde u,\tilde v} \right): =& F\left( v \right) + {h }\left( u \right) + {g^*}\left( q \right) - f_\eta ^*\left( p \right) + \left\langle {\nabla u,p + \alpha q} \right\rangle  \\
&+ \left\langle {\Lambda ,v - \mathcal P u} \right\rangle  + \frac{1}{{2\sigma_1 }}{\left\| {u - \tilde u} \right\|^2}+ \frac{1}{{2\sigma_2 }}{\left\| {v - \tilde v} \right\|^2}        \end{aligned}
\end{equation}
with $F\left( v \right) := \frac{1}{{2\lambda }}{\left\| {W \odot \left( {v - Y} \right)} \right\|^2}$ and $h\left( u \right) := {\mathbb{I}_U}\left( u \right) $.
\begin{condition}
   \label{paras condition}
The parameters $\lambda,\,  \tau,\,  \alpha, \, \sigma_1, \, \sigma_2, \, \eta \,\, {\rm{and}} \,\, \rho$ satisfy
\begin{equation*}
\begin{aligned}
 2-K\tau\alpha >0,\ \frac{{K{\sigma _1}}}{{1 - K\alpha }} < \frac{\eta }{{2K}},\ \frac{1}{{2{\sigma _2}}} - \frac{4}{\rho}\left( {\frac{1}{{{\lambda ^2}}}\left\| {W} \right\|_{\rm{max}} ^4 + \frac{1}{{\sigma _2^2}}} \right) > 0. 
\end{aligned} 
\end{equation*} 
\end{condition}

\begin{lemma}
  \label{split saddle point problem and original model}
   Let $Z^* :=\left(\Lambda^*,u^*,v^*,q^*,p^*,\tilde u^*,\tilde v^*\right)$ be a critical point of the functional $L_\sigma$. Then $\left(u^*,v^*, q^*, p^*,\Lambda^*\right)$ is a critical point of $L$. Moreover, when $\eta =0$, $u^*$ is a critical point of the original model (\ref{equ1}).
\end{lemma}
The corresponding proof can be found in Appendix \ref{proof model split}.\\
Denoting $Z^k :=(\Lambda^k,u^k,v^k,q^k,p^k, u^{k-1}, v^{k-1})$, the following relations hold as follows.
\begin{lemma}
\label{pdhg sufficient decrease condition}
For all $k > 0$, we have 
\begin{eqnarray*}
    \begin{aligned}
{L_\sigma }\big( {{Z^k}} \big) - {L_\sigma }\big( {{Z^{k + 1}}} \big) 
&\ge \tilde C{\big\| {{Z^{k + 1}} - {Z^k}} \big\|^2},
    \end{aligned}
\end{eqnarray*}
where 
\begin{eqnarray*}
\begin{aligned}
    \tilde C:={\rm{min}}\Big\{\frac{1}{2\sigma_1}-\frac{K}{2\xi_2},~{\frac{\eta}{2} - \frac{{K}}{{{\xi _2}}}},~\frac{1}{2\sigma_1}-\frac{K}{2\xi_2}-\frac{K\alpha}{2}, &\\
\frac{1}{{2{\sigma _2}}} - \frac{4}{\rho}\big( {\frac{1}{{{\lambda ^2}}}\left\| {W} \right\|_{\rm{max}} ^4 + \frac{1}{{\sigma _2^2}}} \big)\Big\}&.
\end{aligned} 
\end{eqnarray*}
\end{lemma}
The corresponding proof can be found in Appendix \ref{proof pdhg energy}. 

The introduction of the Lagrange multiplier results in $L_\sigma$ being non-coercive, which precludes easily proving the boundedness of $\left\{{Z^k}\right\}$. Thus, we hereby assume the boundedness of $\left\{{v^k}\right\}$.

\begin{assumption}
\label{seq boumdedness}
The sequence $\left\{{{v^k}}\right\}$ generated by 
Algorithm \ref{Pdhg} is bounded.
\end{assumption}

\begin{lemma}
    \label{pdhg boundedness}
Under Assumption \ref{seq boumdedness}, the  sequence $\left\{ {{Z^k}} \right\}$ generated by Algorithm \ref{Pdhg} is bounded. Furthermore, the auxiliary sequence $\left\{ {{L_\sigma(Z^k)}} \right\}$ is bounded and nonincreasing.     
\end{lemma}

The proof follows directly in a similar manner to that of Lemma \ref{boundedness}. Additionally, the following results are needed.
\begin{lemma}
    \label{pdhg Subgradient upper bound}
For all $k \ge 0$, we obtain 
\begin{eqnarray}
    \begin{aligned}
{\rm{dist}}\big( {0,\partial {L_\sigma }\big( {{Z^{k + 1}}} \big)} \big) \le {C_3}\big\| {{Z^{k + 1}} - {Z^k}} \big\| +\frac{1}{\rho}  \big\| {{Z^{k + 2}} - {Z^{k+1}}} \big\|    
\end{aligned}
\end{eqnarray}
with a positive constant $C_3 :=2 \max \left\{ {K + K\alpha+\frac{1}{\sigma_1},\frac{1}{\tau } + K\alpha ,\frac{1}{\beta } + K ,\frac{1}{\sigma_2} } \right\}$.
\end{lemma}

The proof follows in a similar manner to Lemma \ref{Subgradient upper bound} and thus we omit it here. Ultimately, the aforementioned results allow us to establish the final convergence theorem.
\begin{theorem}
		\label{FSP A finite length property}
Under Assumption \ref{seq boumdedness} and Condition \ref{paras condition}, the sequence $\left\{{Z^k}\right\}$ generated by Algorithm \ref{Pdhg} converges to a critical point $Z^*$ of $L_\sigma$, provided $1 - \frac{1}{{\rho {C_3}}} > 0$.
\end{theorem}

\section{Numerical experiments}\label{numerical experiment}
In this section, we validate the accuracy of our model through simulation data in the first four sections and proceed to test its performance employing real data in the fifth section. All numerical experiments are implemented in MATLAB R2021a on a laptop with 4-cores 3.1 GHz Intel Core and 16 GB RAM.

\subsection{ Experimental settings}

We utilize the fan-beam CT imaging system  for all simulations. The distance from the source to the detector is 949.075 mm, the  source to iso-center distance is 541 mm, and the strip width is 1.024 mm. A complete  360$^\circ$ rotation in a circular orbit comprises a total of 984 projected views, with each view containing 888 bins.

{\color{blue}All the phantoms are composed of soft tissue, bone, and two circular metal objects. The specific linear attenuation coefficients for these objects can be found in reference \cite{hubbell1995tables}. Each reconstructed image comprises $256 \times 256$ pixels. In the NCAT phantom image (depicted in Figure \ref{ncat truth}), the metal components are fabricated from titanium, whereas the iron is used for the other two phantoms (displayed in Figures \ref{truth ano} and \ref{new truth}).}

\begin{figure}[htbp]
\centering
\includegraphics[scale=0.5]{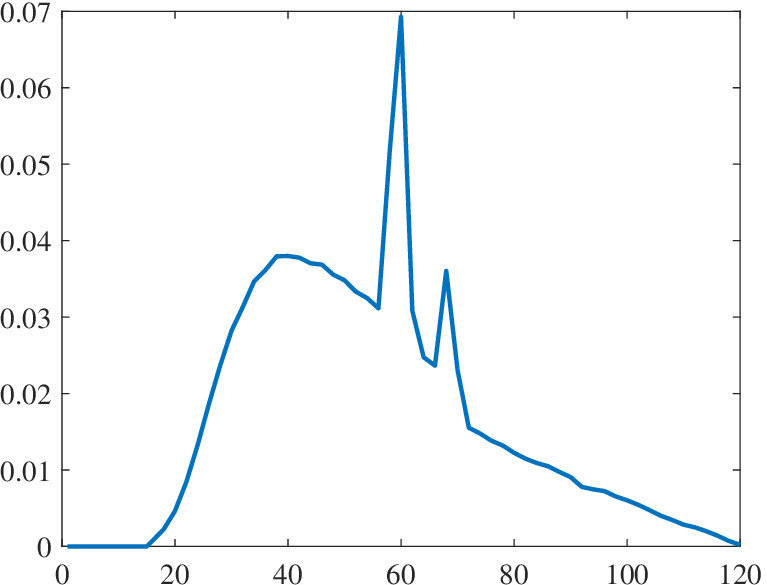}
\caption{ Energy spectrum $\mathcal S_0(E)$}
\label{spec}
\end{figure}

 The simulated $Y$ with the energy spectrum ${{\mathcal S_0}\left( E \right)}$ (as shown in Figure \ref{spec}) is obtained from the projected data  contaminated by Poisson noise in the following way \cite{beer1852bestimmung}
\begin{eqnarray}
\label{noise data}
   Y_{i,j} =  - \log \Big( {\max \left\{ {{\rm{Poissrnd}}\big( {{S_0}{e^{ - {Y_0}_{(i,j)}}}} \big)/{S_0},1/{S_0}} \right\}} \Big)\qquad\forall i,j, 
\end{eqnarray}
where Poissrnd$\big(  \cdot  \big)$ refers to the Poisson noise, $S_0$ is the number of incident photons, and the term $1/S_0$ is to avoid taking the logarithm of 0. We select $S_0 = 10^5$ for the NCAT phantom and $S_0 = 10^9$ for the other two in (\ref{noise data}).

The stopping criterion is defined as follows 
\begin{eqnarray*}
\frac{{\left\| {{u^{k + 1}} - {u^k}} \right\|}}{{\left\| {{u^{k + 1}}} \right\|}} \le tol.    
\end{eqnarray*}
For NCAT and head phantoms, we choose $tol = 9 \times 10^{-5}$, while for the skull, set it to $5 \times 10^{-5}$. For the weighted function $W$, the default values are set to $t=0.94$ and $\varepsilon=10^{-16}$.

\subsection{Tests on various weight functions }
\begin{figure}[htbp]
     \centering
     \begin{subfigure}{0.32\linewidth}
		\centering		\includegraphics[width=0.9\linewidth]{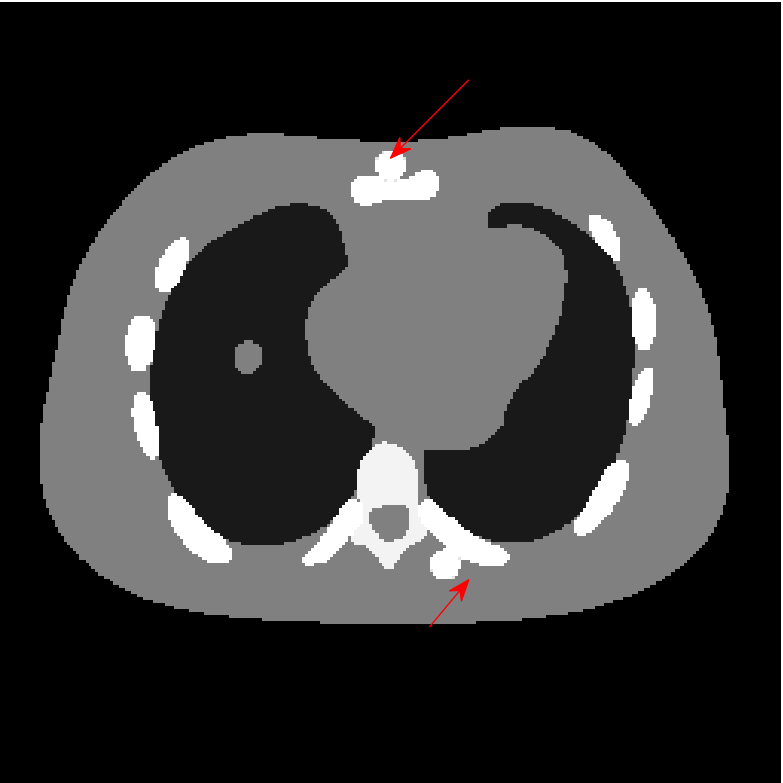}
		\caption{  }
		\label{metaltruth}
	\end{subfigure}
	\centering
	\begin{subfigure}{0.32\linewidth}
		\centering		\includegraphics[width=0.9\linewidth]{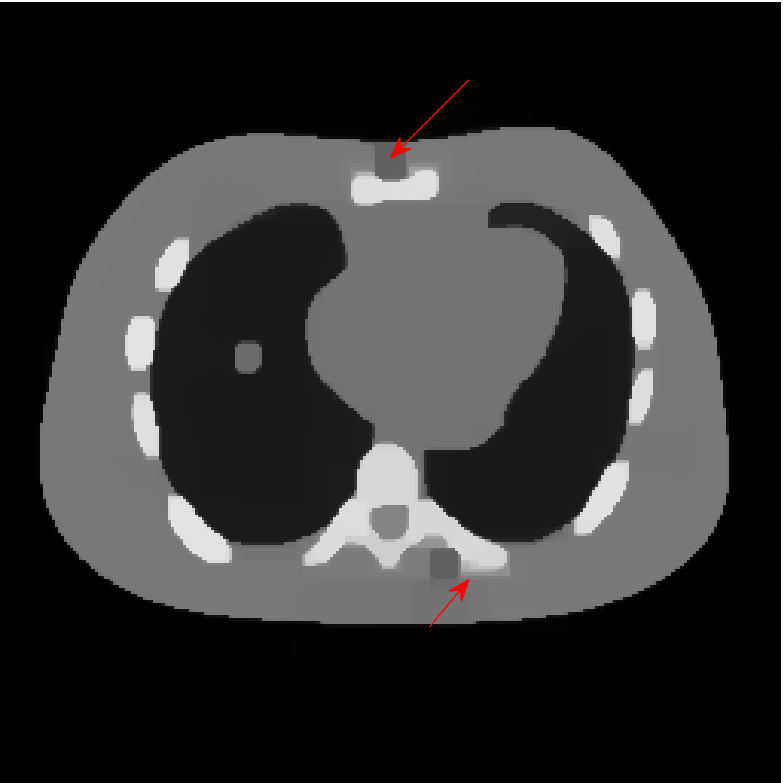}
		\caption{  }
		\label{metalmask1}
	\end{subfigure}
	\centering
	\begin{subfigure}{0.32\linewidth}
		\centering
		\includegraphics[width=0.9\linewidth]{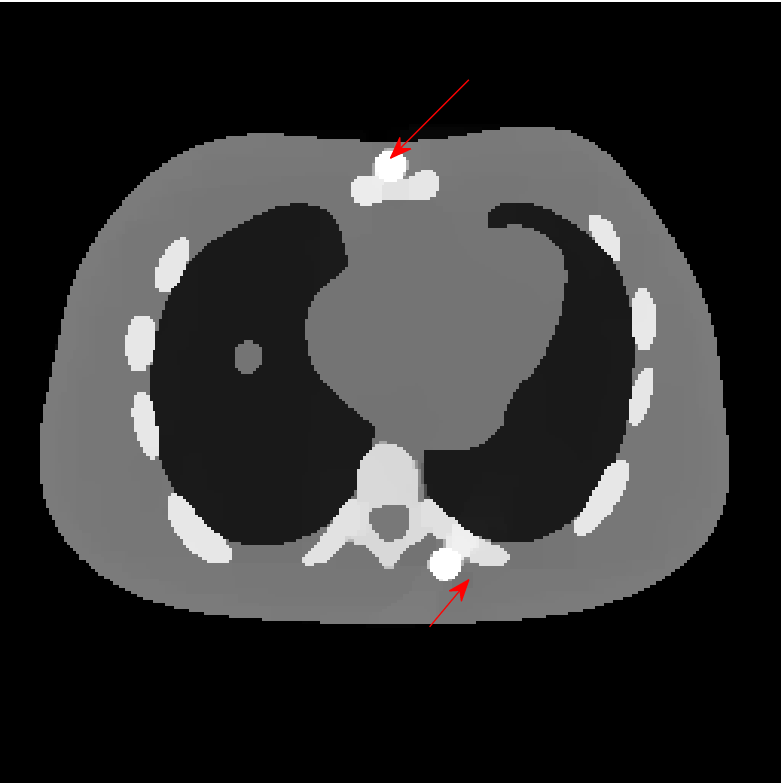}
		\caption{$   $ }
		\label{Yweight1}
	\end{subfigure}

     \centering
	\begin{subfigure}{0.32\linewidth}
		\centering		\includegraphics[width=0.9\linewidth]{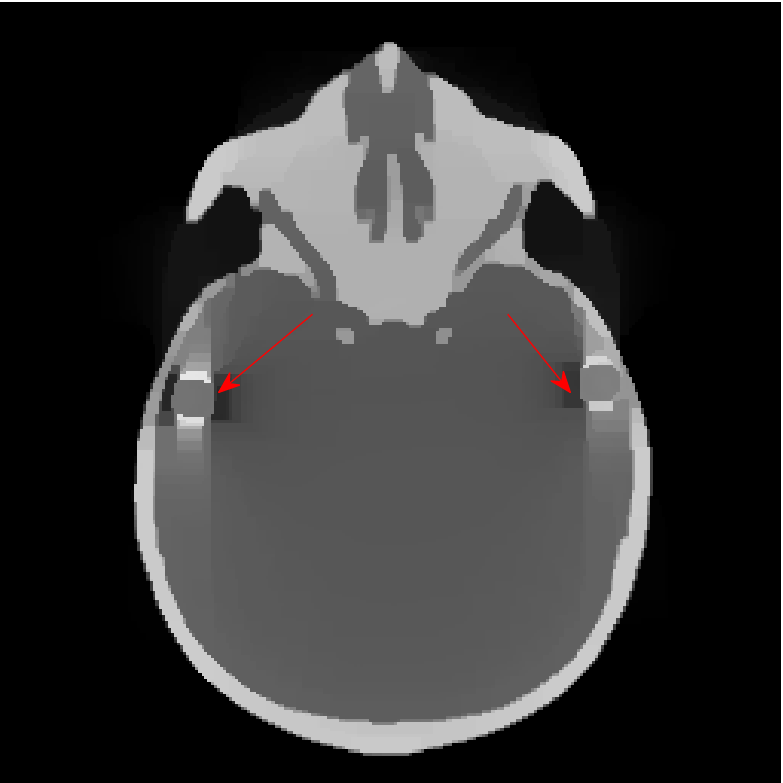}
		\caption{  }
		\label{metalmask}
	\end{subfigure}
	\centering
	\begin{subfigure}{0.32\linewidth}
		\centering
		\includegraphics[width=0.9\linewidth]{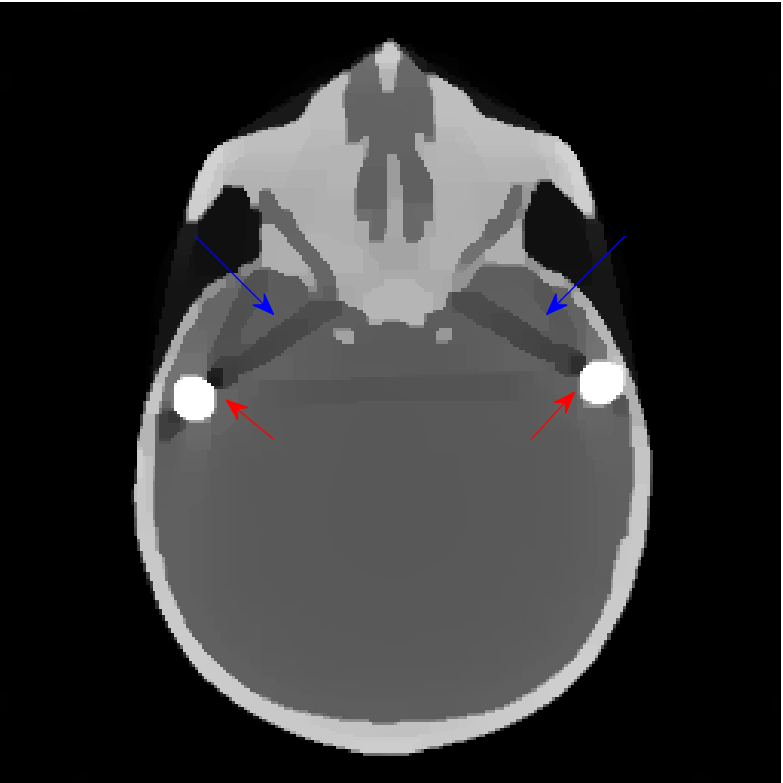}
		\caption{$   $ }
		\label{Yweight}
	\end{subfigure}
 	\centering
	\begin{subfigure}{0.32\linewidth}
		\centering
		\includegraphics[width=0.9\linewidth]{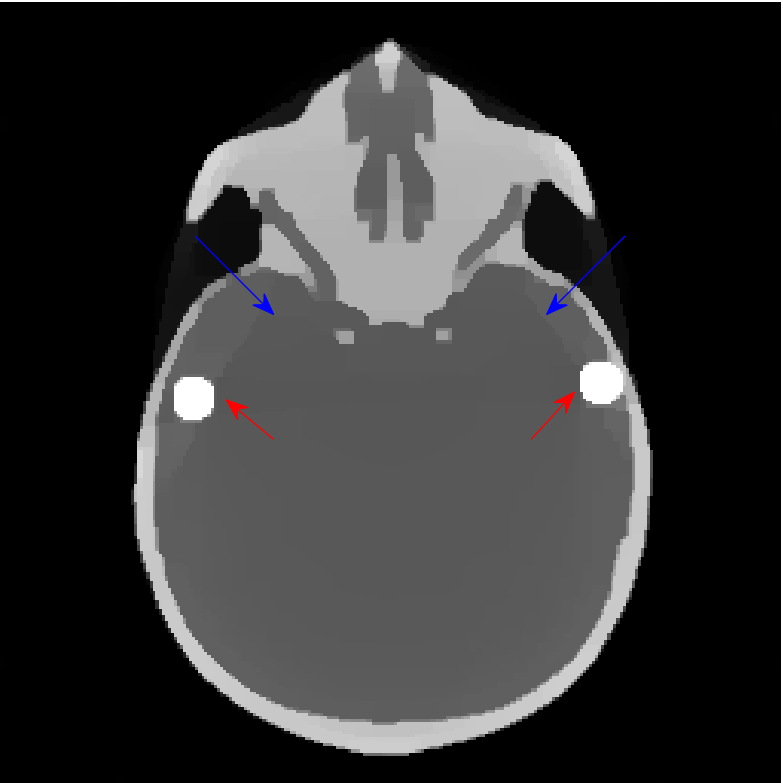}
		\caption{$   $ }
		\label{Yweight2}
	\end{subfigure}
	\caption{Comparison of the different weights. The first row is the reference image and the result of NCAT phantom in models (\ref{mask weight model}) and (\ref{non multiplier }), respectively (4600
HU window, 1300 HU level).  The second row is the reconstruction result of the head phantom in model (\ref{mask weight model}) and the selection of $O_m$ or $\Omega_t$ as the mask region in (\ref{non multiplier }) from left to right (4200 HU window, 1100 HU level). }	\label{weight curve}
\end{figure}

We use the following two phantoms in Figure \ref{weight curve} to demonstrate the efficiency of the proposed weight function $W$. For a fair comparison, we change the regularization of the model (\ref{mask model}) to the AITV:
\begin{eqnarray}
\label{mask weight model}
    {u^*} = \arg \mathop {\min }\limits_u \frac{1}{{2\lambda }}{\left\| {{B_\Omega^c} \odot \left( {\mathcal P u - Y} \right)} \right\|^2} + {\left\| {\nabla u} \right\|_1} - \alpha {\left\| {\nabla u} \right\|_{2,1}}.
\end{eqnarray}

\begin{figure}[htbp]
     \centering
	\begin{subfigure}{0.4\linewidth}
		\centering		\includegraphics[width=0.9\linewidth]{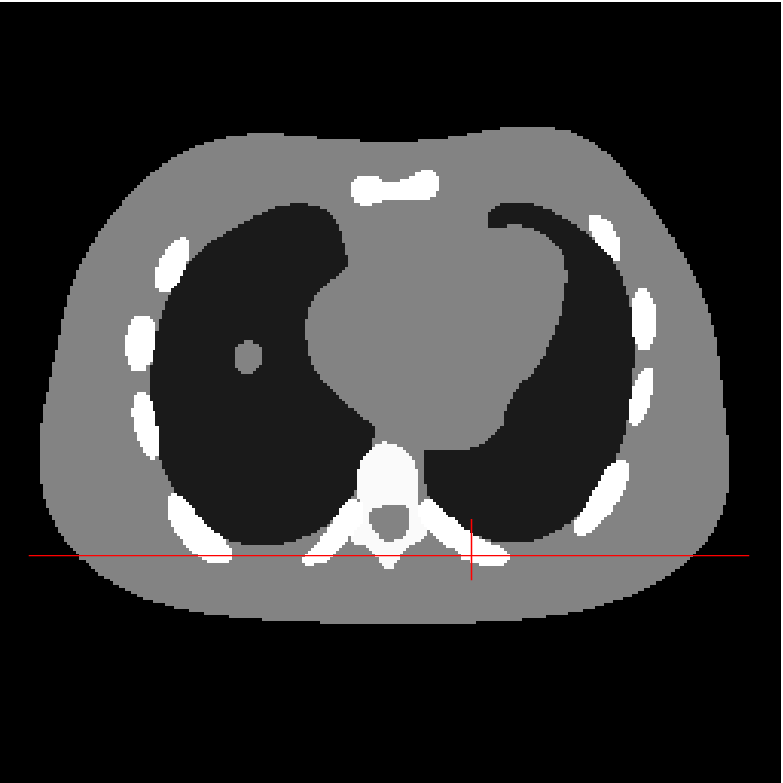}
		\caption{  }
		\label{hbbb}
	\end{subfigure}
	\centering
	\begin{subfigure}{0.4\linewidth}
		\centering
		\includegraphics[width=0.9\linewidth]{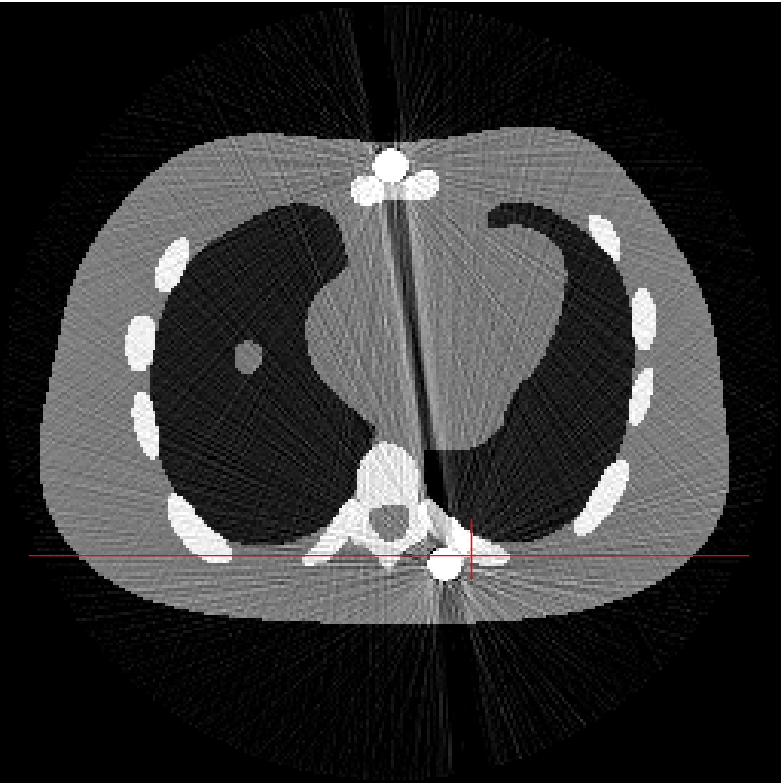}
		\caption{$   $ }
		\label{ma_CTred}
	\end{subfigure}
 
 	\centering
	\begin{subfigure}{0.4\linewidth}
		\centering
		\includegraphics[width=0.9\linewidth]{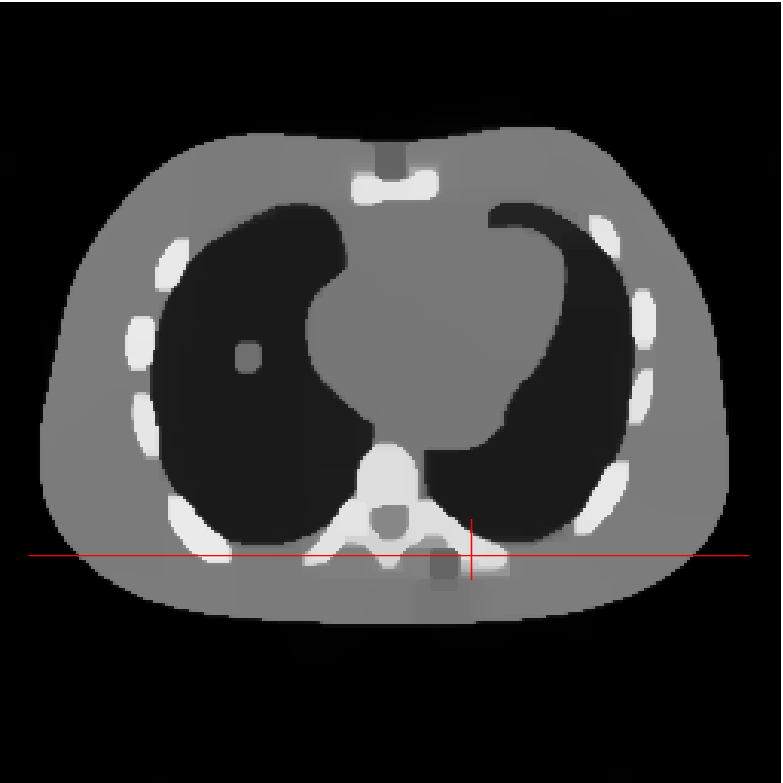}
		\caption{  }
		\label{fssss}
	\end{subfigure} 
 	\centering
	\begin{subfigure}{0.4\linewidth}
		\centering		\includegraphics[width=0.9\linewidth]{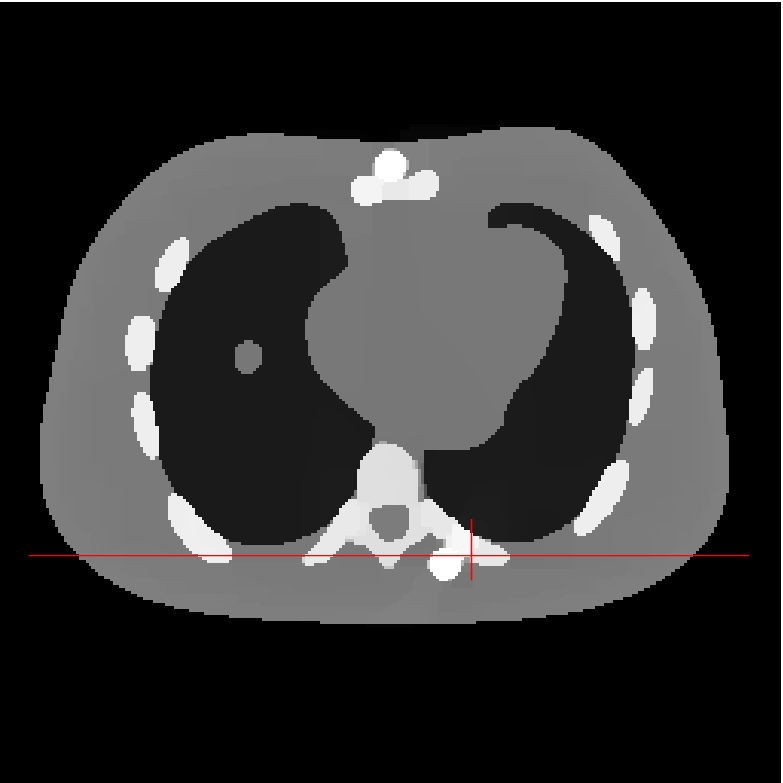}
		\caption{   }
		\label{tvline}
	\end{subfigure} 

 \centering
	\begin{subfigure}{0.46\linewidth}
		\centering		\includegraphics[width=0.9\linewidth]{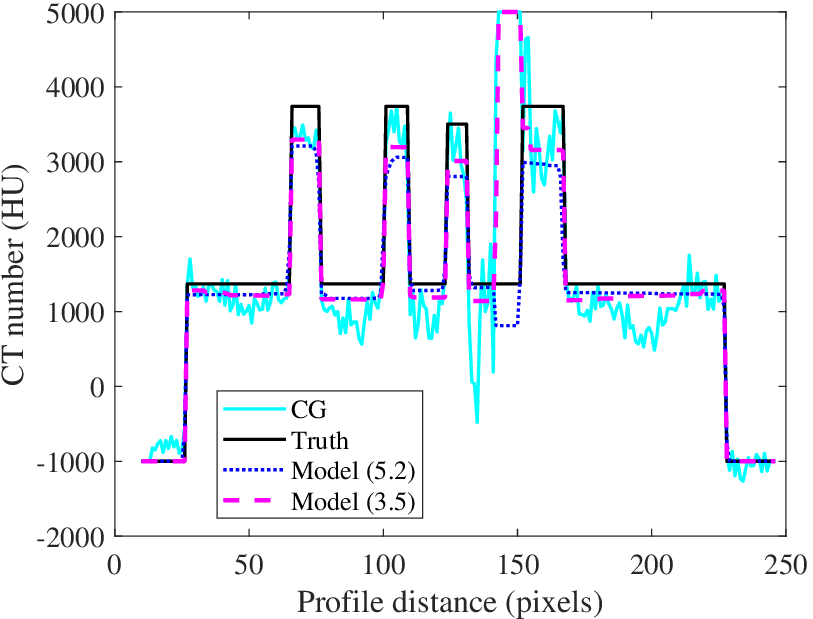}
		\caption{  }
		\label{hor}
	\end{subfigure}
	\centering
	\begin{subfigure}{0.46\linewidth}
		\centering
		\includegraphics[width=0.9\linewidth]{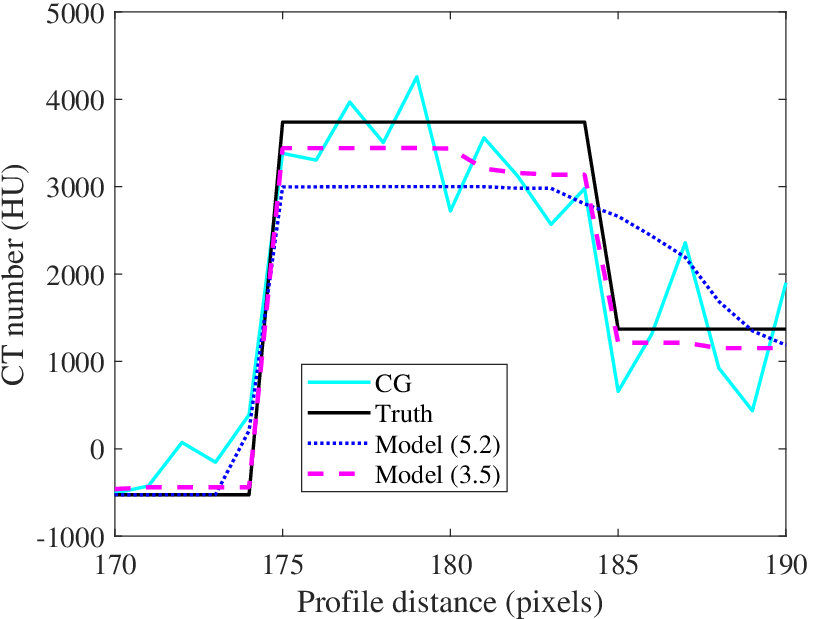}
		\caption{ }
		\label{ver}
	\end{subfigure}
	\caption{The third row represents the variation curve of Hounsfield Units (HU) corresponding to the positions of the red lines in the reconstructed images (\ref{hbbb})-(\ref{tvline}). From left to right, it shows the horizontal profile and the vertical profile. The remaining two rows correspond to the real metal-free image, the CG reconstruction image, the reconstruction result of model (\ref{mask weight model}), and the reconstruction result of model (\ref{non multiplier }) (4600 HU window, 1300 HU level).}	\label{TV-FS HU curve}
\end{figure}

 The results show that model (\ref{mask weight model}) is relatively worse than the proposed model in balancing the reduction of metal artifacts and the sharpness of the image.  In Figures \ref{metalmask1}-\ref{Yweight1}, the image reconstructed by the model (\ref{mask weight model}) has severe diffusion near the metal (see the red arrow below), and almost the metal information is lost because there is no prior related to the metal. On the contrary, the proposed model performs significantly better near the bone. {\color{blue} To provide a more detailed representation, we display the CT values of the two cutlines adjacent to the metal in Figure \ref{TV-FS HU curve}. The results shown in Figure \ref{hor} clearly demonstrate that the proposed model performs better in terms of correcting metal artifacts compared to model (\ref{mask weight model}), as it exhibits a closer proximity to the reference value. In particular, Figure \ref{ver} further highlights the superiority of our method, as it not only demonstrates a closer resemblance to the reference, but also exhibits favorable contrast (near pixel 185 in Figure \ref{ver}).}
 
\begin{figure}[htbp]
     \centering
	\begin{subfigure}{0.32\linewidth}
		\centering		\includegraphics[width=0.9\linewidth]{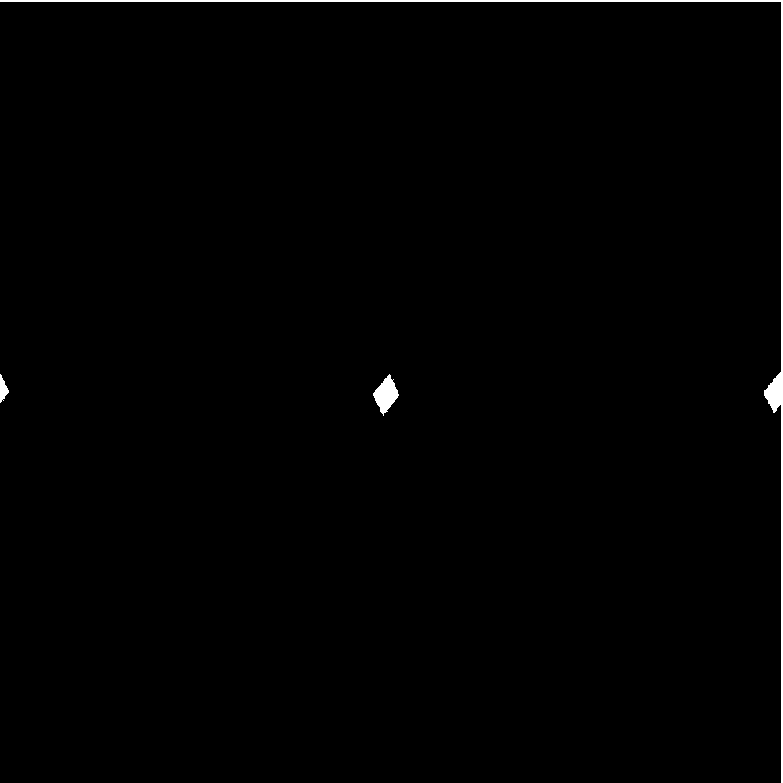}
		\caption{ $O_m$ }
		\label{Om}
	\end{subfigure}
	\centering
	\begin{subfigure}{0.32\linewidth}
		\centering
		\includegraphics[width=0.9\linewidth]{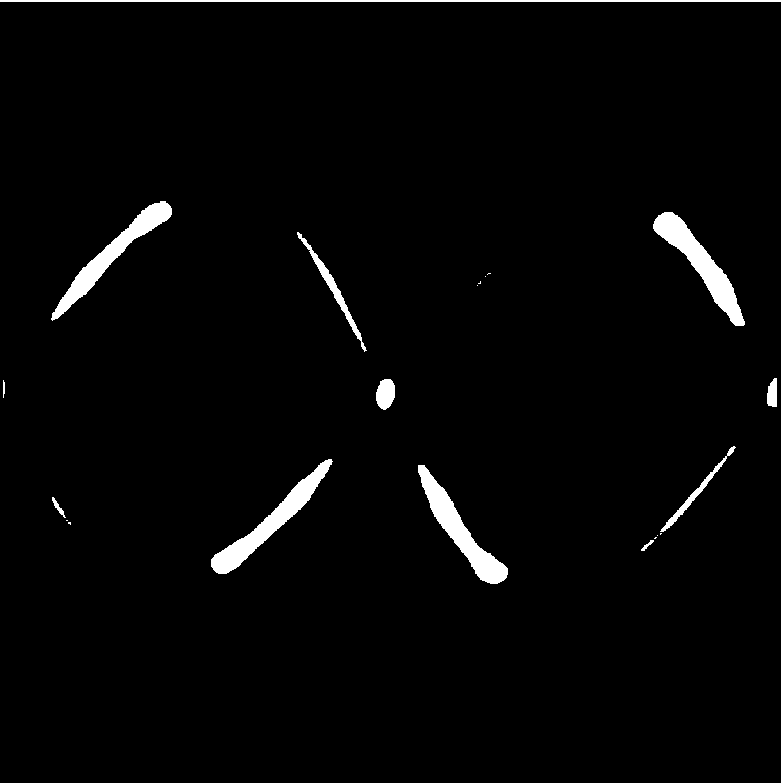}
		\caption{$ O_t  $ }
		\label{Ot}   
	\end{subfigure}
     \centering
	\begin{subfigure}{0.32\linewidth}
		\centering		\includegraphics[width=0.9\linewidth]{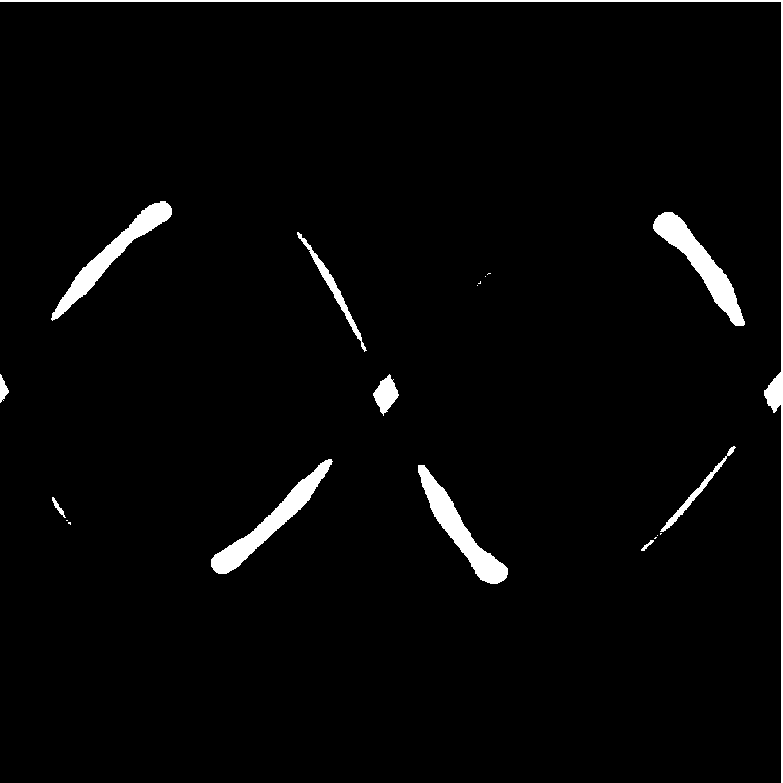}
		\caption{ $\Omega_t$ }
		\label{Omega}   
	\end{subfigure}
	\caption{The mask region of the proposed weight for the head phantom. From left to right are $O_m$, $O_t$ and $\Omega_t$, respectively.}	\label{head omega mask}
\end{figure}

 In the above test, $O_m$ and $O_t$ are almost the same. We will show an example  with the large bone, leading to quite different regions (as shown in Figure \ref{head omega mask}). The results are put to 
 Figures \ref{metalmask}-\ref{Yweight2}, by considering the head phantom. 
Similarly, from Figure \ref{metalmask},  one can see clear artifacts and boundary diffusion simply with the weight function in \eqref{mask weight model} (see red arrows). In Figure \ref{Yweight} obtained by replacing $\Omega_t$ with $O_m$ in the weight function $W$ of the proposed model, one readily sees that  it is difficult to remove artifacts between metal and bone (see the blue arrow), while using the proposed weight function based on $O_t$ and $O_m$,  the metal artifacts are greatly reduced as shown in Figure \ref{Yweight2}.


\subsection{Convergence validation and parameter impact}
\begin{figure}[htbp]
	\centering
	\begin{subfigure}{0.45\linewidth}
		\centering		\includegraphics[width=0.9\linewidth]{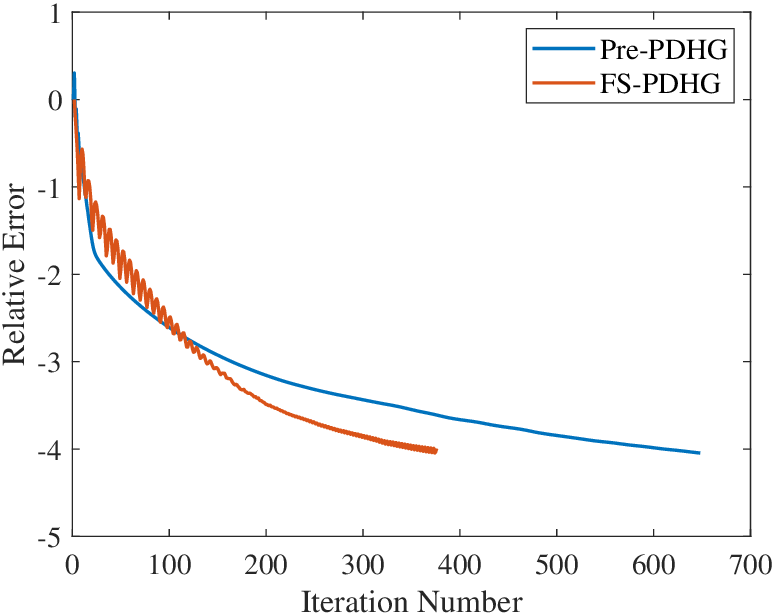}
		\caption{ }
		\label{u err}   
	\end{subfigure}
	\centering
	\begin{subfigure}{0.45\linewidth}
		\centering		\includegraphics[width=0.9\linewidth]{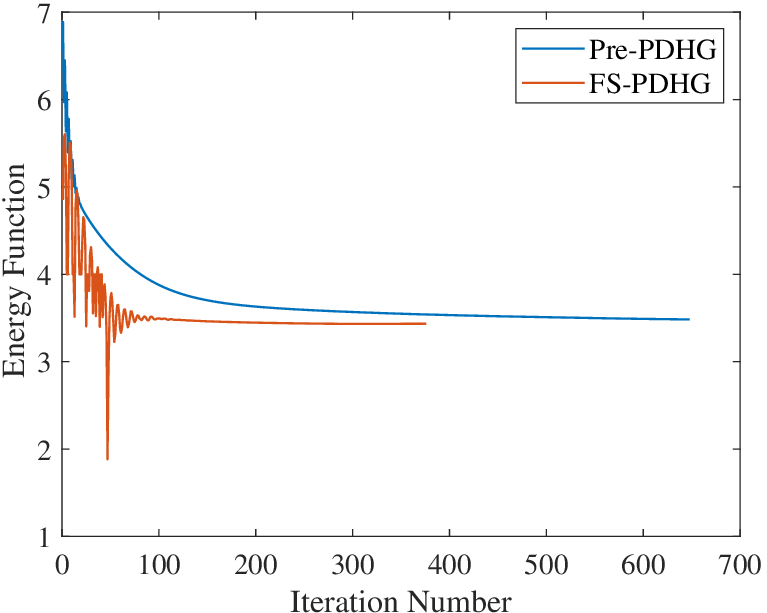}
		\caption{ }
		\label{cost }   
	\end{subfigure}
	\caption{Convergence of Algorithm \ref{Pre Pdhg} and Algorithm \ref{Pdhg} on the NCAT. From left to right are the curves for the relative error and energy function of the proposed algorithms, respectively.}	\label{conv curve}
\end{figure}

\noindent{\bf{Convergence validation}} \quad In order to numerically illustrate the convergence behavior of Algorithm \ref{Pre Pdhg} and Algorithm \ref{Pdhg}, the convergence curves of ${\log _{10}}\Big( {\frac{{\big\| {{u^{k + 1}} - {u^k}} \big\|}}{{\big\| {{u^{k + 1}}} \big\|}}} \Big)$ and the energy function is shown in Figure \ref{conv curve}. One readily sees that though local oscillations appear for the proposed FS-PDHG algorithm, both the relative error and the energy keep decreasing during iterations, which is consistent with the  convergence analysis in section \ref{convergence guarantee}.  Moreover, it is evident that the FS-PDHG algorithm exhibits faster convergence compared to the Pre-PDHG algorithm since the former reaches the desired tolerance with much fewer iterations than the latter.  

\begin{figure}[htbp]
     \centering
	\begin{subfigure}{0.45\linewidth}
		\centering		\includegraphics[width=0.9\linewidth]{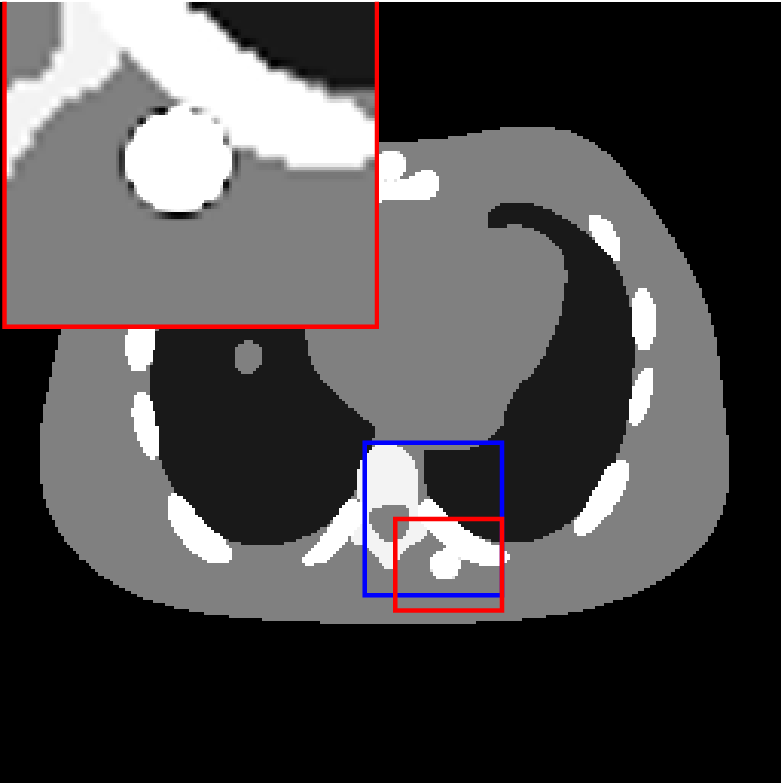}
		\caption{original }
		\label{truth}   
	\end{subfigure}
	\centering
	\begin{subfigure}{0.45\linewidth}
		\centering
		\includegraphics[width=0.9\linewidth]{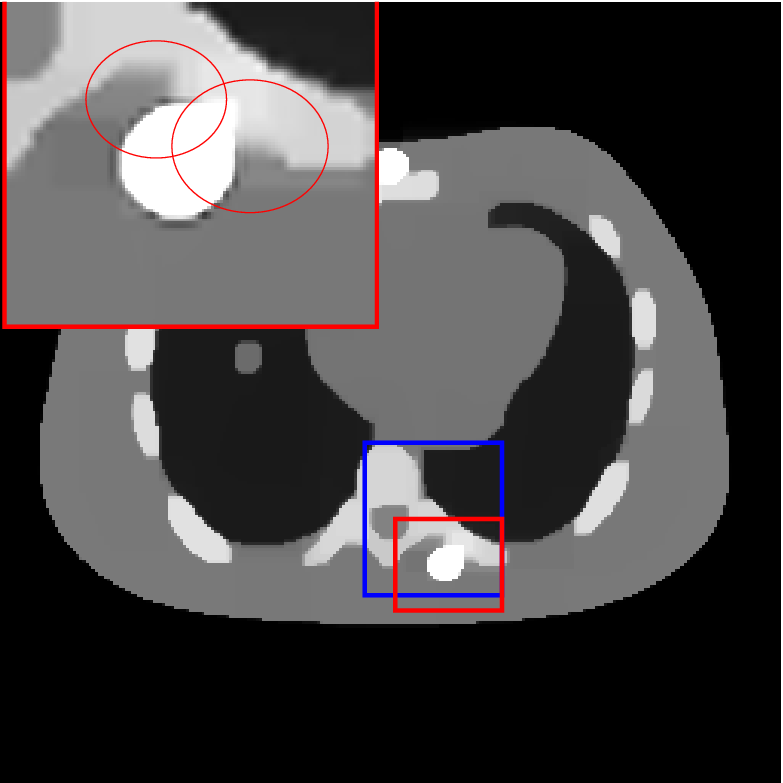}
		\caption{$\alpha =0$ }
		\label{alpha=0}   
	\end{subfigure}
 
	\centering
	\begin{subfigure}{0.45\linewidth}
		\centering
		\includegraphics[width=0.9\linewidth]{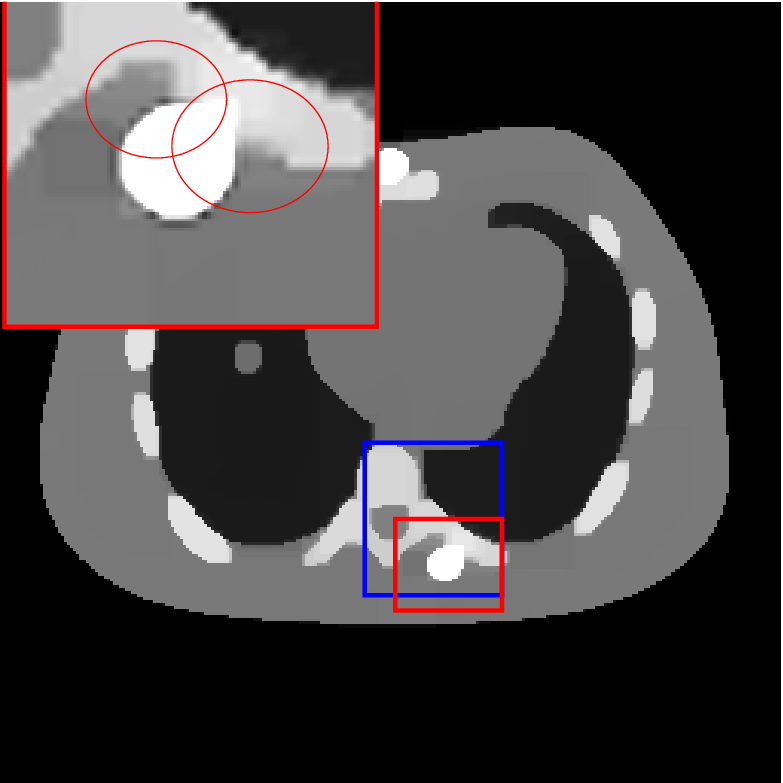}
		\caption{ $\alpha =0.25$}
		\label{alpha=0.25 }   
	\end{subfigure}
    \centering
	\begin{subfigure}{0.45\linewidth}
		\centering
		\includegraphics[width=0.9\linewidth]{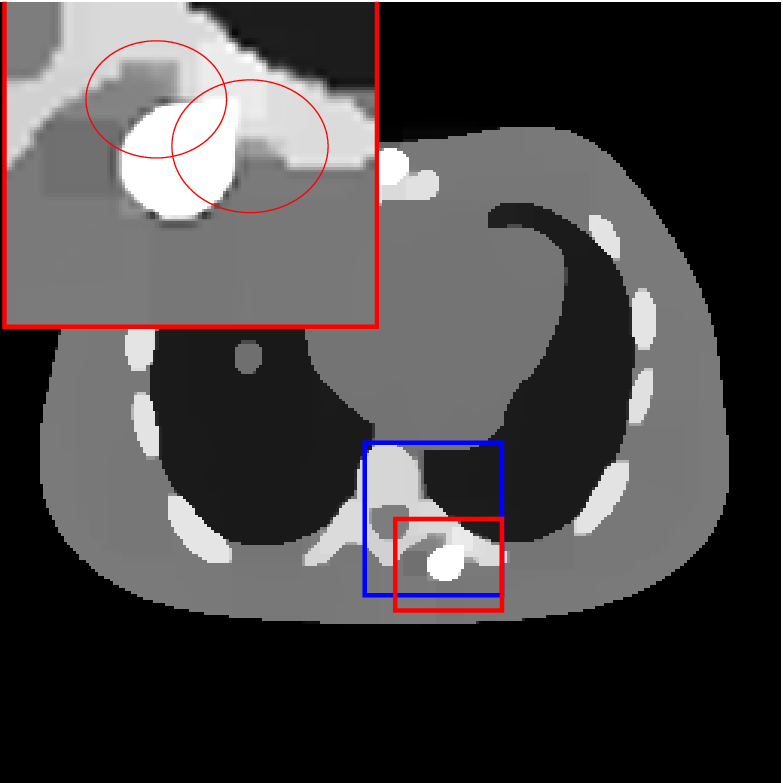}
		\caption{ $\alpha =0.5$}
		\label{alpha=0.5 }   
	\end{subfigure}
 
	\begin{subfigure}{0.45\linewidth}
		\centering
		\includegraphics[width=0.9\linewidth]{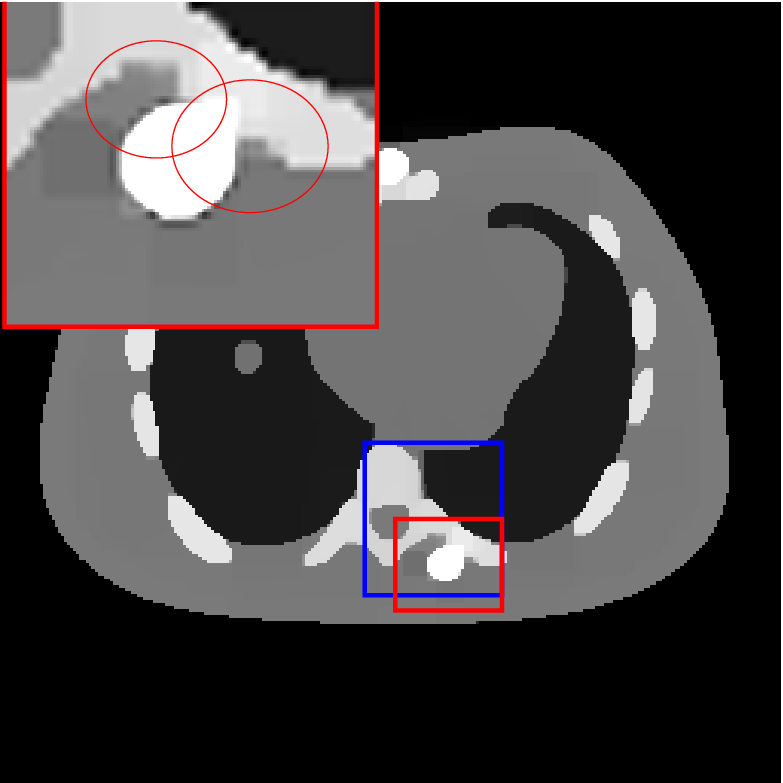}
		\caption{$\alpha =0.75$}
		\label{alpha=0.75}   
	\end{subfigure} 
    \centering
	\begin{subfigure}{0.45\linewidth}
		\centering
		\includegraphics[width=0.9\linewidth]{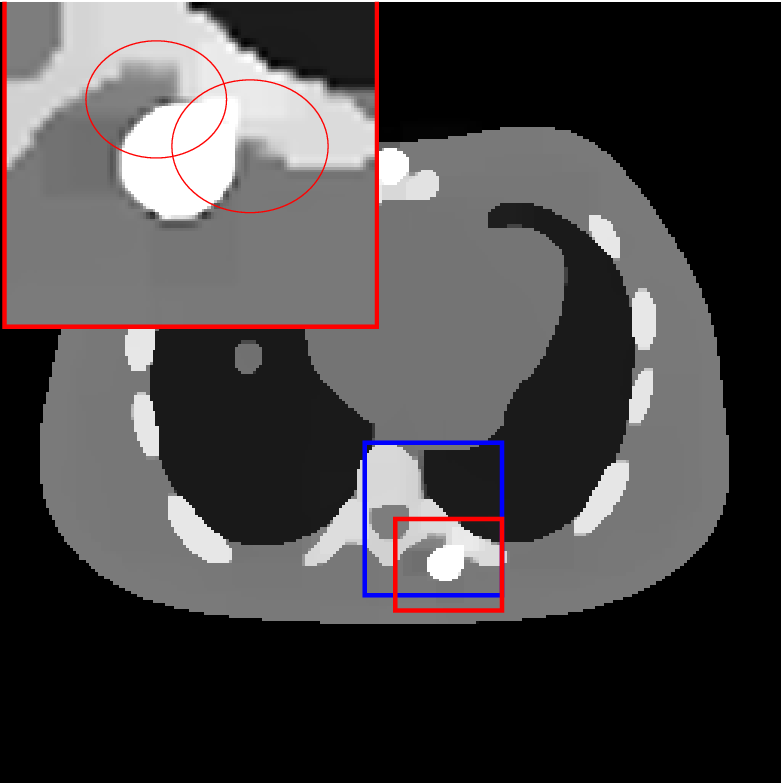}
		\caption{$\alpha =1$}
		\label{alpha=1}   
	\end{subfigure}
	\caption{{\color{blue}Comparison results with different parameters $\alpha$ (4600 HU window, 1300 HU level).}}	\label{alpha curve image}
\end{figure}

\begin{figure}[htbp]
	\centering
	\begin{subfigure}{0.45\linewidth}
		\centering		\includegraphics[width=0.9\linewidth]{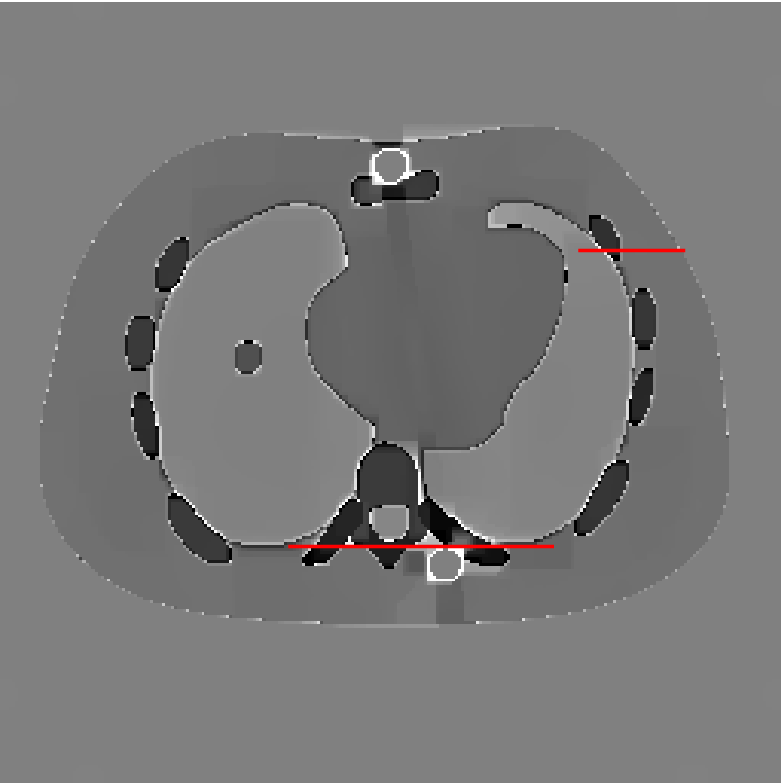}
		\caption{ $ $ }
		\label{diff 0}   
	\end{subfigure}
	\centering
	\begin{subfigure}{0.45\linewidth}
		\centering		\includegraphics[width=0.9\linewidth]{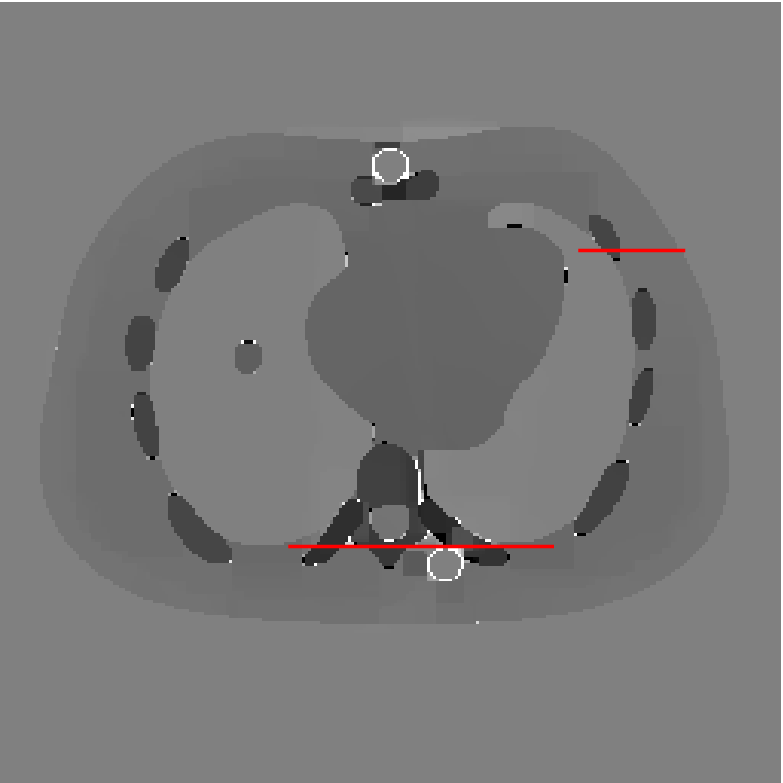}
		\caption{ $ $ }
		\label{diff 75 }   
	\end{subfigure}

 \centering
	\begin{subfigure}{0.45\linewidth}
		\centering		\includegraphics[width=0.9\linewidth]{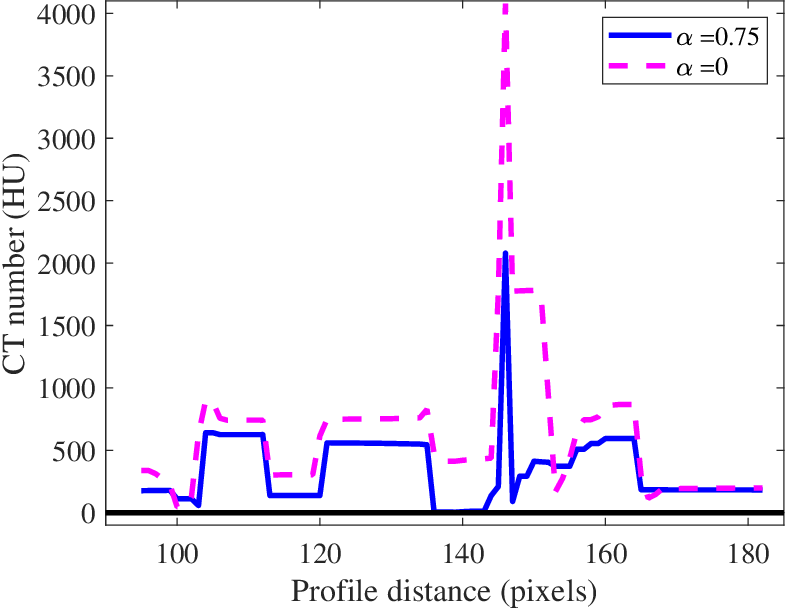}
		\caption{ $ $ }
		\label{diff long}   
	\end{subfigure}
	\centering
	\begin{subfigure}{0.45\linewidth}
		\centering		\includegraphics[width=0.9\linewidth]{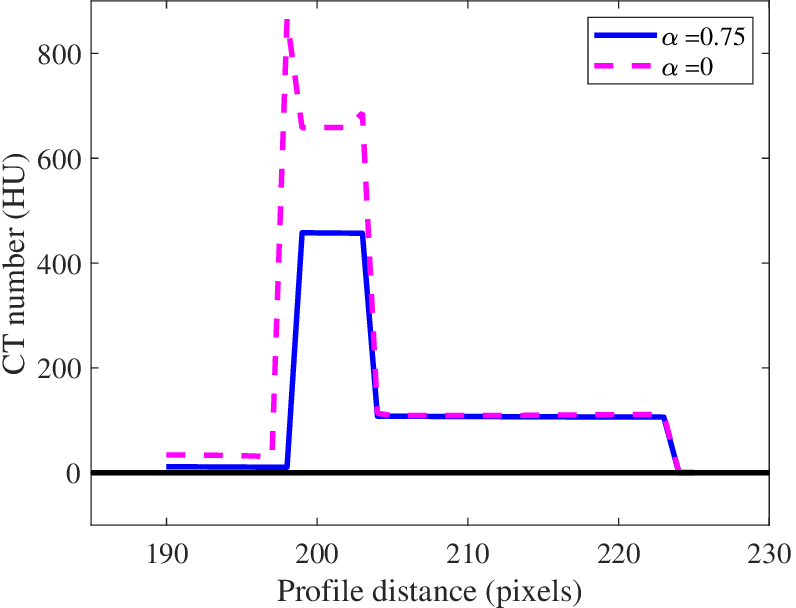}
		\caption{ $  $ }
		\label{diff  }   
	\end{subfigure}
	\caption{The first row represents the difference between the reconstruction results with parameters $\alpha=0$ and $\alpha=0.75$ and the reference image, and the second row displays the CT value curve highlighting the absolute value of the difference along the marked red line from left to right (2000 HU window, 0 HU
level).}	\label{diff images}
\end{figure}

\begin{figure}[htbp]
     \centering
	\begin{subfigure}{0.45\linewidth}
		\centering		\includegraphics[width=0.9\linewidth]{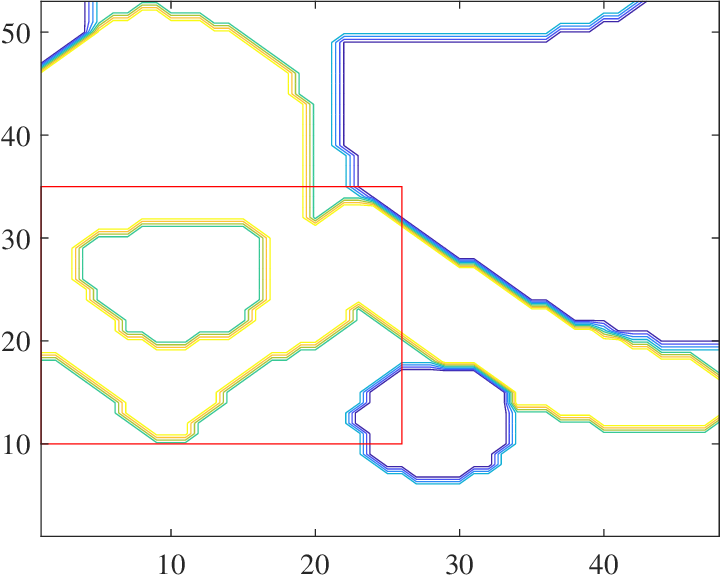}
		\caption{original }
		\label{truth con}   
	\end{subfigure}
	\centering
	\begin{subfigure}{0.45\linewidth}
		\centering
		\includegraphics[width=0.9\linewidth]{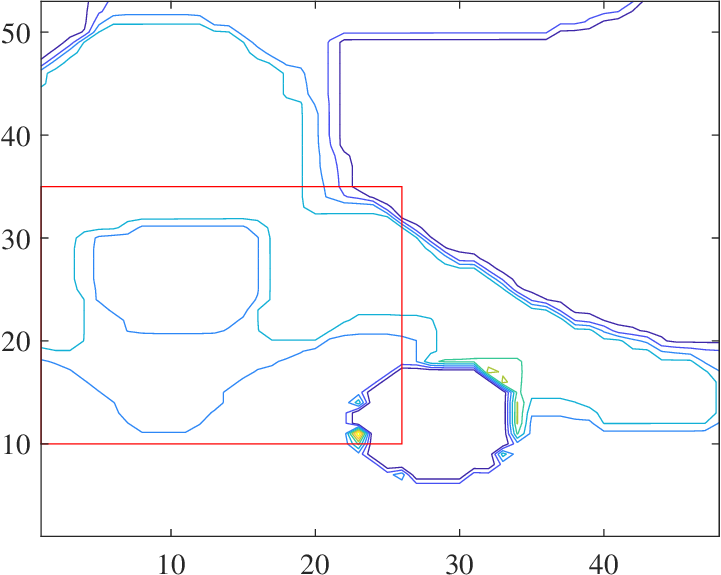}
		\caption{$\alpha =0$ }
		\label{alpha=0 con}   
	\end{subfigure}
 
	\centering
	\begin{subfigure}{0.45\linewidth}
		\centering
		\includegraphics[width=0.9\linewidth]{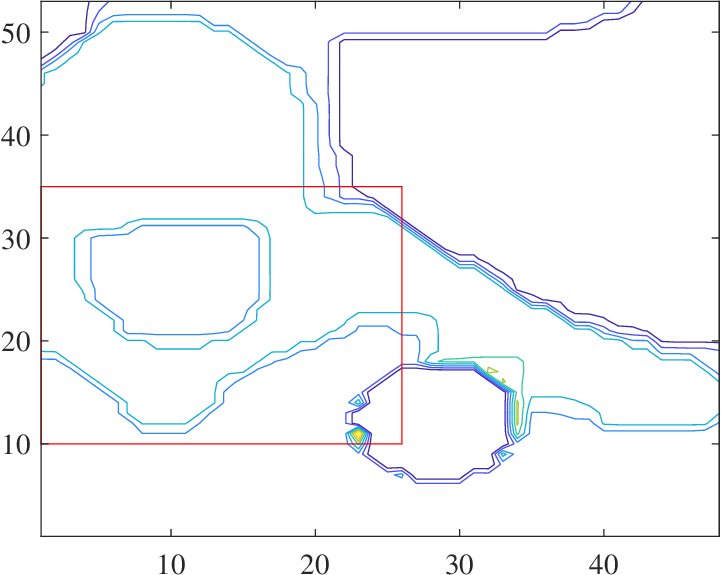}
		\caption{ $\alpha =0.25$}
		\label{alpha=0.25 con}   
	\end{subfigure}
    \centering
	\begin{subfigure}{0.45\linewidth}
		\centering
		\includegraphics[width=0.9\linewidth]{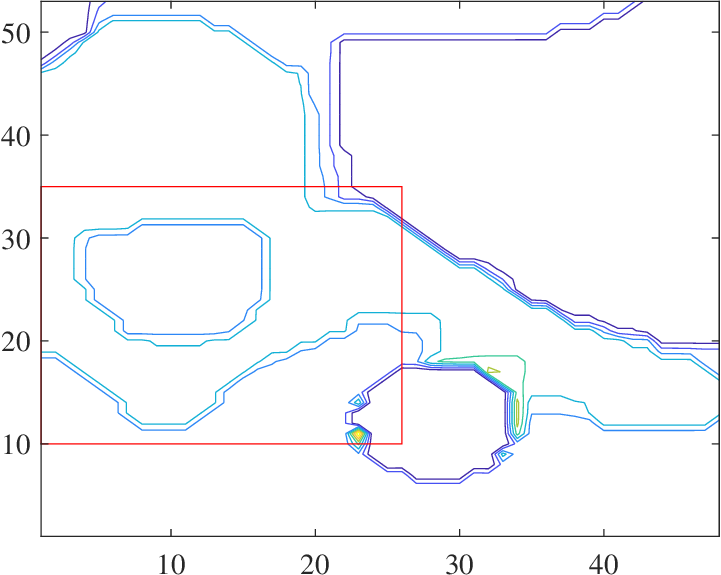}
		\caption{ $\alpha =0.5$}
		\label{alpha=0.5 con}   
	\end{subfigure}
	\begin{subfigure}{0.45\linewidth}
		\centering
		\includegraphics[width=0.9\linewidth]{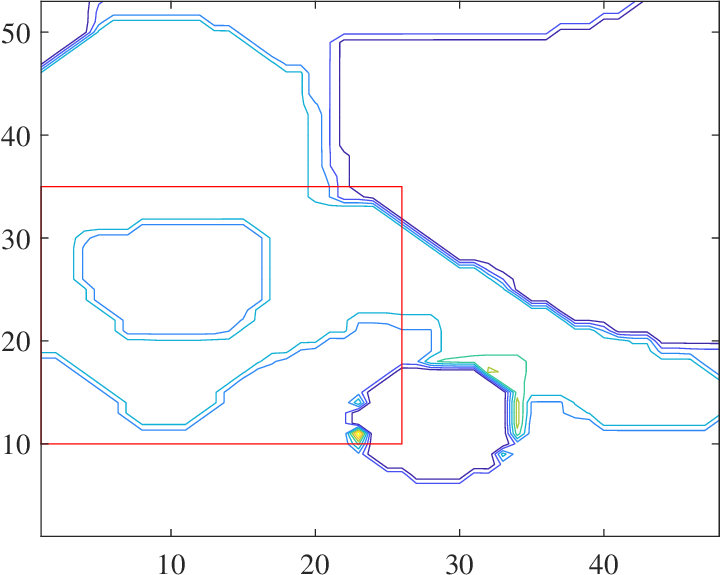}
		\caption{$\alpha =0.75$}
		\label{alpha=0.75 con}   
	\end{subfigure} 
    \centering
	\begin{subfigure}{0.45\linewidth}
		\centering		\includegraphics[width=0.9\linewidth]{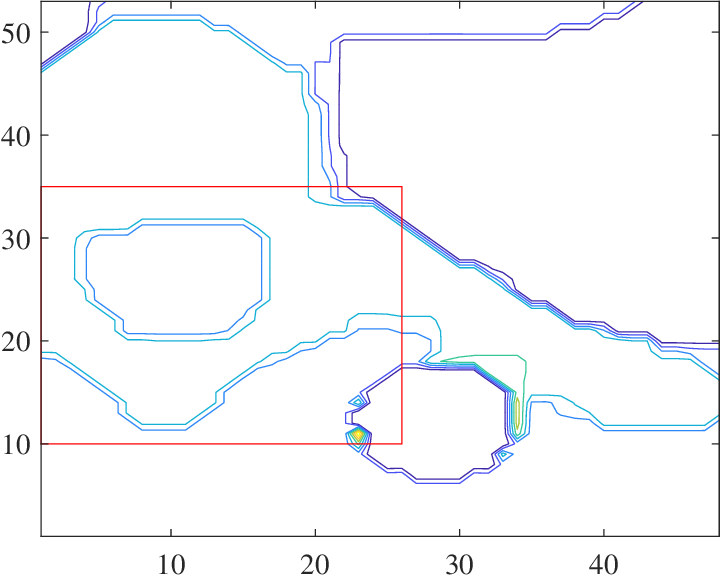}
		\caption{$\alpha =1$}
		\label{alpha=1 con}   
	\end{subfigure}
	\caption{Comparison of contours with different parameters $\alpha$.}	\label{alpha curve}
\end{figure}

\noindent{\bf Parameter impact}\quad We first test the performance of the proposed algorithms with different nonconvex regularization parameter $\alpha$ based on the NCAT. Table \ref{ssim alpha} records the reconstruction errors  (normalized $L^2$ norm between the recovery result and the truth), SSIM and PSNR values with respect to the different values of $\alpha$, where one readily sees that the proposed algorithms with proper parameter $\alpha=0.75$ produce best results in terms of these three metrics, and the nonconvex regularization obtain the PSNR improvement of 1.5dB compared with the convex model by setting $\alpha=0.$  We also put the visual results in  Figure \ref{alpha curve image}.  One can see obvious diffusion near the metal when $\alpha=0$. However, if $\alpha>0$, particularly between the metal and bone, there is an improved boundary recovery, as depicted by the red circle in Figure \ref{alpha=0}. The blur effect is noticeably reduced at $\alpha \ne 0$. Among the different cases, $\alpha =0.25$ exhibits relatively weaker improvement, as shown in the circled area below Figure \ref{alpha=0.25 }.
Comparing the results of $\alpha=0.5$ with $\alpha=0.75$, the recovered edges of the latter are sharper. The reconstructed images of $\alpha=0.75$ and $\alpha=1$ are very close, but the former yields better quantitative values according to Table \ref{ssim alpha}. In order to further highlight the details, we display the contours (see Figure \ref{alpha curve}) framed in blue within the region where the metals are removed. The results within the red rectangle further demonstrate the effectiveness of such nonconvex regularization. Based on the experimental tests in Table \ref{ssim alpha} and Figure \ref{alpha curve image}, we set $\alpha =0.75$ as the default, unless otherwise specified. {\color{blue}
Figures \ref{diff 0}-\ref{diff 75 } demonstrate the disparity between the reference image and the outcomes obtained with $\alpha=0$ and $\alpha=0.75$, respectively. It is apparent that the boundary of the reconstruction with $\alpha=0$ is not accurately depicted, whereas the results obtained with $\alpha=0.75$ significantly enhance the boundary. By examining the absolute HU value curve (shown in Figures \ref{diff long}-\ref{diff  }) along the red line in Figures \ref{diff 0}-\ref{diff 75 }, the smaller the differences obtained with $\alpha=0.75$, the better the reconstruction results. In contrast, the results obtained with $\alpha=0$ exhibit noticeable jumps near the boundary, further indicating inaccurate boundary. }
 
\begin{table}
\caption{\label{ssim alpha}Reconstruction error, SSIM and PSNR index of the NCAT phantom reconstructed by different $\alpha$, i.e. $\alpha=0, 0.25, 0.5, 0.75$ and $1$ respectively.}
\begin{indented}
\item[]\begin{tabular}{@{}llll}
\br
$\alpha$       & Reconstruction error         & SSIM         & PSNR      \\
\mr
0               & 0.1314                & 0.9766       & 28.1629 \\
        0.25            & 0.1103                & 0.9850       & 28.8709 \\
        0.5             & 0.1035                & 0.9883       & 29.4360 \\
        0.75            & {\bf{0.1003}}                & {\bf{0.9892}}       & {\bf{29.7163}} \\
        1               & 0.1006                & 0.9891       & 29.6967 \\
\br
\end{tabular}
\end{indented}
\end{table}

Next, we test the impact of other parameters for the proposed algorithms. To show how the algorithms vary w.r.t the parameters, we change the parameters in the following manner $c_0 \times 2^{k0+r0 \cdot l}$, where  $c_0$, $k_0$ and $r_0$ are constants dependent on the parameters of interest,  and $l$ represents the change frequency. Unless otherwise specified, the default values for $k_0$ and $r_0$ are both set to 1.

\begin{figure}[htbp]
\begin{subfigure}{0.3\linewidth}
		\centering		\includegraphics[width=0.9\linewidth]{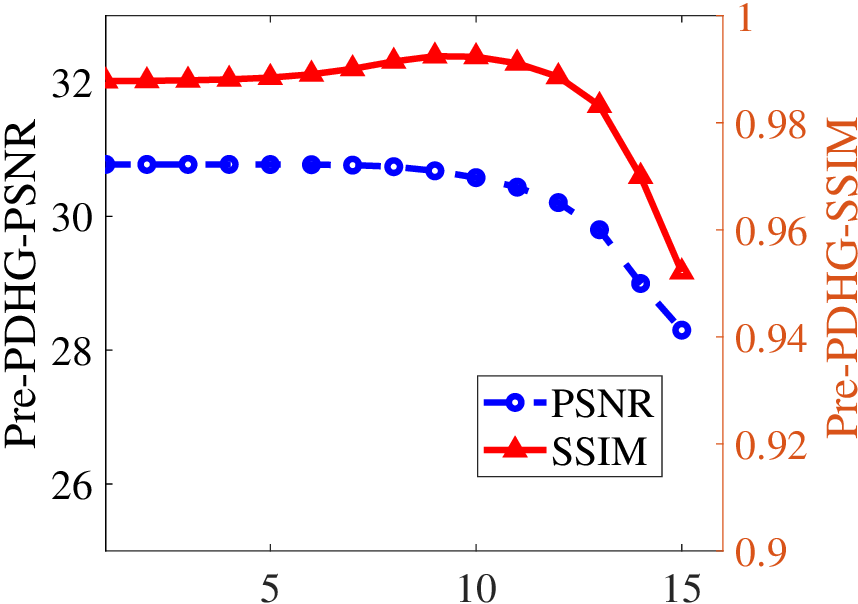}
		\caption{ $\eta$ }
		\label{etaPre}   
	\end{subfigure}
     \centering
	\begin{subfigure}{0.3\linewidth}
		\centering		\includegraphics[width=0.9\linewidth]{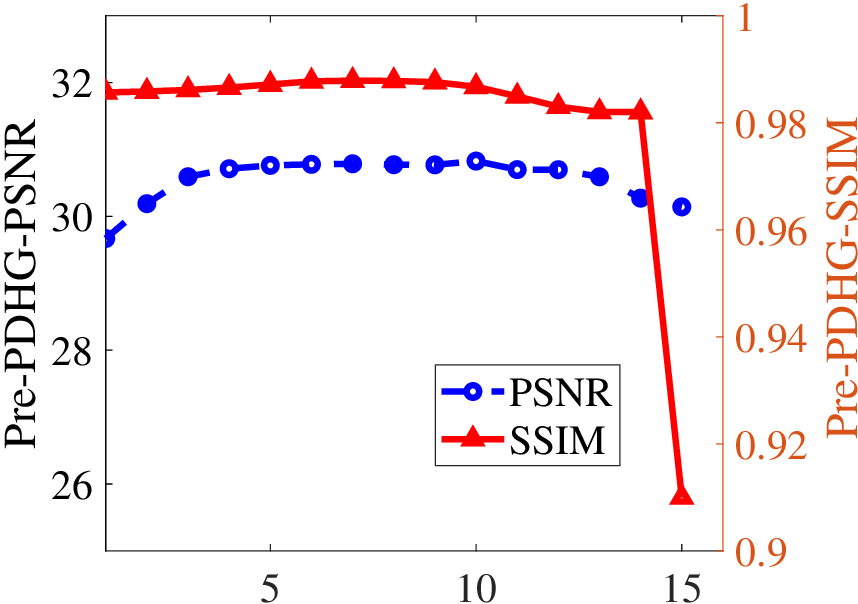}
		\caption{ $\tau$}
		\label{tauPre}   
	\end{subfigure}
	\centering
	\begin{subfigure}{0.3\linewidth}
		\centering
		\includegraphics[width=0.9\linewidth]{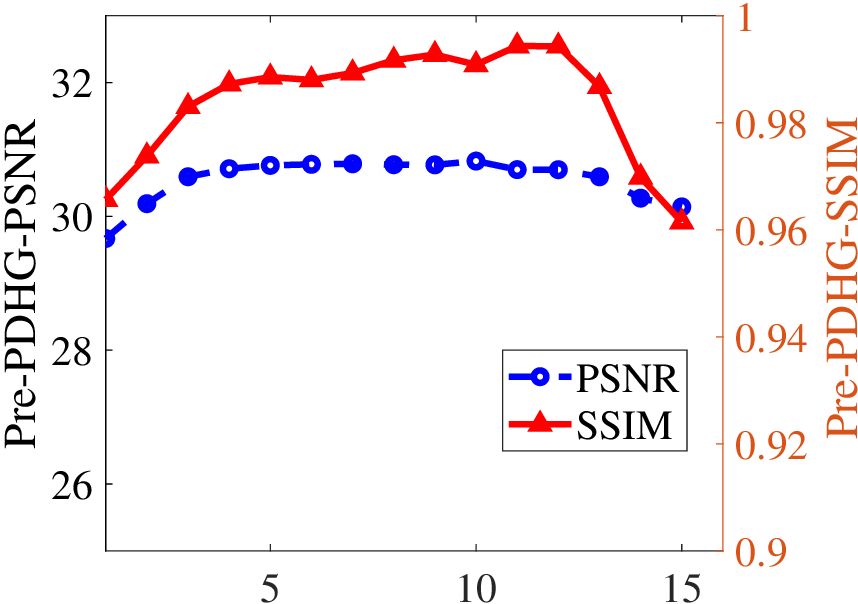}
		\caption{ $\beta$}
		\label{betaPre}   
	\end{subfigure}

      \begin{subfigure}{0.3\linewidth}
 	\centering
		\includegraphics[width=0.9\linewidth]{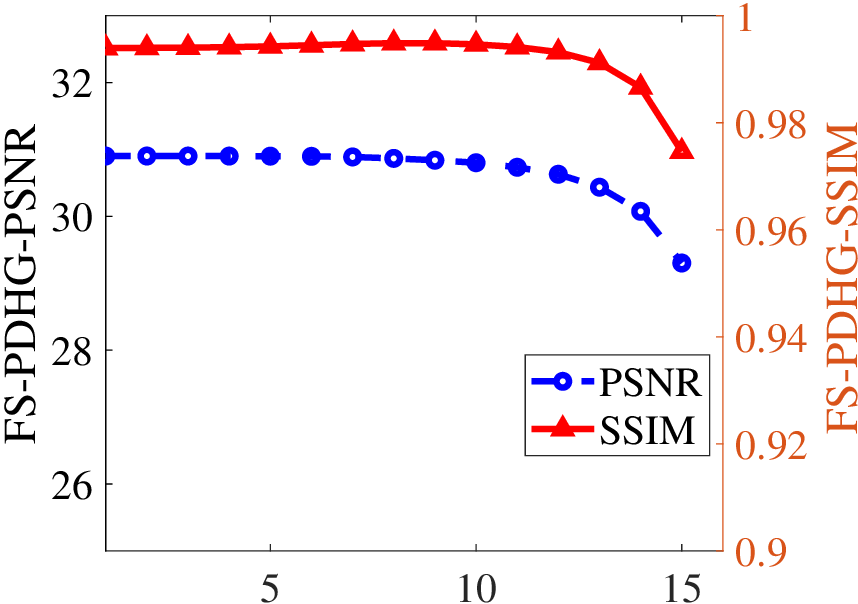}
		\caption{ $\eta$ }
		\label{etaFS}   
	\end{subfigure}
     \centering
	\begin{subfigure}{0.3\linewidth}
		\centering		\includegraphics[width=0.9\linewidth]{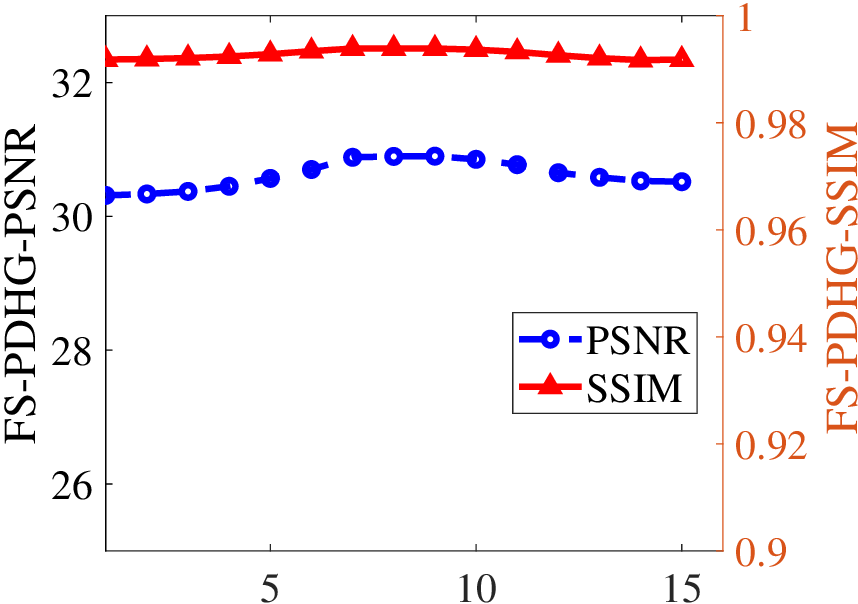}
		\caption{ $\tau$}
		\label{tauFS}   
	\end{subfigure}
	\centering
	\begin{subfigure}{0.3\linewidth}
		\centering
		\includegraphics[width=0.9\linewidth]{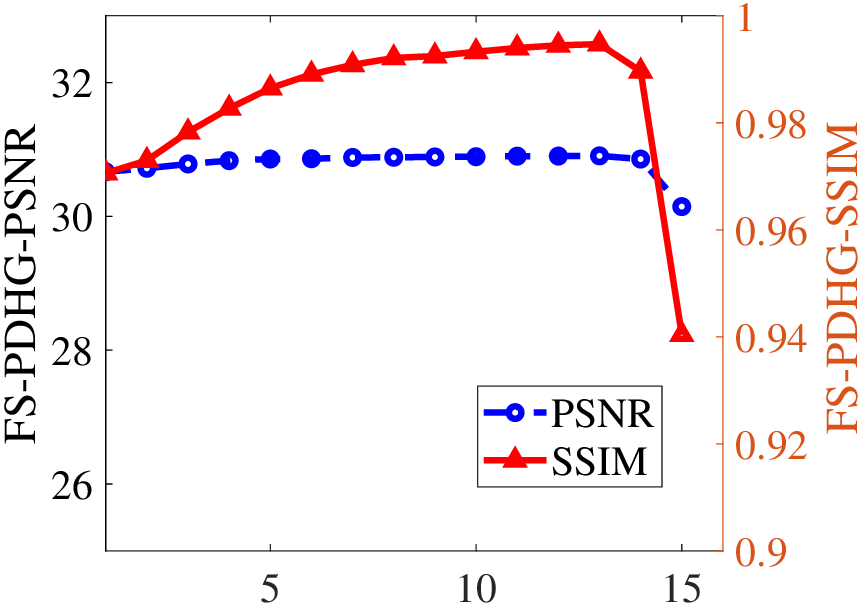}
		\caption{ $\beta$}
		\label{betaFS}   
	\end{subfigure}
 
    \centering
	\begin{subfigure}{0.3\linewidth}
		\centering
		\includegraphics[width=0.9\linewidth]{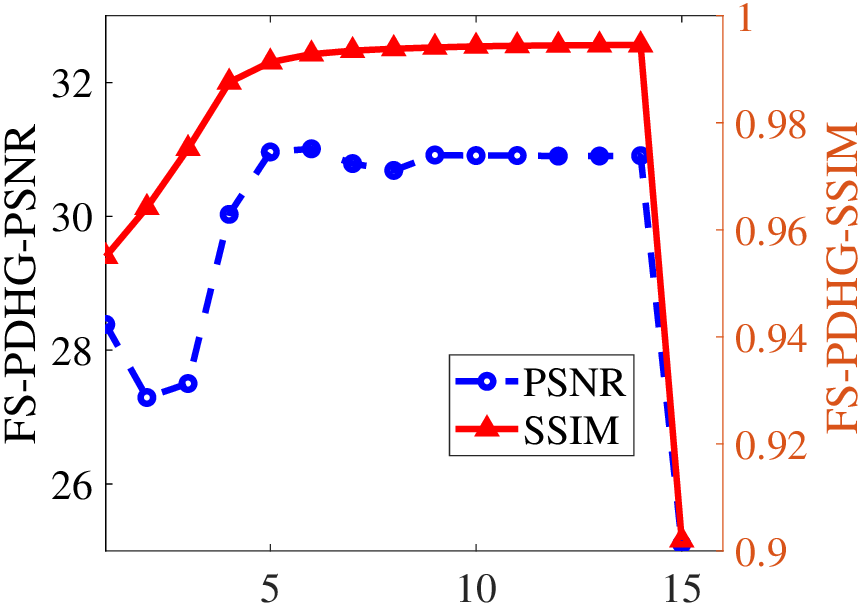}
		\caption{$\rho$ }
		\label{rho}   
	\end{subfigure}
	\begin{subfigure}{0.3\linewidth}
		\centering
		\includegraphics[width=0.9\linewidth]{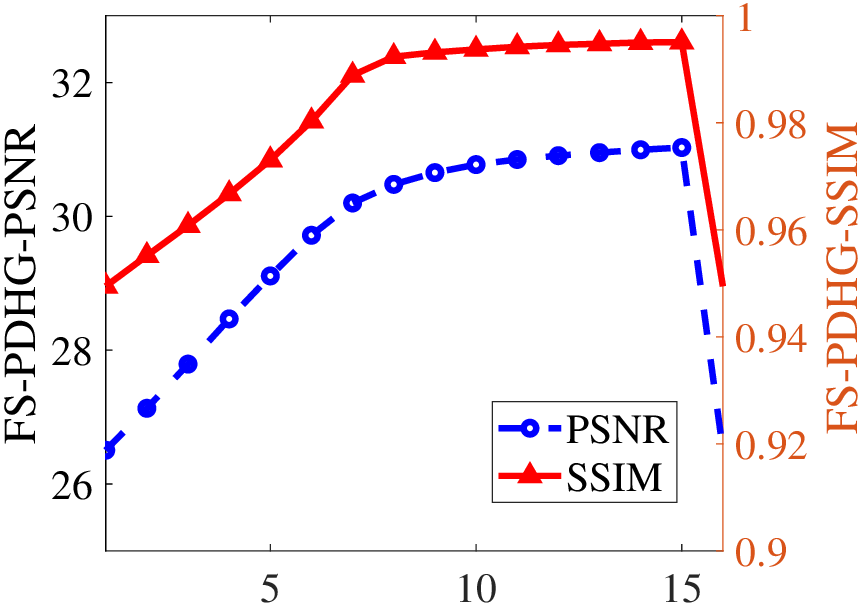}
		\caption{$\sigma_1$ }
		\label{sigma1}   
	\end{subfigure} 
    \centering
	\begin{subfigure}{0.3\linewidth}
		\centering
		\includegraphics[width=0.9\linewidth]{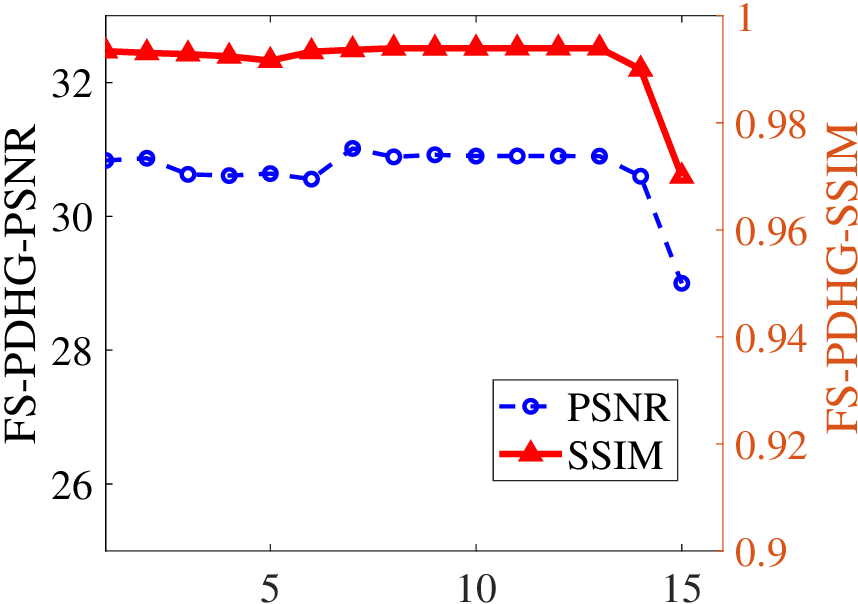}
		\caption{$\sigma_2$ }
		\label{sigma2}   
	\end{subfigure}
	\caption{Parameter optimization results for the NCAT phantom. The first row depicts the PSNR and SSIM values obtained under varying $\eta,\tau,$ and $\beta$ utilizing the Pre-PDHG algorithm. The second and third rows (Figures (d)-(i)) illustrate the performance of the FS-PDHG algorithm under adjustments to the parameters $\eta,\tau,\beta,\rho,\sigma_1,$ and $\sigma_2$, respectively.}	
 \label{paracurve image}
\end{figure}

Experimental tests in Figures \ref{etaPre} and \ref{etaFS} show that $\eta$ with $c_0=5 \times 10^{-6}$ is relatively insensitive when it is small, but the results demonstrate a threshold effect, whereby performance declines beyond a determinable upper bound. Besides, choose a proper $\eta =10^{-4}$ as defaulted hereafter.  Now let us consider the algorithm parameters, which include $\gamma$, $\tau$ and $\beta$. 
For the proposed Pre-PDHG, firstly, it can be observed by \eqref{gamma condition} that $\gamma$ is highly dependent on  $\lambda$ and $\mathcal P$, thus we do not provide a specific range. For $\tau$ and $\beta$ with $k_0=-5,c_0=5\times 10^{-3}$, 
as depicted in Figures \ref{tauPre}, \ref{tauFS}, \ref{betaPre} and \ref{betaFS}, the values do not exhibit significant variation with respect to these two parameters, empirically remaining within the approximate ranges of 0.01 and 5, respectively. Similarly, the parameter results for  $\tau$ and $\beta$ using the FS-PDHG algorithm exhibit comparable trends to those previously described. However, excessively small $\beta$ values may induce oscillations. Therefore, we adjust its range to center around 50, setting the default $k_0=1$, while maintaining the same empirically-derived $\tau$. 
For the remaining FS-PDHG parameters, we uniformly set $r_0$ to 0.25. The parameter $\rho$, with associated constants $c_0 = 1.2 \times 10^{-5}$ and $k_0 = 5.5$. Analysis indicates that excessively small $\rho$ results in slow convergence, while excessively large $\rho$ induces oscillations and failure to converge, as depicted in Figure \ref{rho}. 
Based on the analysis, the recommended viable range for $\rho$ is [0.001, 0.01], with a default value of 0.003. Utilizing the same settings as $\rho$, we conduct an exploration of the effects induced by varying $\sigma_1$, as depicted in Figure \ref{sigma1}. The results demonstrate that $\sigma_1$ exhibits relatively low sensitivity within a certain range. However excessive magnitudes lead to convergence failure. Empirically, we select $\sigma_1$  in the range of [0.001, 0.01]. Regarding $\sigma_2$ (see figure \ref{sigma2}), the associated constants are defined as $c_0 = 10^{-1}$ and $k_0 = 10$, respectively. Estimation reveals an optimal $\sigma_2$  centered around 300.

\subsection{Reconstructed results  and comparisons} 
\begin{figure}[htbp]
     \centering
	\begin{subfigure}{0.32\linewidth}
		\centering
		\includegraphics[width=0.9\linewidth]{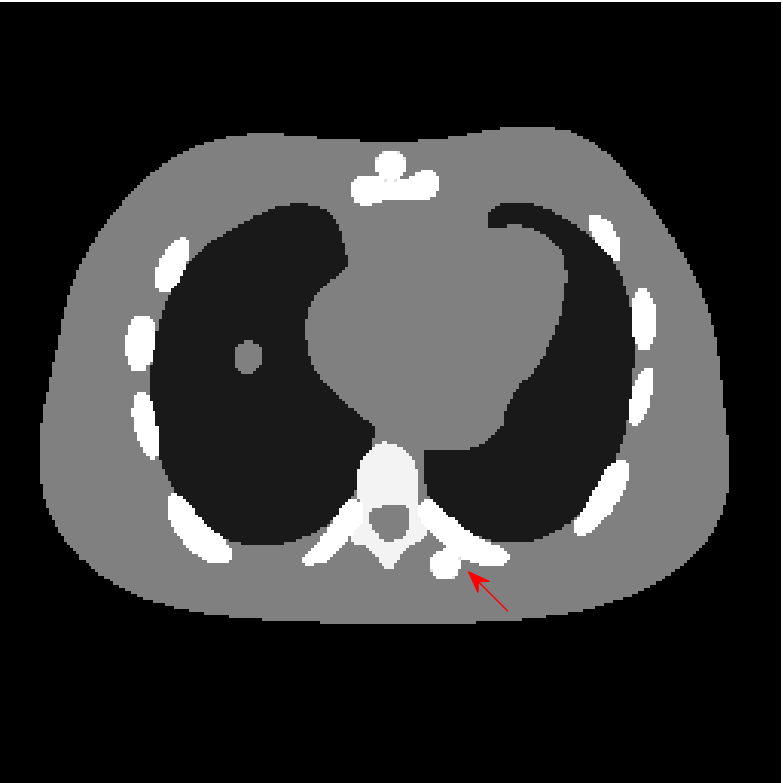}
		\caption{Original}
		\label{ncat truth}   
	\end{subfigure}
	\centering
	\begin{subfigure}{0.32\linewidth}
		\centering
		\includegraphics[width=0.9\linewidth]{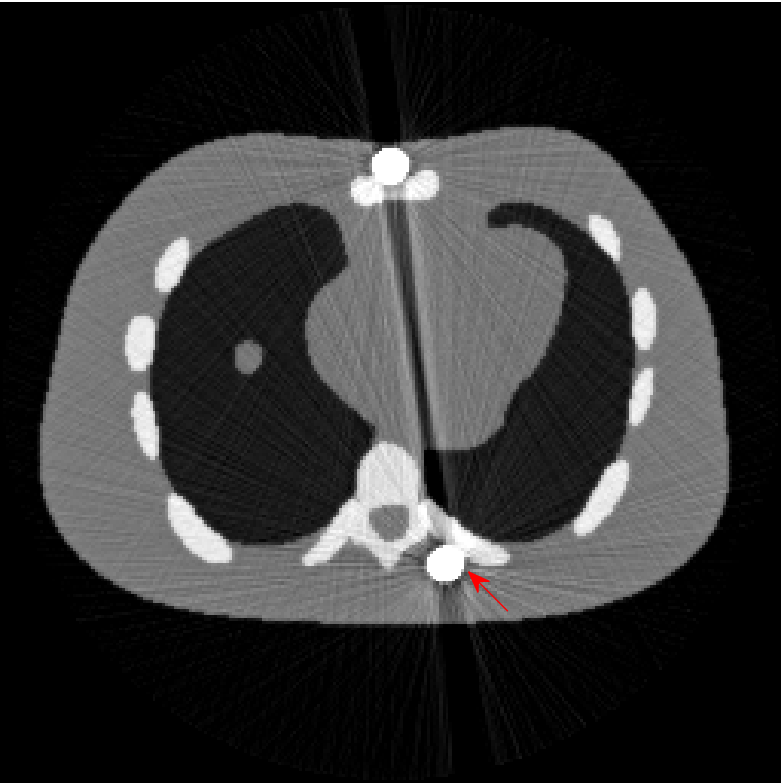}
		\caption{CG}
		\label{cg ncat}   
	\end{subfigure} 
    \centering
	\begin{subfigure}{0.32\linewidth}
		\centering
		\includegraphics[width=0.9\linewidth]{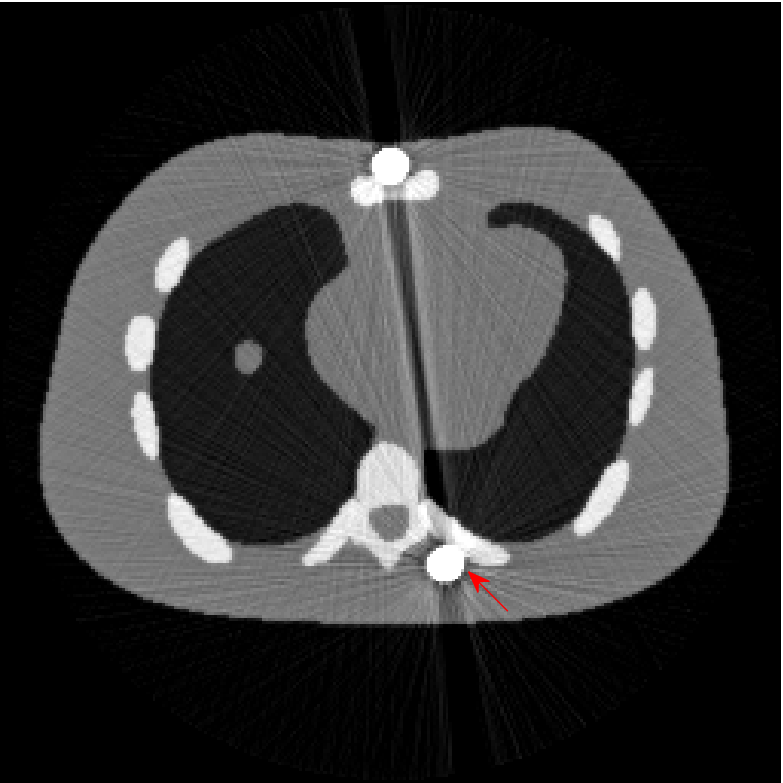}
		\caption{BCMAR}
		\label{corr ncat}   
	\end{subfigure}	
 
	\centering
	\begin{subfigure}{0.32\linewidth}
		\centering
		\includegraphics[width=0.9\linewidth]{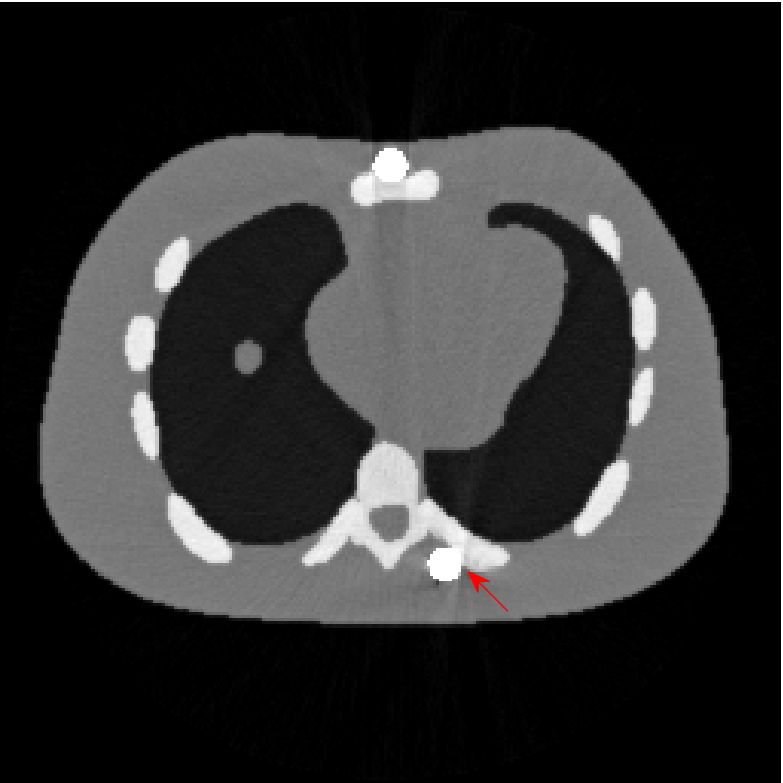}
		\caption{NMAR}
		\label{nmar ncat}   
	\end{subfigure}
 \centering
	\begin{subfigure}{0.32\linewidth}
		\centering
		\includegraphics[width=0.9\linewidth]{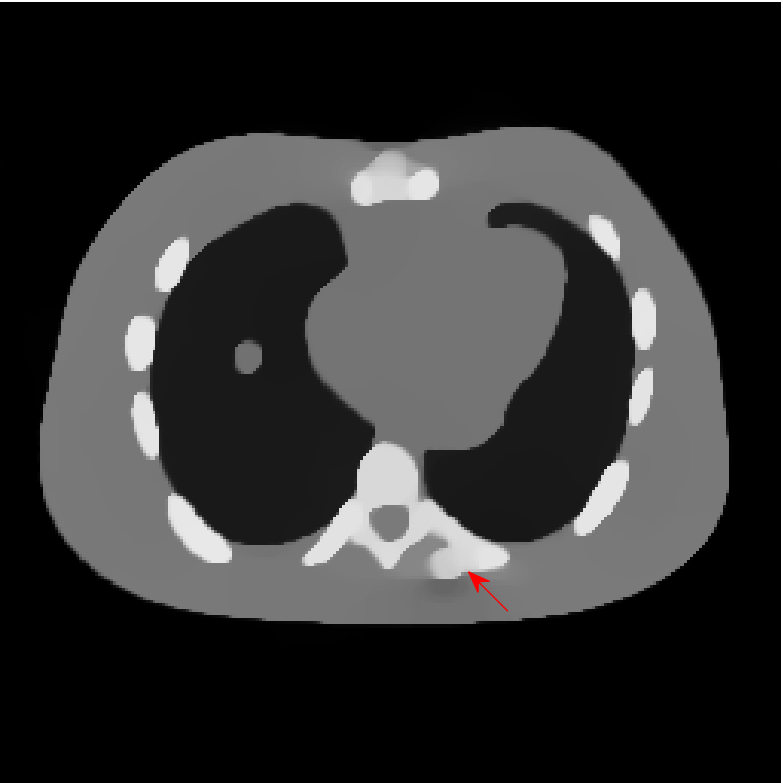}
		\caption{TV-MAR}
		\label{tv M ncat}   
	\end{subfigure}
 \centering
	\begin{subfigure}{0.32\linewidth}
		\centering
		\includegraphics[width=0.9\linewidth]{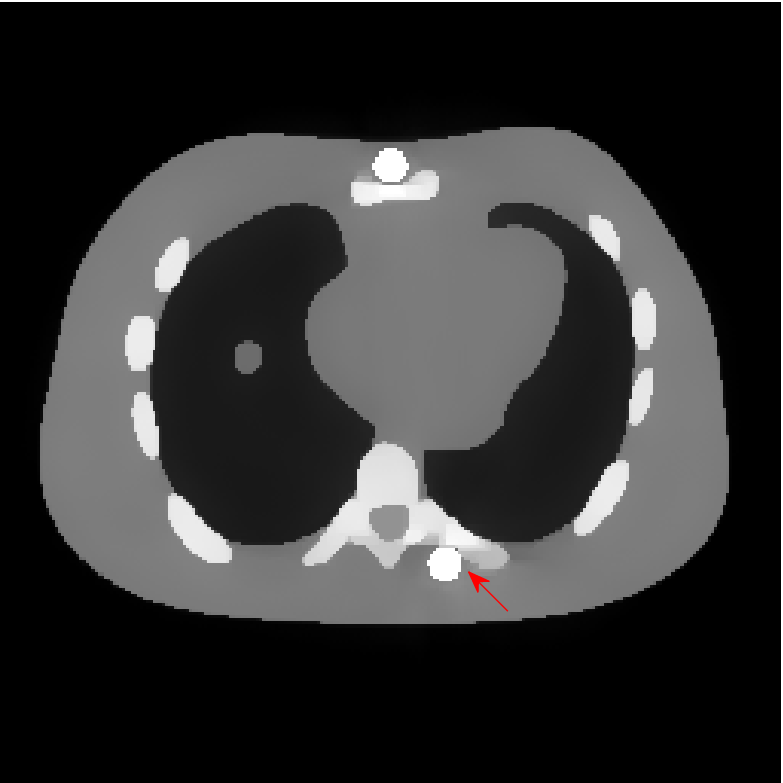}
		\caption{TV-TV inpainting}
		\label{TVTV ncat}   
	\end{subfigure}
 
	\begin{subfigure}{0.32\linewidth}
		\centering
		\includegraphics[width=0.9\linewidth]{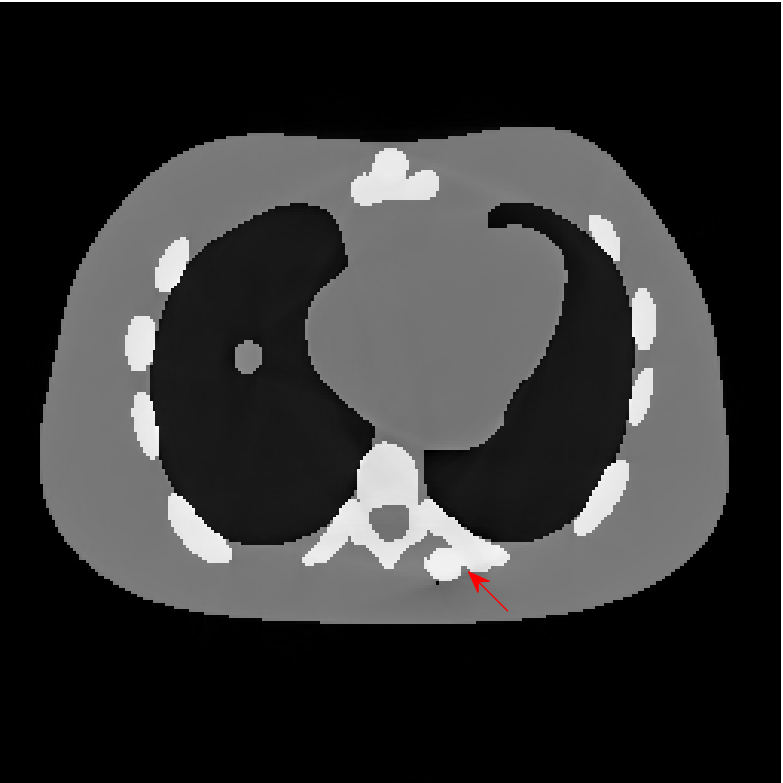}
		\caption{Reweighted JSR}
		\label{RJSR ncat}   
	\end{subfigure} 
 \centering
	\begin{subfigure}{0.32\linewidth}
		\centering
		\includegraphics[width=0.9\linewidth]{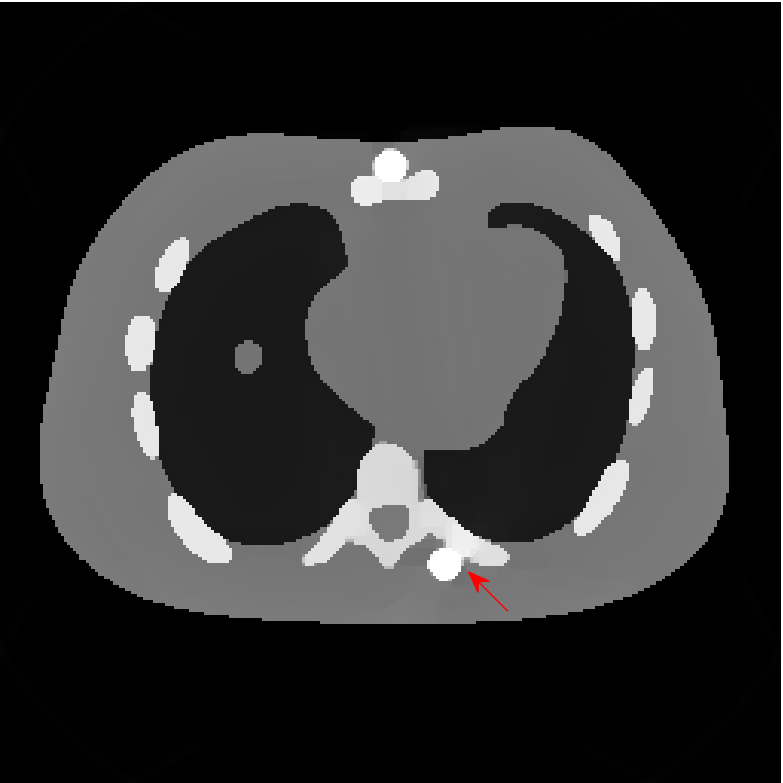}
		\caption{Pre-PDHG}
		\label{pre ncat}   
	\end{subfigure}
     \centering
	\begin{subfigure}{0.32\linewidth}
		\centering		\includegraphics[width=0.9\linewidth]{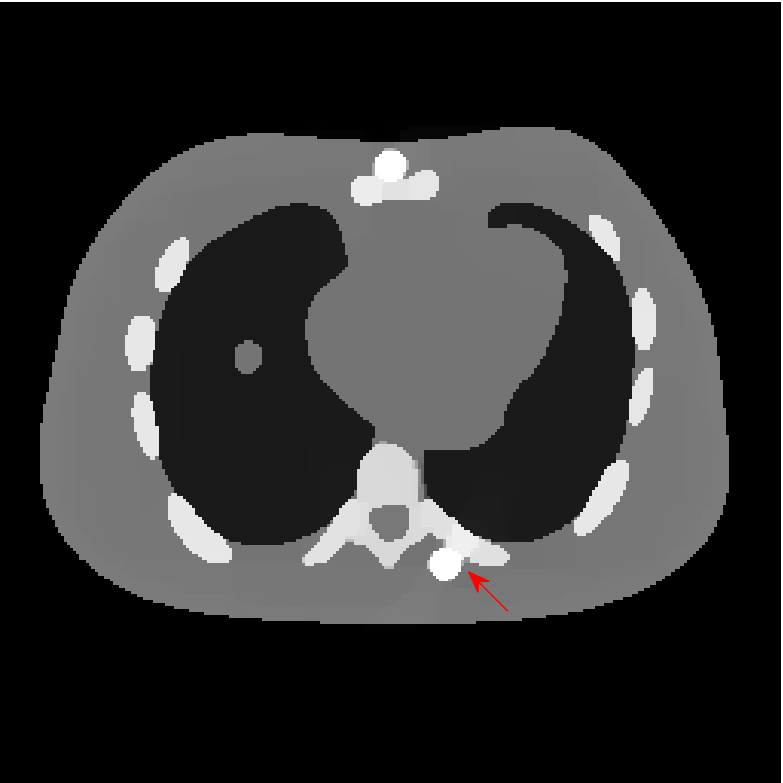}
		\caption{FS-PDHG}
		\label{split ncat}   
	\end{subfigure}
	\caption{Comparison of reconstructed NCAT phantom. (a) Phantom image, (b) CG-reconstruction of the measured projection, (c)-(i) corrected images utilizing BCMAR, NMAR, TV-MAR, TV-TV inpainting, Reweighted JSR, Pre-PDHG and FS-PDHG, respectively (4600 HU window, 1300 HU level).}	\label{ncat}
\end{figure}

\begin{figure}[htbp]
     \centering
	\begin{subfigure}{0.32\linewidth}
		\centering
		\includegraphics[width=0.9\linewidth]{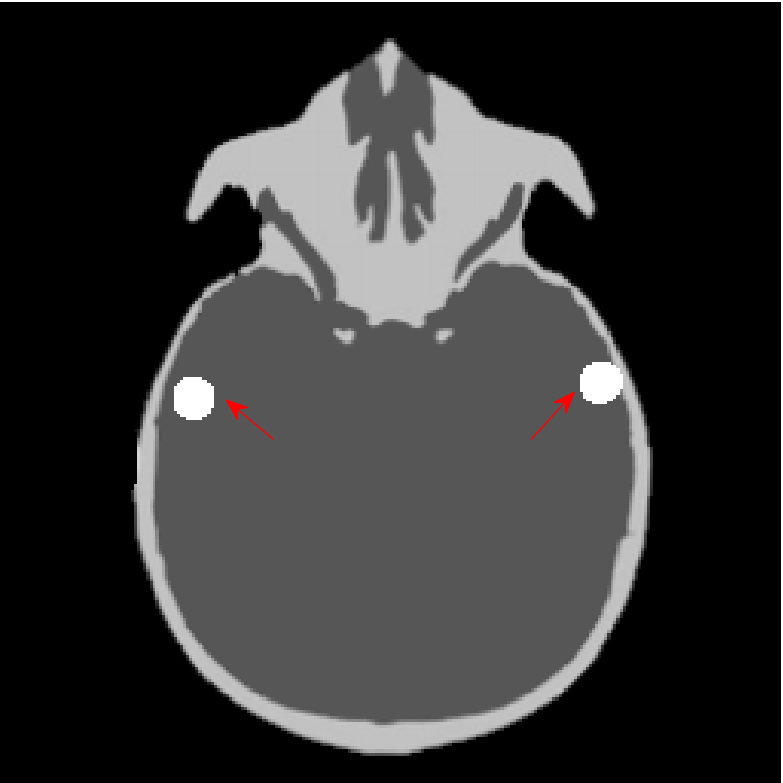}
		\caption{orginal}
		\label{truth ano}   
	\end{subfigure}
	\centering
	\begin{subfigure}{0.32\linewidth}
		\centering
		\includegraphics[width=0.9\linewidth]{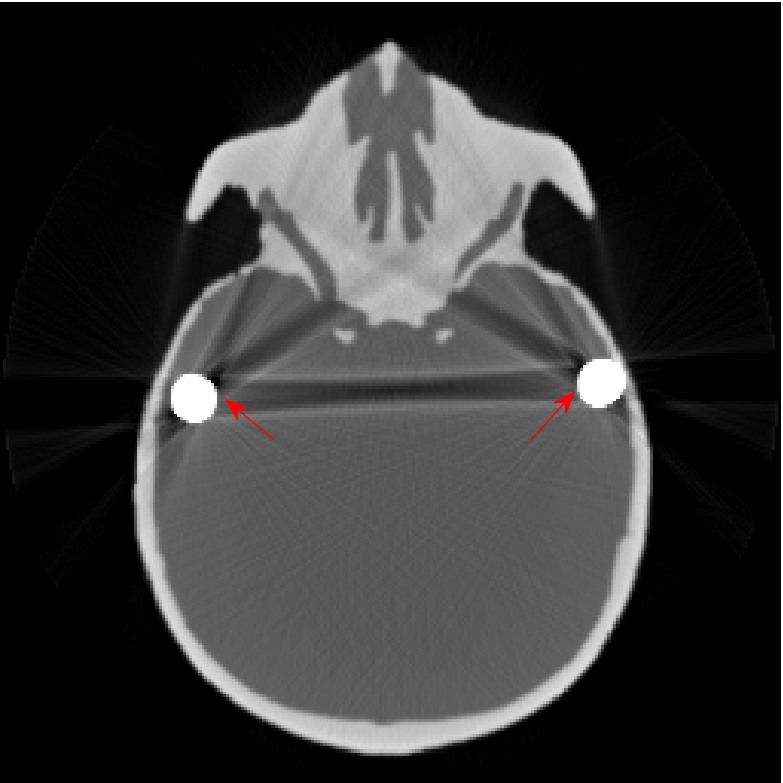}
		\caption{CG}
		\label{cg ano}   
	\end{subfigure}
    \centering
	\begin{subfigure}{0.32\linewidth}
		\centering
		\includegraphics[width=0.9\linewidth]{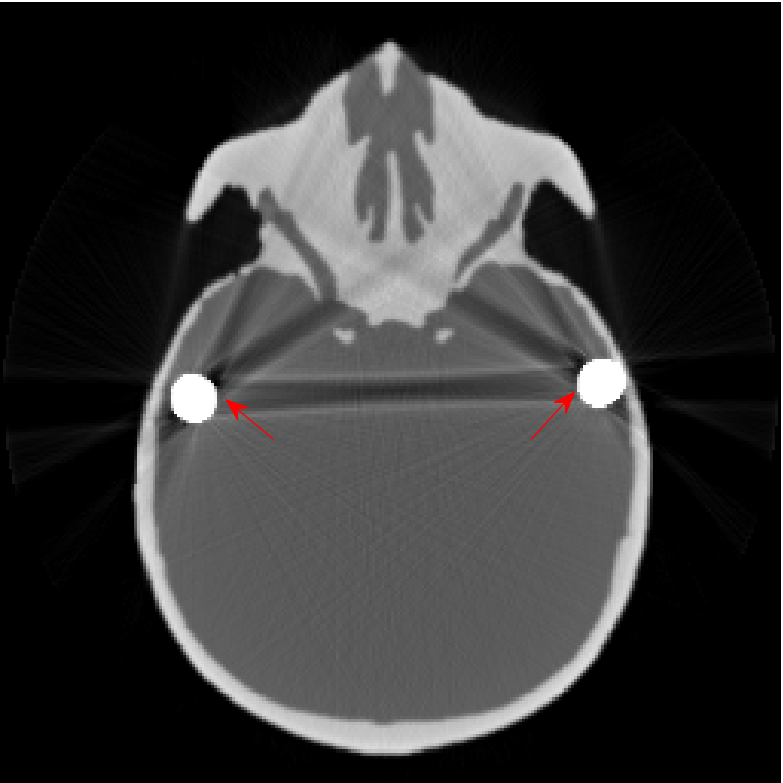}
		\caption{BCMAR}
		\label{corr ano}   
	\end{subfigure}
 
	\centering
	\begin{subfigure}{0.32\linewidth}
		\centering
		\includegraphics[width=0.9\linewidth]{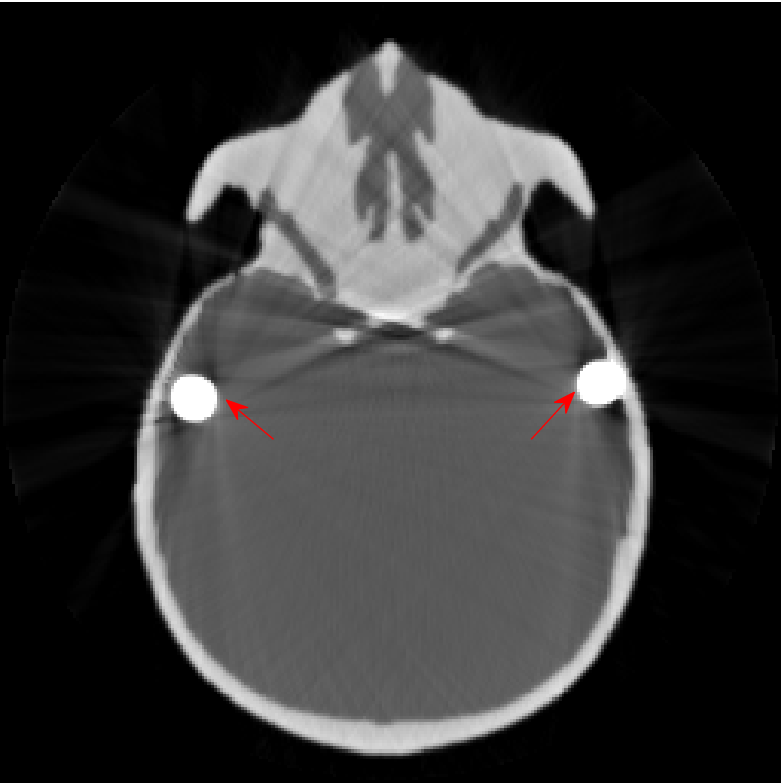}
		\caption{NMAR}
		\label{nmar ano}   
	\end{subfigure}
		 \centering
	\begin{subfigure}{0.32\linewidth}
	\includegraphics[width=0.9\linewidth]{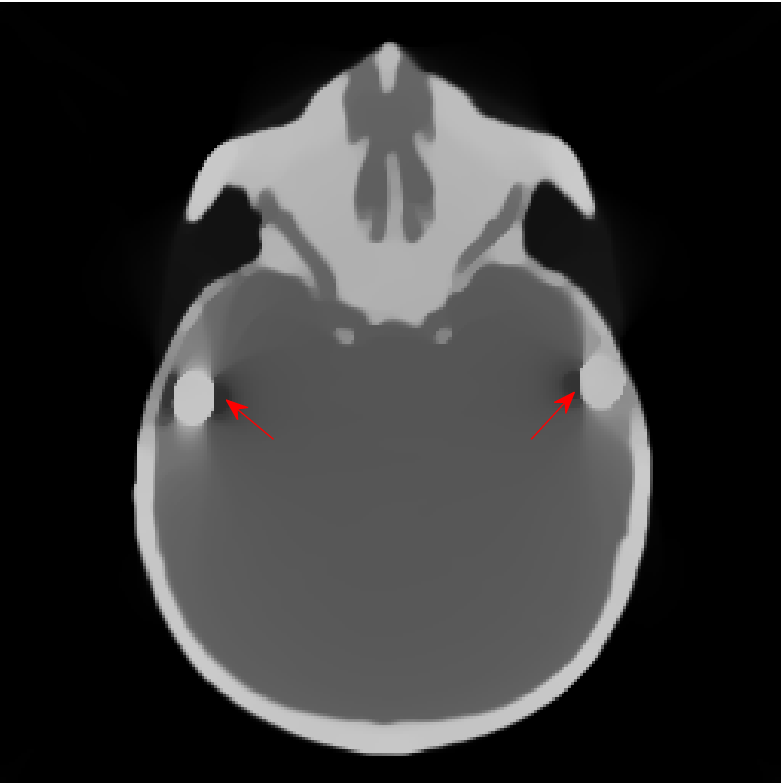}
		\caption{TV-MAR}
		\label{tv M ano}   
	\end{subfigure} 
 	\centering
	\begin{subfigure}{0.32\linewidth}
		\centering
	\includegraphics[width=0.9\linewidth]{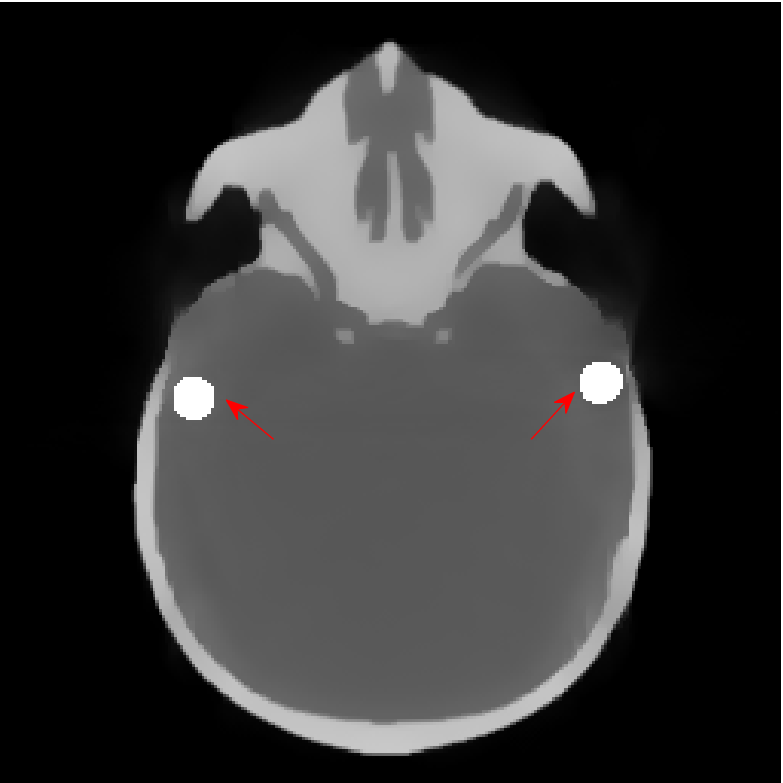}
		\caption{TV-TV inpainting}
		\label{tvtv ano}   
	\end{subfigure}
 
	\begin{subfigure}{0.32\linewidth}
		\centering		\includegraphics[width=0.9\linewidth]{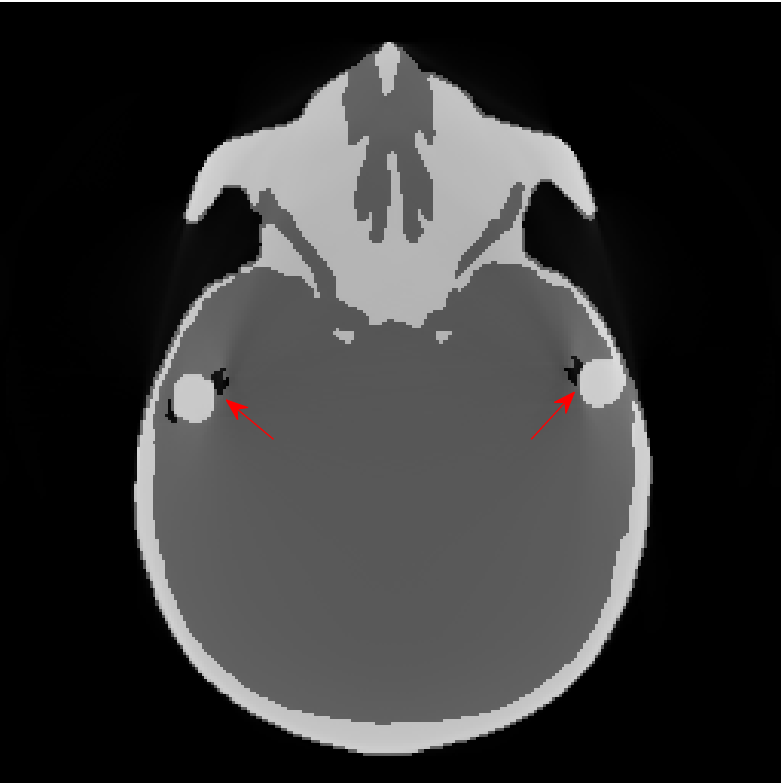}
		\caption{Reweighted JSR}
		\label{RJSR ano}   
	\end{subfigure} 
	\centering
	\begin{subfigure}{0.32\linewidth}
		\centering
		\includegraphics[width=0.9\linewidth]{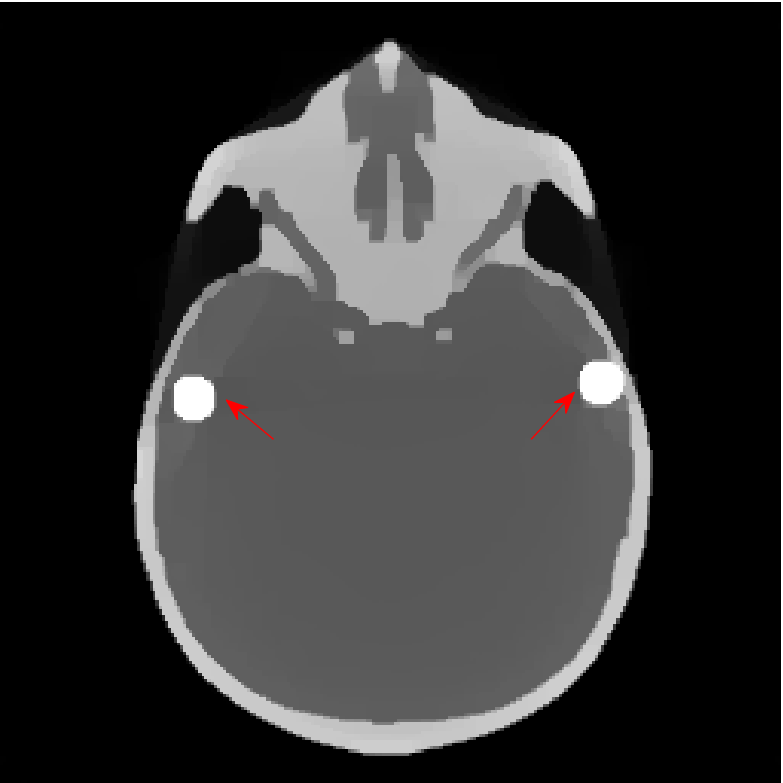}
		\caption{Pre-PDHG}
		\label{pre ano}   
	\end{subfigure}
     \centering
	\begin{subfigure}{0.32\linewidth}
		\centering
		\includegraphics[width=0.9\linewidth]{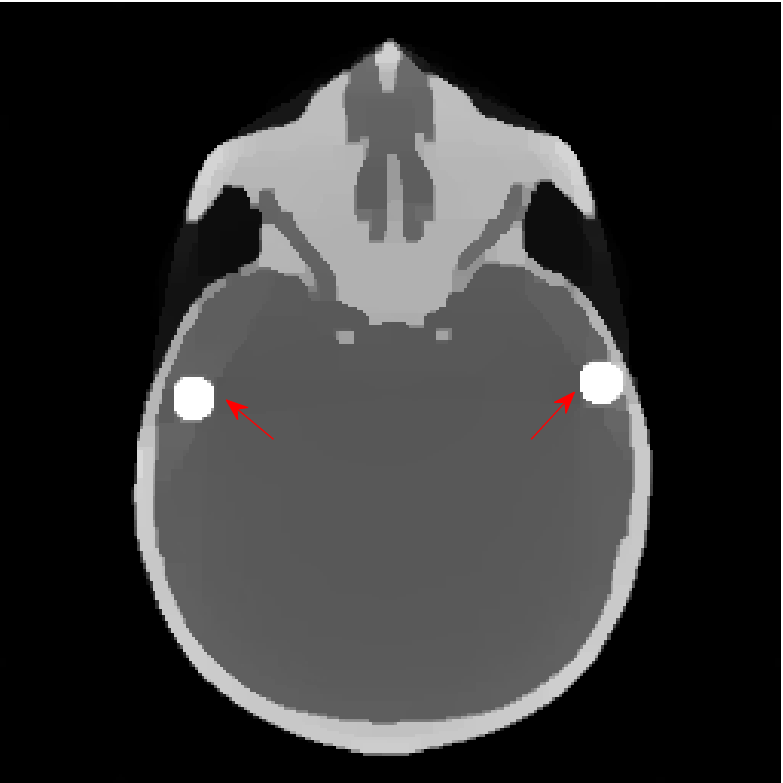}
		\caption{FS-PDHG}
		\label{split ano}   
	\end{subfigure}
	\caption{Image reconstruction of the head phantom by ground truth, CG, BCMAR, NMAR, TV-MAR, TV-TV inpainting, Reweighted JSR, Pre-PDHG and FS-PDHG (4200 HU window, 1100 HU level).}	\label{another}
\end{figure}

\begin{figure}[htbp]
	\centering
	\begin{subfigure}{0.32\linewidth}
		\centering
		\includegraphics[width=0.9\linewidth]{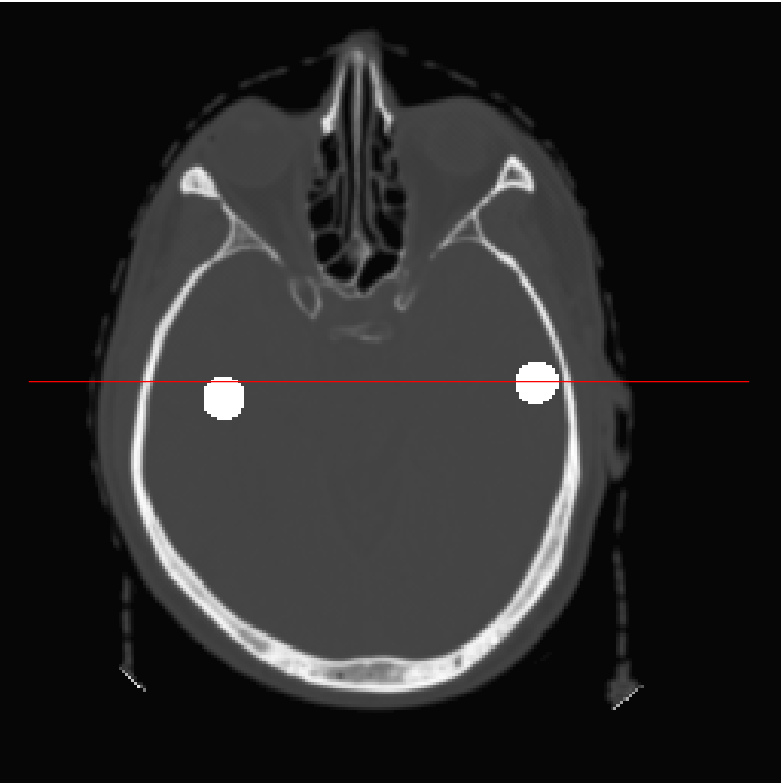}
		\caption{original}
		\label{new truth}   
	\end{subfigure}
 \centering
	\begin{subfigure}{0.32\linewidth}
		\centering
		\includegraphics[width=0.9\linewidth]{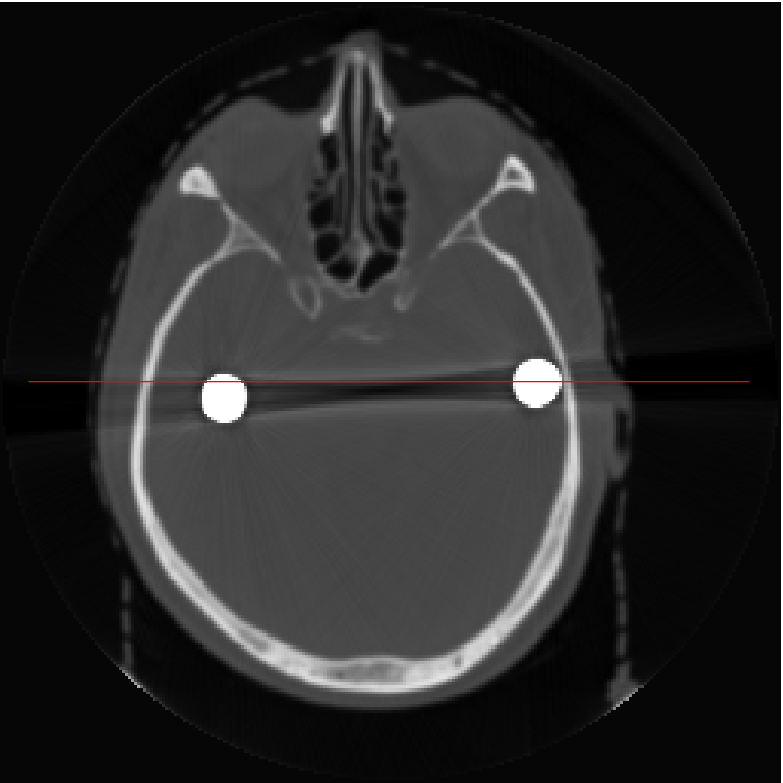}
		\caption{CG}
		\label{cg new}   
	\end{subfigure}
     \centering
	\begin{subfigure}{0.32\linewidth}
		\centering
		\includegraphics[width=0.9\linewidth]{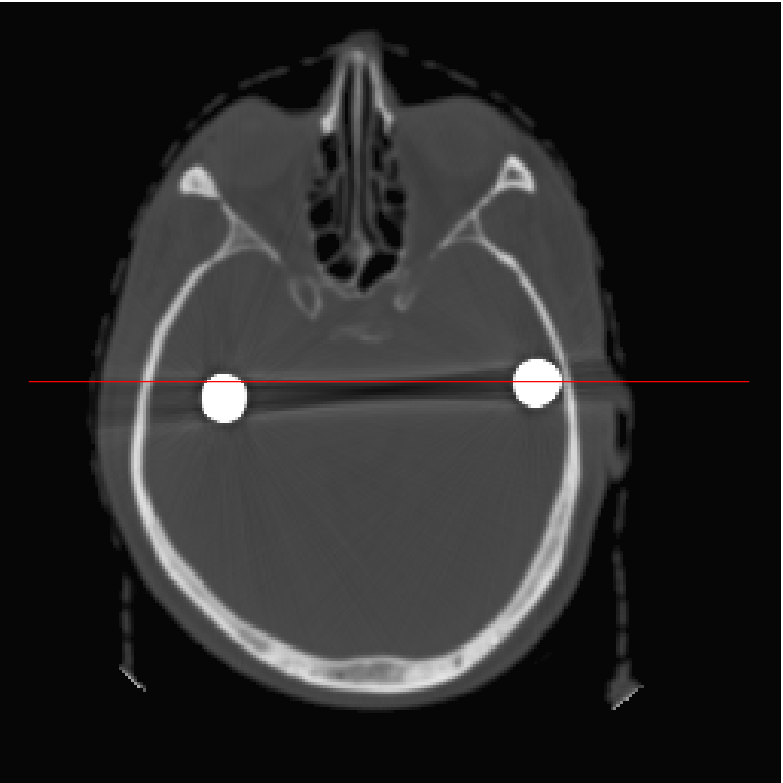}
		\caption{BCMAR}
		\label{corr new}   
	\end{subfigure}	
 
	\centering
	\begin{subfigure}{0.32\linewidth}
		\centering
		\includegraphics[width=0.9\linewidth]{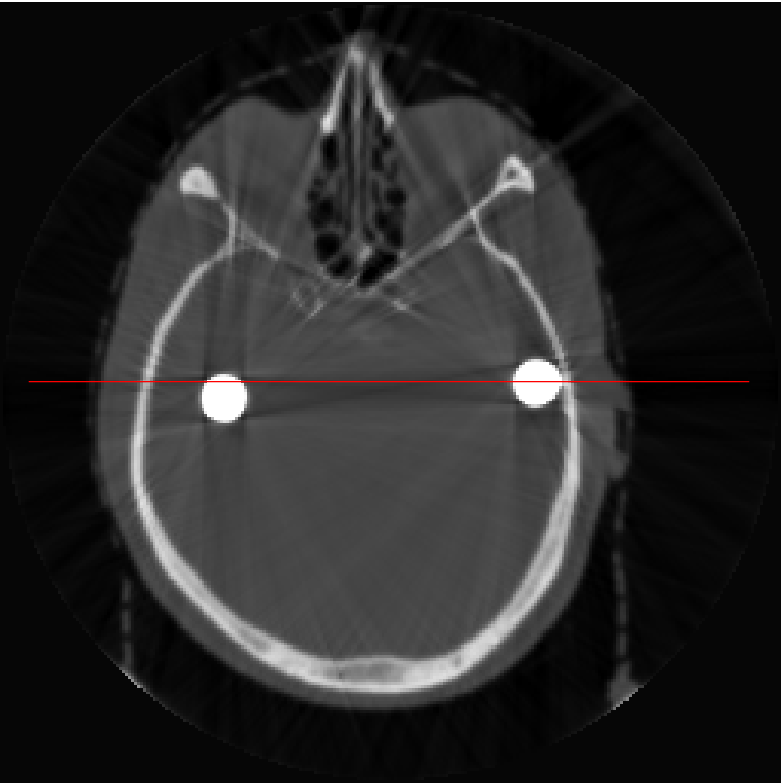}
		\caption{NMAR}
		\label{nmar new}   
	\end{subfigure}
 \centering
	\begin{subfigure}{0.32\linewidth}
		\centering
		\includegraphics[width=0.9\linewidth]{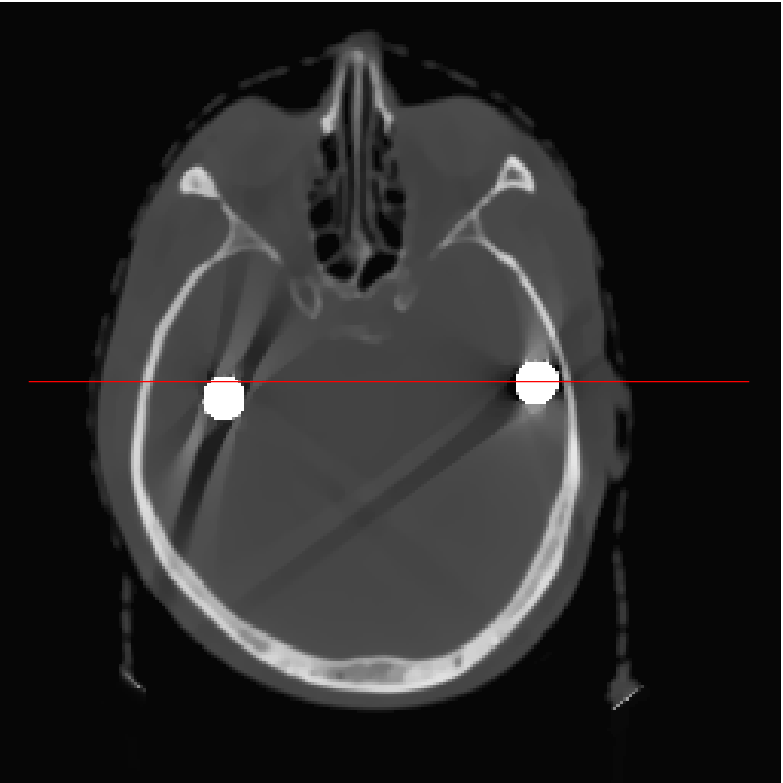}
		\caption{TV-MAR}
		\label{tv M new}   
	\end{subfigure} 
 \centering
	\begin{subfigure}{0.32\linewidth}
		\centering
		\includegraphics[width=0.9\linewidth]{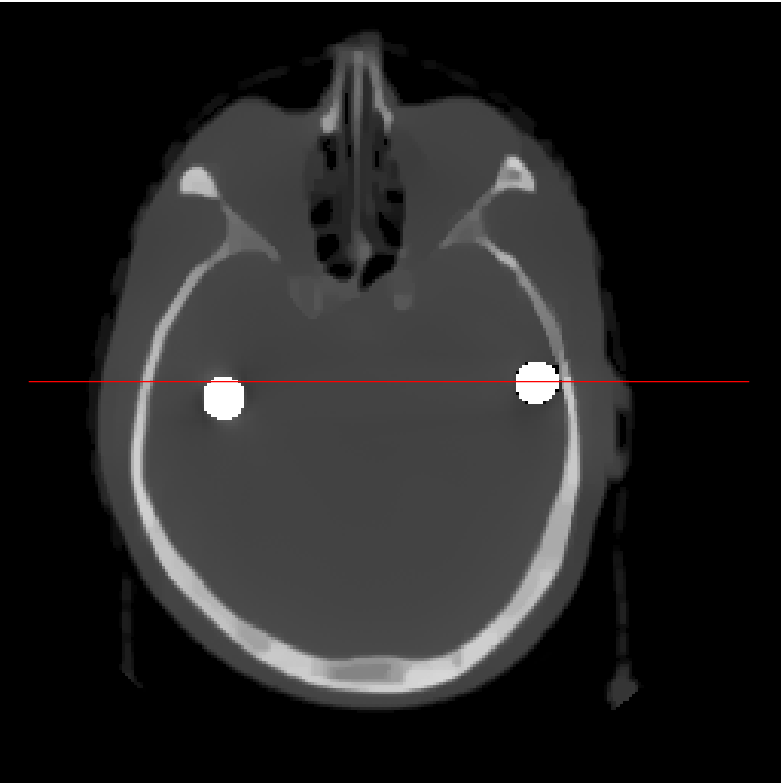}
		\caption{TV-TV inpainting}
		\label{TVTV}   
	\end{subfigure}	
 
 \centering
 	\begin{subfigure}{0.32\linewidth}
		\centering
		\includegraphics[width=0.9\linewidth]{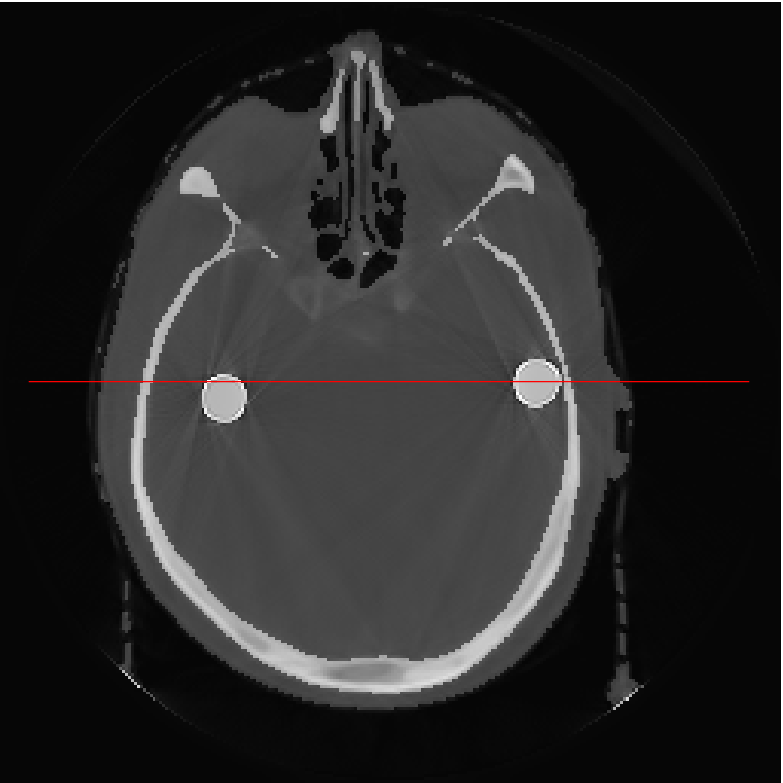}
		\caption{Reweighted JSR}
		\label{RJSR new}   
	\end{subfigure} 
	\begin{subfigure}{0.32\linewidth}
		\centering
		\includegraphics[width=0.9\linewidth]{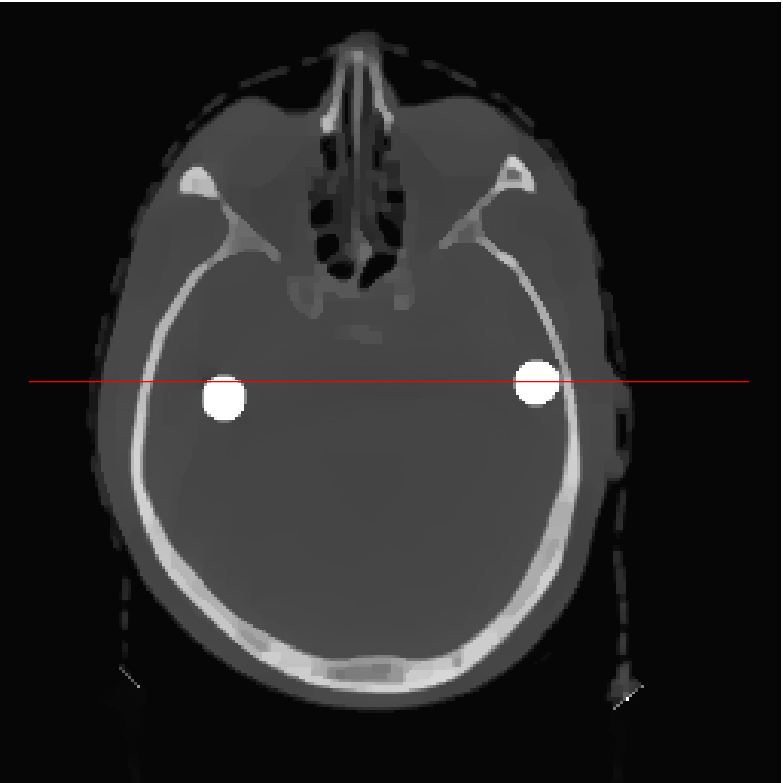}
		\caption{Pre-PDHG}
		\label{pre new}   
	\end{subfigure}
     \centering
	\begin{subfigure}{0.32\linewidth}
		\centering
		\includegraphics[width=0.9\linewidth]{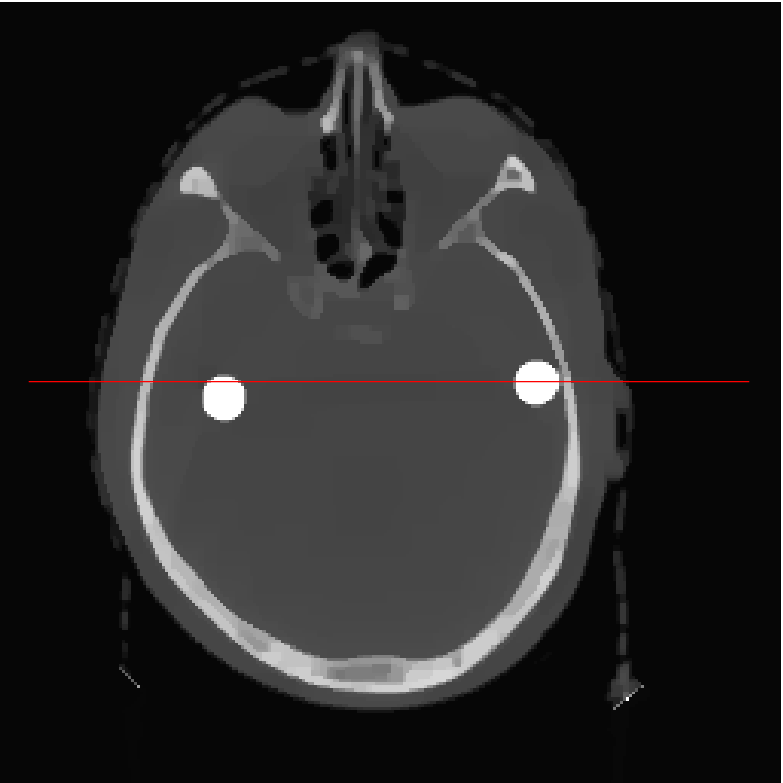}
		\caption{FS-PDHG}
		\label{split new}   
	\end{subfigure}
\caption{Comparison of numerical simulation results using the skull phantom (2800 HU window, 400 HU level).}
\label{new}
\end{figure}

In this subsection, we further show the results by the proposed algorithms and compare them with those by CG, BCMAR \cite{park2015metal}, NMAR \cite{meyer2010normalized}, TV-MAR (\ref{mask model}) (solving by \cite{chambolle2011first}), {\color{blue}TV-TV inpainting \cite{zhang2016computed}}, and Reweighted JSR \cite{zhang2018reweighted}. {\color{blue}The code of  Reweighted JSR is provided by the authors of \cite{zhang2018reweighted}, while the remaining codes have been developed by ourselves.} 
Results and comparisons between the proposed model and the above methods are shown in Figures \ref{ncat}, \ref{another} and \ref{new}, respectively. 

\begin{figure}[htbp]
     \centering
	\begin{subfigure}{0.32\linewidth}
		\centering		\includegraphics[width=0.9\linewidth]{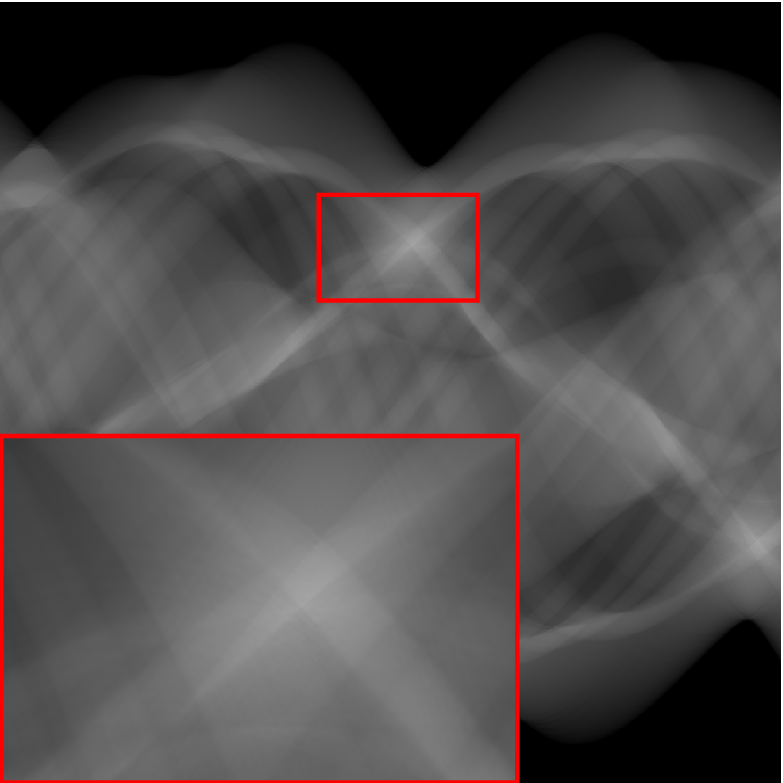}
		\caption{  }
		\label{hbsino}   
	\end{subfigure}
	\centering
	\begin{subfigure}{0.32\linewidth}
		\centering
		\includegraphics[width=0.9\linewidth]{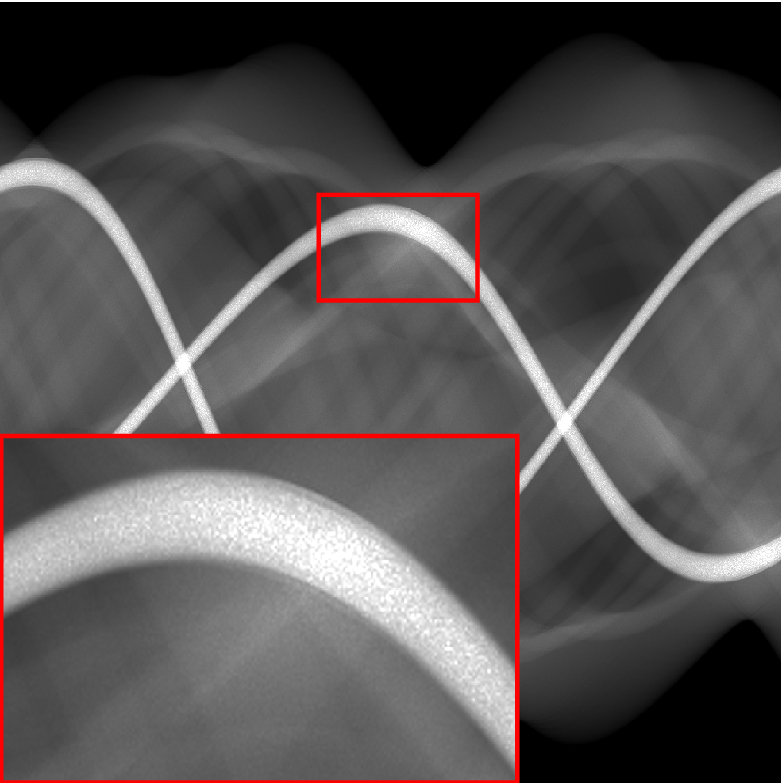}
		\caption{  }
		\label{Ysino}   
	\end{subfigure}
     \centering
	\begin{subfigure}{0.32\linewidth}
		\centering		\includegraphics[width=0.9\linewidth]{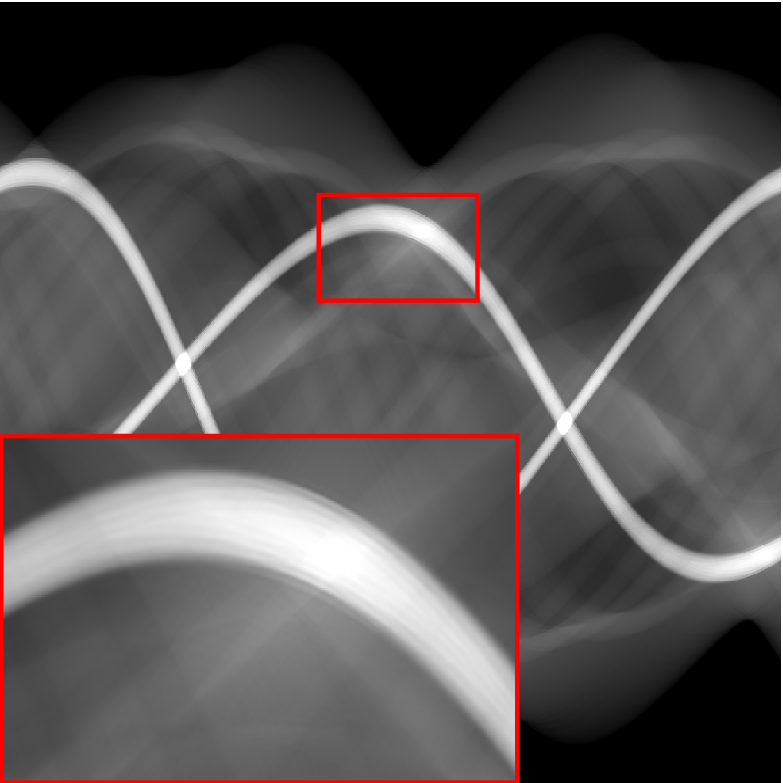}
		\caption{   }
		\label{Pusino}   
	\end{subfigure}
		\caption{(a) Reference sinogram. (b) Measured sinogram. (c) The sinogram result from FS-PDHG.}\label{NCAT sino}
\end{figure}

Figure \ref{cg ncat} illustrates that the NCAT phantom reconstructed via the CG approach contains substantial metal artifacts and noise. As shown in Figure \ref{corr ncat}, the BCMAR corrector \cite{park2015metal} fails to adequately address artifacts when high noise levels are present. While the NMAR approach (Figure \ref{nmar ncat}) admirably suppresses metal artifacts, blurriness around the metal regions exceeds that of TV-MAR. As depicted in Figure \ref{tv M ncat}, the TV-MAR reconstructed result manifests negligible metal artifacts and noise. However, the metal  components exhibit fusion with adjoining structures. {\color{blue}The TV-TV inpainting method yields improved results by reducing the excessive smoothing. However, there is still some diffusion present in the final output.} The reconstructed images obtained utilizing Reweighted JSR model and the proposed model not only display negligible artifacts, but also remarkably preserve salient details. {\color{blue}In Figure \ref{NCAT sino}, we present the sinograms of the reconstructed results \ref{split ncat}. Particularly in the magnified region, it is apparent that the proposed algorithm successfully suppresses the noise on the metal and bone traces depicted in Figure \ref{Pusino}.}

Figure \ref{another} displays the reconstructed images of the head phantom. This particular phantom poses greater challenges relative to the NCAT phantom, as it contains additional bone structures, yielding more prolific artifacts between the metals and bones alongside metal artifacts. The prior methods exhibit analogous consequences as the preceding examples. However, as observed in Figure \ref{nmar ano}, the outcomes of NMAR do not effectively recover regions near the bone, with artifacts still existing around the metals. Compared to the NCAT phantom, the enhanced artifacts arising between metals and bones detrimentally impact segmentation accuracy, engendering undesired structures in Figure \ref{RJSR ano} (structures around the red arrow) by Reweighted JSR. Nevertheless, akin to the NCAT phantom, the proposed method aptly corrects the metal artifacts as shown in Figures \ref{pre ano}-\ref{split ano}.  

\begin{figure}[htbp]
	\centering
	\begin{subfigure}{0.45\linewidth}
		\centering		\includegraphics[width=0.9\linewidth]{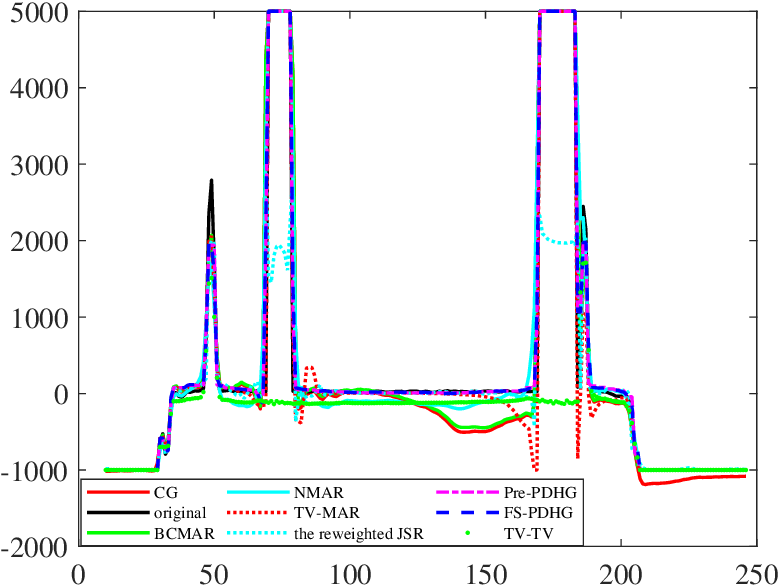}
		\caption{ $  $ }
		\label{entire}   
	\end{subfigure}
	\centering   
 \begin{subfigure}{0.45\linewidth}
		\centering		\includegraphics[width=0.9\linewidth]{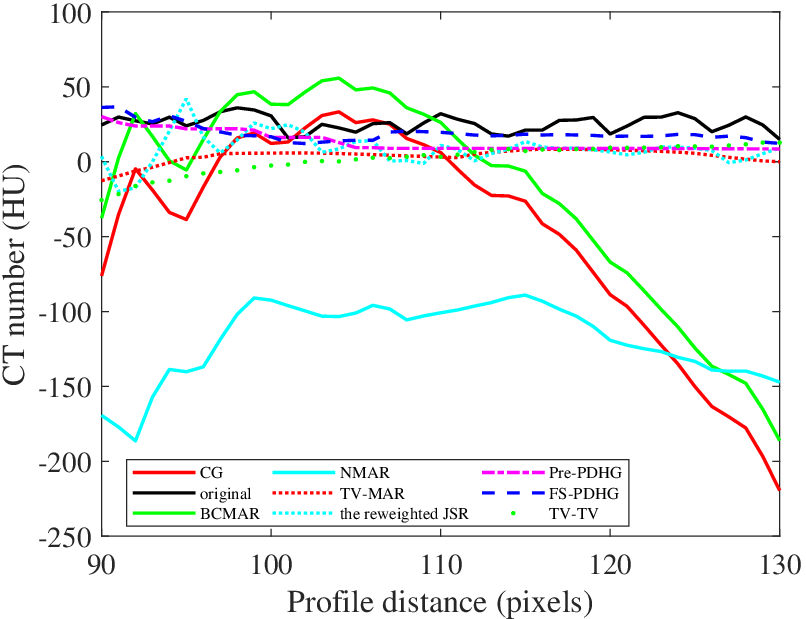}
		\caption{ $  $ }
		\label{local}   
	\end{subfigure}
    \caption{(a) The CT value curve corresponding to the red line representing the MAR results in Figure \ref{new}. (b) The local results within the range of [90, 130] in Figure (\ref{entire}).}
    \label{Skull HU curve}
\end{figure}

Figure \ref{new} exhibits the reconstructed skull phantom with more textures using various models. The result reconstructed by TV-MAR (Figure \ref{tv M new}) emerges significant artifacts. The inaccurate  pre-segmentation result leads to the correction of metal artifacts by Reweighted JSR (Figure \ref{RJSR new}) being insufficiently accurate. Analogous to prior phantom examples, the proposed algorithms persistently exhibit optimal comprehensive performance. {\color{blue}Furthermore, Figure \ref{Skull HU curve} illustrates the CT value distribution of the region of interest, providing further evidence that the results by our proposed method are more faithful to the ground truth. }

The quantitative assessments of reconstructed image quality, encompassing reconstruction error, SSIM, PSNR values, and computation time, are delineated in Table \ref{ssim psnr}. 
One readily discerns the enhancement conferred by Reweighted JSR and the proposed algorithms in terms of PSNR and SSIM, as compared to TV-MAR.
In comparison to Reweighted JSR, the proposed algorithms (Pre-PDHG and FS-PDHG) can reconstruct competitive outcomes in terms of PSNR and SSIM, given accurate pre-segmentation (as in the case of NCAT). For less precise pre-segmentation, they can reconstruct significantly superior results, garnering PSNR escalations of 0.8dB and 1.8dB for the head and skull cases, respectively.
 Moreover, in comparison to Reweighted JSR, our proposed FS-PDHG conspicuously curtails computational expenditure,  reducing the total computational cost to at most one-third.

\begin{table}\small
\caption{Reconstruction error, SSIM index, PSNR value, and the CPU time (sec.) for a single image of the NCAT, head and skull phantoms reconstructed by different algorithms, i.e., CG, BCMAR \cite{park2015metal}, NMAR \cite{meyer2010normalized}, TV-MAR (using the primal-dual algorithm), TV-TV inpainting \cite{zhang2016computed}, Reweighted JSR model \cite{zhang2018reweighted} and ours.}  
\begin{subtable}[t]{5in}
\renewcommand{\thesubtable}{(a)}
\centering
    \begin{tabular}{@{}lllll}
        \br
        Model           & Reconstruction error         & SSIM      & PSNR      &  Time\\
        \mr
        CG               & 0.1586                & 0.7107    &26.0206           &---\\
        BCMAR        & 0.1557                & 0.7124    &26.1511           & ---\\
        NMAR             & 0.1144                & 0.8921    &28.2581           & ---\\
        TV-MAR    & 0.1080                & 0.9758    &28.7206           &207\\
        TV-TV inpainting    & 0.1060                & 0.9789   & 28.8581          &245   \\
        Reweighted JSR             & 0.0829                & 0.9831    & 30.6438   &492\\
        Pre-PDHG         & 0.0839                & 0.9932    &30.8273           &438\\
        FS-PDHG       &{\bf{0.0825} }                 &{\bf{0.9947}}     &{\bf{30.9797}}           &{\bf{238}}\\
        \br       
    \end{tabular}
    \caption{NCAT}
    \label{NCAT}
\end{subtable}
\quad 
\begin{subtable}[t]{5in}
\renewcommand{\thesubtable}{(b)}
\centering
\begin{tabular}{@{}lllll}
        \br
        Model           & Reconstruction error         & SSIM            &PSNR       &  Time \\
        \mr
        CG               & 0.1530             &  0.7339       &27.8421      &---        \\
        BCMAR        & 0.1503              &  0.7388        & 28.0049      &--- \\
        NMAR             & 0.1308             & 0.7608       &29.0103    &--- \\
        TV-MAR    & 0.1332                & 0.9114          &28.6113   &164\\
        TV-TV inpainting    & 0.1217                & 0.9220   & 28.9188          &355   \\
        Reweighted JSR             & 0.1199                & 0.9085          &29.5254            &525\\
        Pre-PDHG         & 0.1109                & 0.9470          &30.2593  &454\\
        FS-PDHG       & {\bf{0.1093}}                & {\bf{0.9486} }          &{\bf{30.3540}}     &{\bf{344}}      \\
        \br      
    \end{tabular}
    \caption{Head}
    \label{head}
\end{subtable}
\quad 
\begin{subtable}[t]{5in}
\renewcommand{\thesubtable}{(c)}
\centering
\begin{tabular}{@{}lllll}
        \br
        Model           & Reconstruction error         & SSIM      & PSNR      &  Time \\
        \mr
        CG               & 0.1892                & 0.8511    &29.4922            &---     \\
        BCMAR        & 0.1808                & 0.9394    & 29.8744          & --- \\
        NMAR             & 0.1648                & 0.8550    &30.2571           & --- \\
        TV-MAR    & 0.1745                & 0.9565    & 29.8293          &290   \\
        TV-TV inpainting    & 0.1700                & 0.9575   & 29.8788          &400   \\
        Reweighted JSR             & 0.1676                & 0.9500    & 30.5401   &1104            \\
        Pre-PDHG         & 0.1323                & 0.9675    &32.1448           &1079 \\
        FS-PDHG       &{\bf{0.1299}}                &{\bf{0.9681}}     &{\bf{32.2954}}           &{\bf{388}}           \\
        \br      
    \end{tabular}
    \label{skull}
    \caption{Skull}
    \end{subtable}
\label{ssim psnr}
\end{table}
{\color{blue}\subsection{Experiments on real data}
To further illustrate the validity, we conducted a CT scan of a tooth enclosed in a plastic bottle with two titanium rods positioned on one side (Figure \ref{real}). The projection data are from a microscopic CT scanner independently developed by Capital Normal University.  Each projection view comprised 540 angles, with 752 detector bins employed. The distance from the rotation center to the ray source was 200mm, while the distance from the rotation center to the detector was 890mm. Additionally, the space between each detector bin was set at 0.1mm. The reconstructed image consisted of 256$\times$256 pixels.

Figure \ref{real comparison} displays reconstructed images obtained using various methods, including CG, BCMAR, NMAR, TV-MAR,  Reweighted JSR and the proposed FS-PDHG. These images, shown in this section, are presented in the grayscale range of [0, 0.3]. NMAR performs better than CG and BCMAR, in preserving image details and correcting artifacts, albeit with the presence of new artifacts near the metal. Furthermore, the results obtained from variation-based models exhibit significantly less noise. Our proposed method demonstrates the best overall performance in terms of preserving details and rectifying artifacts, especially near the teeth pointed by  red arrows in Figure \ref{real comparison}.
\begin{figure}
	\centering
	\begin{subfigure}{0.45\linewidth}
		\centering
		\includegraphics[width=0.9\linewidth]{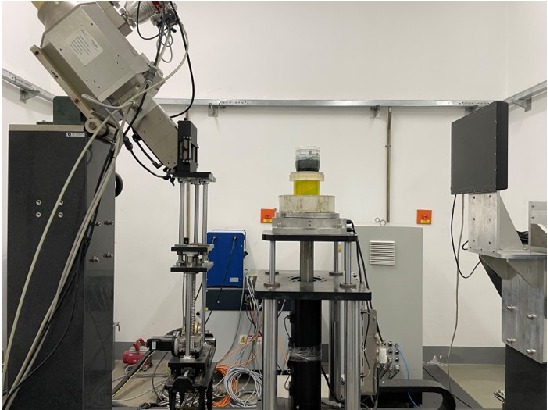}
		\caption{  }
		\label{real setting}   
	\end{subfigure}
     \centering
	\begin{subfigure}{0.45\linewidth}
		\centering		\includegraphics[width=0.9\linewidth]{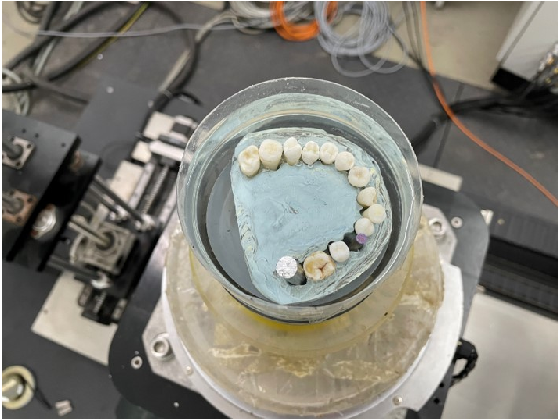}
		\caption{   }
		\label{tooth}   
	\end{subfigure}
		\caption{(a) Scanning device. (b) Tooth with metal on one side.}\label{real}
  \end{figure}
  
\begin{figure}
	\centering
	\begin{subfigure}{0.32\linewidth}
		\centering
		\includegraphics[width=0.9\linewidth]{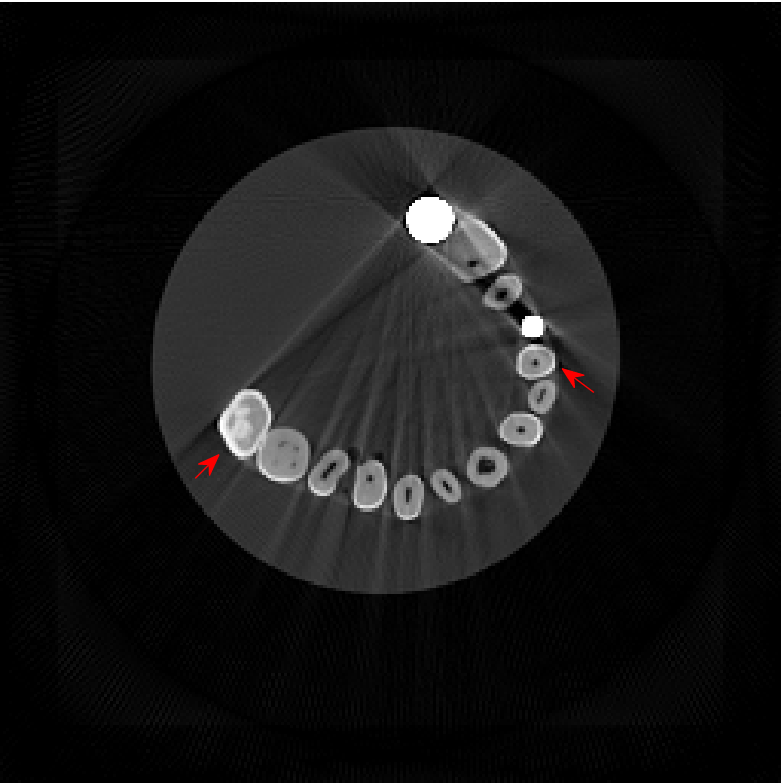}
		\caption{ CG }
		\label{real CG}   
	\end{subfigure}
 \centering
	\begin{subfigure}{0.32\linewidth}
		\centering
		\includegraphics[width=0.9\linewidth]{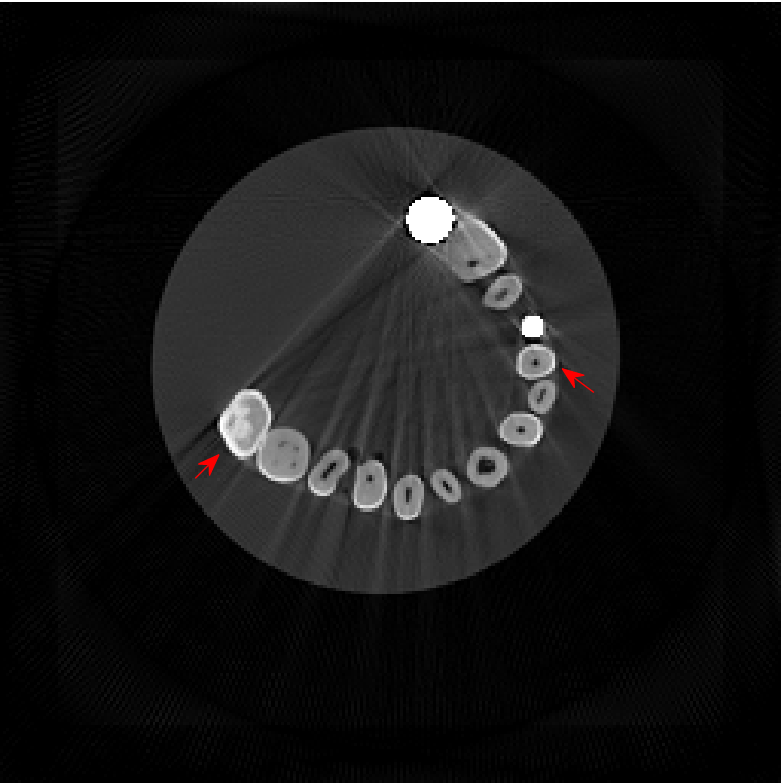}
		\caption{ BCMAR }
		\label{real BAMAR}   
	\end{subfigure}
     \centering
	\begin{subfigure}{0.32\linewidth}
		\centering		\includegraphics[width=0.9\linewidth]{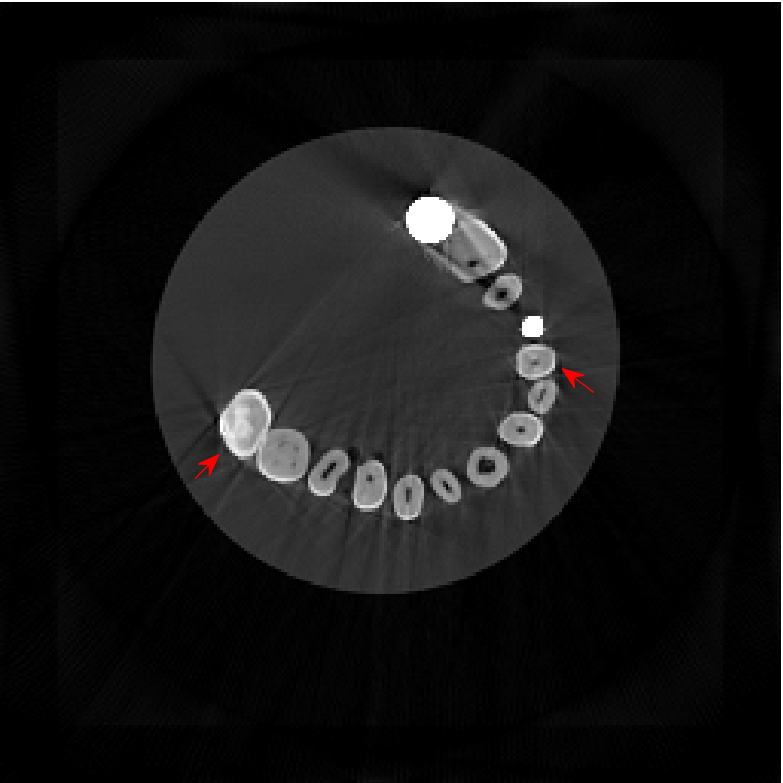}
		\caption{ NMAR  }
		\label{real NMAR}   
	\end{subfigure}

     \centering
	\begin{subfigure}{0.32\linewidth}
		\centering		\includegraphics[width=0.9\linewidth]{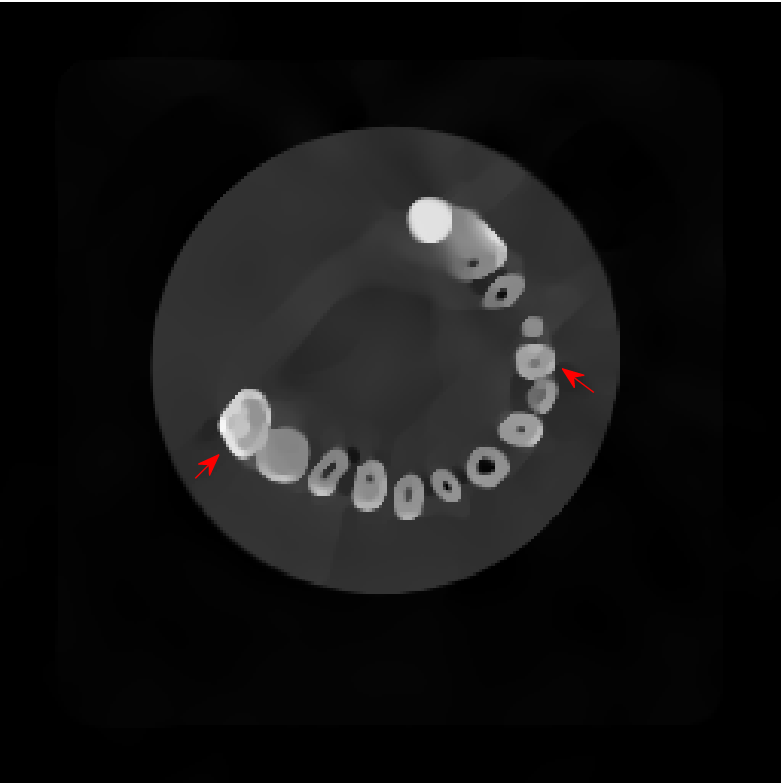}
		\caption{ TV-MAR  }
		\label{real TV}   
	\end{subfigure}
 \centering
	\begin{subfigure}{0.32\linewidth}
		\centering
		\includegraphics[width=0.9\linewidth]{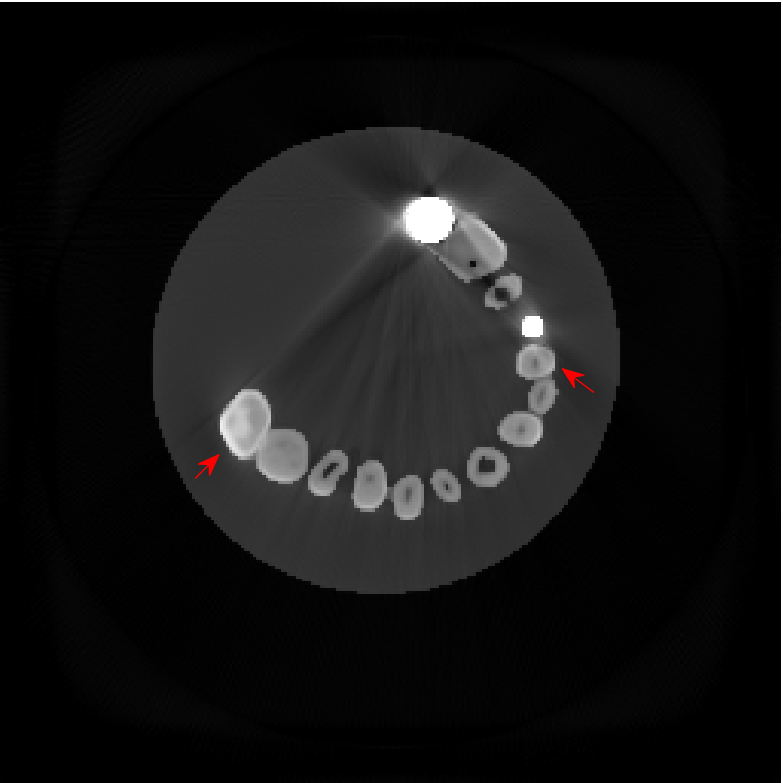}
		\caption{ Reweighted JSR }
		\label{real RJSR}   
	\end{subfigure}
 \centering
	\begin{subfigure}{0.32\linewidth}
		\centering
		\includegraphics[width=0.9\linewidth]{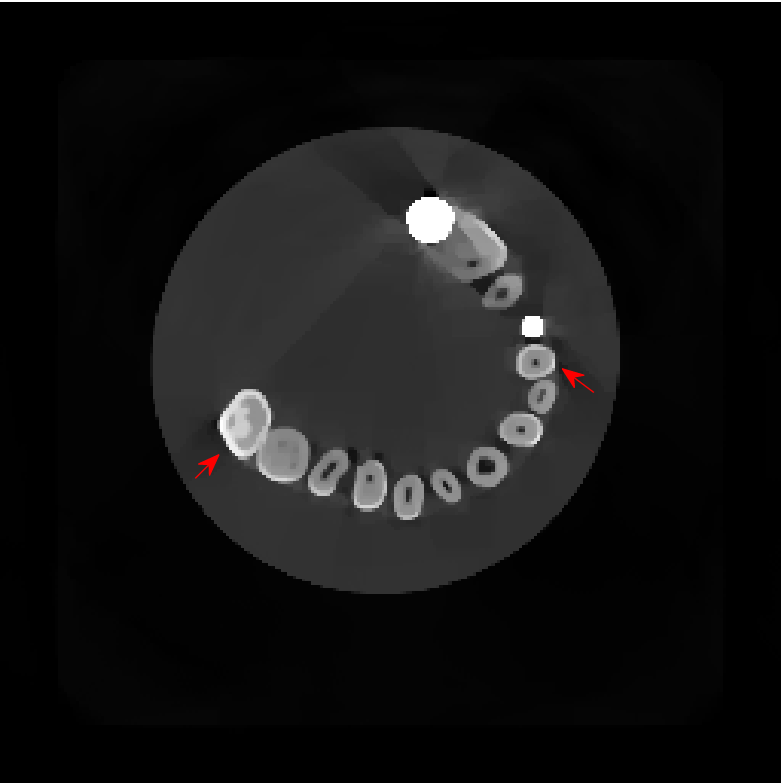}
		\caption{ FS-PDHG }
		\label{real FS}   
	\end{subfigure}
		\caption{Reconstruction results of the real dental bone (Gray value range $[0, 0.3]$).}\label{real comparison}
\end{figure}}

\section{Conclusions}\label{conclusion}
In this paper, a nonlinear optimization model combining a weighted box-constraint has been proposed for metal artifact removal, whereby non-convex regularization is incorporated into the MAR framework to enhance edge contrast. Instead of directly treating metal trace as missing data, adaptive weights are considered to promote effective reconstruction of projection structures of distinct regions and inhibit diffusion effect, constituting a high-fidelity image reconstruction model. Novel convergent algorithms including  preconditioning and full splitting  primal-dual algorithms have been developed. To validate the proposed model, both simulated and real experiments  have been conducted, demonstrating the superiority of the presented approaches.

\section*{Data availability policy}
The data that support the findings of this study are available upon request from the authors.
\section*{Acknowledgments}
The work is partially supported by the National Natural Science Foundation of China (Nos. 12271404, 62201384, 12301545, and 11871372), the Tianjin Research Innovation Project for Postgraduate Students 2022SKY257 and PHD Program 52XB2013 of Tianjin Normal University.

We would like to express our sincere gratitude to the editor and the two reviewers for their invaluable comments, which significantly enhanced the quality of the paper. Additionally, we would like to express our appreciation  to Professor Wang Chao from the School of Science, Southern University of Science and Technology for constructive suggestions.

\begin{appendices}
\setcounter{equation}{0} 
\renewcommand{\theequation}{\thesection.\arabic{equation}} 
\section{Proof of Lemma \ref{saddle point problem and original model}}
\label{proof model pre}
\begin{proof}
A direct differentiation calculus gives that
\begin{eqnarray*}
   \partial {L_\lambda }( X ) = 
   \begin{cases}
     \partial G\big( u \big) - {\rm{div}}\big( {p + \alpha q} \big) + \mathcal M_\lambda\big(u - \tilde u\big),\\
\partial {g^*}\big( q \big) + \alpha \nabla u,\\
 - \partial f_\eta ^*\big( p \big) + \nabla u,\\
\mathcal M_\lambda\big( \tilde u-u\big).  
   \end{cases}
\end{eqnarray*}
Since $X^*$ is a critical point of $L_\lambda$, then $0 \in \partial {L_\lambda }\big( {{X^*}} \big)$, i.e.
\begin{eqnarray*}
    \begin{cases}
    0 \in \partial G\big( {{u^*}} \big) - {\rm{div}}\big( {{p^*} + \alpha {q^*}} \big) + \mathcal M_\lambda\big(u^* - \tilde u^*\big),\\
0 \in \partial {g^*}\big( {{q^*}} \big) + \alpha \nabla {u^*},\\
0 \in  - \partial f_\eta ^*\big( {{p^*}} \big) + \nabla {u^*},\\
0 =\mathcal M_\lambda\big( {{\tilde u}^*} - {u^*}\big).      \end{cases}
\end{eqnarray*}
Therefore,  $0 \in \partial {L_{PD}}\big( {{u^*},{q^*},{p^*}} \big)$ with positive definite operator $\mathcal M_\lambda$.

If $\eta =0$, we conclude that
\begin{eqnarray}\label{saddle optimality condition}
\eqalign{{\rm{div}}\big( {{{p}^{*}} + \alpha {{q}^{*}}} \big) \in \partial G\big( {{u^{*}}} \big),\cr
\nabla {{ u}^{*}} \in \partial {f^*}\big( {{p^{* }}} \big),\cr
 - \alpha \nabla {{ u}^{*}}  \in \partial {g^*}\big( {{q^{*}}} \big).} 
\end{eqnarray}
Due to the closed convex functions $f$ and $g$, the conclusions hold
\begin{eqnarray*}
 {p^*} \in \partial f\big( {\nabla {u^*}} \big), \, {q^*} \in \partial g\big( { - \alpha \nabla {u^*}} \big).   
\end{eqnarray*}
Hence we can deduce that
\begin{eqnarray}
    \label{original model condition}
    \begin{aligned}
 0&={\rm{div}}\big( {{p^*} + \alpha {q^*}} \big) + \nabla ^T{p^*} + \alpha\nabla ^T{q^*} \\
    &\in \partial G\big( {{u^*}} \big) + \nabla^T\big[ {\partial f\big( {\nabla {u^*}} \big)} \big] - \alpha \nabla^T\big[ {\partial g\big( {\nabla {u^*}} \big)} \big].       
    \end{aligned}
\end{eqnarray}
Thus we get the desired lemma by (\ref{original model condition}).
\end{proof}
\section{Proof of Lemma \ref{sufficient decrease condition}}
\label{proof energy}
\setcounter{equation}{0} 
\begin{proof}
First, we consider the subproblem in Step 1 of (\ref{Pre iter}) as
\begin{eqnarray}
\label{energy}
    \begin{aligned}
{L_\lambda }&\big( {{X^k}} \big) - {L_\lambda }\big( {{u^{k + 1}},{q^k},{p^k},{u^{k - 1}}} \big)
\\\mathop  \ge \limits^{\big( {\ref{u sub}} \big)}  & \big\| {E_u^{k + 1}} \big\|_{{\mathcal M_\lambda }}^2  + \frac{1}{2}{\big\| {E_u^{k}} \big\|_{{\mathcal M_\lambda }}^2}
 - \frac{{\big\| {{u^{k + 1}} - {u^{k - 1}}} \big\|_{{\mathcal M_\lambda }}^2}}{2}.
    \end{aligned}
\end{eqnarray}

For the subproblem in Step 3 of (\ref{Pre iter}), similarly, we have
\begin{eqnarray}
\label{qqqqqqqq}
    \begin{aligned}
{L_\lambda }&\big( {{u^{k + 1}},{q^k},{p^k},{u^{k - 1}}} \big) - {L_\lambda }\big( {{u^{k + 1}},{q^{k + 1}},{p^k},{u^{k - 1}}} \big)\\
\mathop  \ge \limits^{\big( {\ref{q sub}} \big)}& \frac{1}{\tau } {\big\| {E_q^{k + 1}} \big\|^2}
+ \alpha \big\langle {\nabla \big( {{{\bar u}^{k + 1}} - {u^{k + 1}}} \big),{q^{k + 1}} - {q^k}} \big\rangle.
    \end{aligned}
\end{eqnarray}

For the subproblem in Step 4 of (\ref{Pre iter}),
\begin{eqnarray}
\label{ppppppp}
    \begin{aligned}
{L_\lambda }&\big( {{u^{k + 1}},{q^{k + 1}},{p^k},{u^{k - 1}}} \big) - {L_\lambda }\big( {{u^{k + 1}},{q^{k + 1}},{p^{k + 1}},{u^{k - 1}}} \big)  \\
\mathop\ge \limits^{\big( {\ref{y sub}} \big)}&\frac{{{{\big\| {E_p^k} \big\|}^2}}}{{2\beta }} + \big( {\frac{\eta }{2} + \frac{1}{{2\beta }}} \big){\big\| {E_p^{k + 1}} \big\|^2}+ \big\langle {\nabla \big( {{{\bar u}^k} - {u^{k + 1}}} \big),{p^{k + 1}} - {p^k}} \big\rangle \\
&- \frac{{{{\big\| {{p^{k + 1}} - {p^{k - 1}}} \big\|}^2}}}{{2\beta }}\ge\frac{\eta }{2}{\big\| {E_p^{k + 1}} \big\|^2}+ \big\langle {\nabla \big( {{{\bar u}^k} - {u^{k + 1}}} \big),{p^{k + 1}} - {p^k}} \big\rangle.
    \end{aligned}
\end{eqnarray}

Moreover,
\begin{eqnarray}
\label{Energy}
    \begin{aligned}
{L_\lambda }&\big( {{u^{k + 1}},{q^{k + 1}},{p^{k + 1}},{u^{k - 1}}} \big) - {L_\lambda }\big( {{X^{k + 1}}} \big)\\
=& \frac{{\big\| {{u^{k + 1}} - {u^{k - 1}}} \big\|_{{\mathcal M_\lambda }}^2}}{2}
- \frac{{\big\| {E_u^{k+1}} \big\|_{{\mathcal M_\lambda }}^2}}{2}.
    \end{aligned}
\end{eqnarray}

Summing up (\ref{energy}) and (\ref{Energy}), and using $2ab \le \xi {a^2} + \frac{{{b^2}}}{\xi }$, for any ${\xi } > 0$, it follows
\begin{eqnarray}
\label{energy estimation}
    \begin{aligned}
 {L_\lambda }&\big( {{X^k}} \big) - {L_\lambda }\big( {{X^{k + 1}}} \big) \\
\ge&\big\langle {\nabla \big( {{u^k} - {u^{k + 1}} + {u^k} - {u^{k - 1}}} \big),{p^{k + 1}} - {p^k}} \big\rangle +  \frac{1}{\tau } {\big\| {E_q^{k + 1}} \big\|^2} + \frac{\eta }{2}{\big\| {E_p^{k + 1}} \big\|^2}\\
&+ \alpha \big\langle {\nabla \big( {{u^{k + 1}} - {u^k}} \big),{q^{k + 1}} - {q^k}} \big\rangle + \frac{{\big\| E_u^{k+1} \big\|_{{\mathcal M_\lambda }}^2}}{2} + \frac{{\big\| E_u^{k} \big\|_{{\mathcal M_\lambda }}^2}}{2} \\
\ge&\big( {\frac{2}{\tau}-K\alpha}\big)\frac{{\big\| {E_q^{k + 1}} \big\|^2}}{2} + \big( {\eta-2K\xi_1}\big)\frac{{\big\| {E_p^{k + 1}} \big\|^2}}{2} + \frac{{\big\| E_u^{k} \big\|_{{\mathcal M_\lambda } - \frac{K}{\xi_1}\mathcal I}^2}}{2}\\
&+ \frac{{\big\| E_u^{k+1} \big\|_{{\mathcal M_\lambda } - \big(\frac{K}{\xi_1}+K\alpha\big)\mathcal I}^2}}{2}
    \end{aligned}
\end{eqnarray}
with $\xi_1 < \frac{\eta}{2K}$. Due to Condition \ref{para condition}, the operator ${\mathcal M_\lambda } - \big(\frac{K}{\xi_1}+K\alpha\big)\mathcal I$ is positive definite, such that ${\big\| E_u^{k+1} \big\|_{{\mathcal M_\lambda } - \big(\frac{K}{\xi_1}+K\alpha\big)\mathcal I}^2} \ge C_1 {\big\| {E_u^{k + 1}} \big\|^2}$ with the minimum eigenvalue $C_1$. 

Rearranging the above equation gives the  desired Lemma.
\end{proof}

\section{Proof of Lemma \ref{boundedness}}
\label{proof bound}
\begin{proof}
Obviously, the  sequence $\big\{ {{X^k}} \big\}$ is bounded. One can readily prove the nonincrease of (\ref{energy function}) by Lemma \ref{sufficient decrease condition}, and 
\begin{eqnarray}
    \begin{aligned}
{L_\lambda}&\big( {{X^k}} \big)\\
=& {G }\big( {{u^k}} \big) + {g^*}\big( {{q^k}} \big) - f_\eta ^*\big( {{p^k}} \big) + \big\langle {\nabla {u^k},{p^k} + \alpha {q^k}} \big\rangle  + \frac{{\big\| {{E_u^k} } \big\|_{{\mathcal M_\lambda }}^2}}{2}\\
 \le& {G }\big( {{u^0}} \big) + {g^*}\big( {{q^0}} \big) - f_\eta ^*\big( {{p^0}} \big) + \big\langle {\nabla {u^0},{p^0} + \alpha {q^0}} \big\rangle,     
    \end{aligned}
\end{eqnarray}
where we have introduced $u^{-1}=u^0$.

On the other hand,
\begin{eqnarray}
    \begin{aligned}
{L_\lambda }&\big( {{X^k}} \big)\\
=&{G }\big( {{u^k}} \big) + {g^*}\big( {{q^k}} \big) - f_\eta ^*\big ( {{p^k}} \big) + \big\langle {\nabla {u^k},{p^k} + \alpha {q^k}} \big\rangle  +\frac{{\big\| {{E_u^k} } \big\|_{{\mathcal M_\lambda }}^2}}{2}\\
 \ge &G\big ( {{u^k}} \big)  + {g^*}\big ( {{q^k}} \big) - f_\eta ^*\big ( {{p^k}} \big)  +\frac{{\big\| {{E_u^k}} \big\|_{{\mathcal M_\lambda }}^2}}{2}-  {\frac{{K\alpha {{\big\| {{u^k}} \big\|}^2}}}{2} + \frac{{{{\big\|K{{p^k} + \alpha {q^k}} \big\|}^2}}}{{2\alpha}}} \\
 \ge& G\big ( {{u^k}} \big ) + {g^*}\big ( {{q^k}} \big ) - f ^*\big ( {{p^k}} \big )-\frac{\eta}{2}  +\frac{{\big\| {{E_u^k} } \big\|_{{\mathcal M_\lambda }}^2}}{2}- \frac{Kc^2\alpha+K\big (1+\alpha\big )^2}{{2\alpha }}\\
 >&- \infty. \nonumber
    \end{aligned}
\end{eqnarray}
\end{proof}
\section{Proof of Lemma \ref{Subgradient upper bound}}
\label{proof subgradient}
\begin{proof}
We first estimate the upper bound of the partial derivative w.r.t. to $u$. By (\ref{u optimality condition}),
\begin{eqnarray*}
{\rm{div}}\big( {{p^k}-p^{k + 1}} \big) + \alpha {\rm{div}}\big( {{q^k}} - {q^{k + 1}} \big) \in {\partial _u}{L_\lambda }\big( {{X^{k + 1}}} \big).
\end{eqnarray*}
Hence we need to estimate the bound of the left-hand-side term of the above equation. Readily
\begin{eqnarray}
\label{u sub part}
    \begin{aligned}
{\rm{dist}}\big( {0,{\partial _u}{L_\lambda }\big( {{X^{k + 1}}} \big)} \big) \le K\big\| {{E_p^{k + 1}} } \big\| + K\alpha \big\| {{E_q^{k + 1}} } \big\|.      \end{aligned}
\end{eqnarray}

Then consider the derivative w.r.t. $q$. Similarly,
\begin{eqnarray*}
\begin{aligned}
&\alpha \big( {\nabla {u^{k + 1}} - \nabla {{\bar u}^{k + 1}}} \big) + \frac{{{q^k} - {q^{k + 1}}}}{\tau }\\
 = &\alpha \nabla \big( {{u^k} - {u^{k + 1}}} \big) + \frac{{{q^k} - {q^{k + 1}}}}{\tau } \in {\partial _q}{L_\lambda }\big( {{X^{k + 1}}} \big).
\end{aligned}
\end{eqnarray*}
We have
\begin{eqnarray}
\label{q sub part}
    \begin{aligned}
{\rm{dist}}\big( {0,{\partial _q}{L_\lambda }\big( {{X^{k + 1}}} \big)} \big) \le K\alpha \big\| {{E_u^{k + 1}}} \big\| + \frac{1}{\tau }{\big\| {{E_q^{k + 1}}} \big\|}.   
    \end{aligned}
\end{eqnarray}

For the derivative w.r.t. $p$,
\begin{eqnarray*}
    \begin{aligned}
&\nabla {u^{k + 1}} - \nabla {{\bar u}^{k + 1}} - \frac{{{p^k} - {p^{k + 1}}}}{\beta }\\
 = &\nabla \big( {{u^k} - {u^{k + 1}}} \big) - \frac{{{p^k} - {p^{k + 1}}}}{\beta } \in {\partial _p}{L_\lambda }\big( {{X^{k + 1}}} \big).
    \end{aligned}
\end{eqnarray*}
Then
\begin{eqnarray}
\label{y sub part}
    \begin{aligned}
{\rm{dist}}\big( {0,{\partial _p}{L_\lambda}\big( {{X^{k + 1}}} \big)} \big)  \le K\big\| {{E_u^{k + 1}} }\big\| + \frac{1}{\beta }{\big\| {{E_p^{k + 1}} } \big\|}.
    \end{aligned}
\end{eqnarray}

Last we consider w.r.t. ${\tilde u}$
\begin{eqnarray}
\label{uwan sub part}
    \begin{aligned}
{\rm{dist}}\big( {0,{\nabla _{\tilde u}}{L_\sigma }\big( {{X^{k + 1}}} \big)} \big) = {\big\| {\mathcal M_\lambda}\big({{u^{k + 1}} - {u^k}} \big)\big\|}\le {\big\| {{\mathcal M_\lambda }} \big\| }\big\| {{E_u^{k + 1}} } \big\|.
    \end{aligned}
\end{eqnarray}

Since 
\begin{eqnarray}
    \begin{aligned}
{\rm{dist}}\big( {0,\partial {L_\lambda }\big( {{X^{k + 1}}} \big)} \big) &\le   {\rm{dist}}\big( {0,{\partial _u}{L_\lambda }\big( {{X^{k + 1}}} \big)} \big) + {\rm{dist}}\big( {0,{\partial _q}{L_\lambda }\big( {{X^{k + 1}}} \big)} \big)\\
&  +{\rm{dist}}\big( {0,{\partial _p}{L_\lambda }\big( {{X^{k + 1}}} \big)} \big)+{\rm{dist}}\big( {0,{\nabla _{\tilde u}}{L_\lambda }\big( {{X^{k + 1}}} \big)} \big),  \nonumber      
    \end{aligned}
\end{eqnarray}
we finally conclude Lemma  \ref{Subgradient upper bound} with the given $C_2$ by summing up (\ref{u sub part})-(\ref{uwan sub part}).    
\end{proof}

\section{Proof of Lemma \ref{split saddle point problem and original model}}
\label{proof model split}
\begin{proof}
Similar to the proof of Lemma \ref{saddle point problem and original model}, if $Z^*=\big(\Lambda^*,u^*,v^*,q^*,p^*,\tilde u^*,\tilde v^*\big)$ is a critical point of the function $L_\sigma$,  $\big(u^*,v^*, q^*, p^*,\Lambda^*\big)$ is a critical point of $L$. Therefore details are omitted here.

At the same time, we have
\begin{eqnarray*}
    \begin{cases}
    0 \in \partial h\big( {{u^*}} \big) - {{\mathcal P}^T}{\Lambda ^*} - {\rm{div}}\big( {{p^*} + \alpha {q^*}} \big),\\
0 = \nabla F\big( {{v^*}} \big) + {\Lambda ^*},\\
0 \in \partial {g^*}\big( {{q^*}} \big) + \alpha \nabla {u^*},\\
0 \in  - \partial f_\eta ^*\big( {{p^*}} \big) + \nabla {u^*},\\
0 = {v^*} - \mathcal P{u^*}.    
    \end{cases}
\end{eqnarray*}
Owing to $\eta=0$, we can obtain by closed convex functions $f$ and $g$
\begin{eqnarray*}
   {p^*} \in \partial f\big( {\nabla {u^*}} \big), \, {q^*} \in \partial g\big( { - \alpha \nabla {u^*}} \big). 
\end{eqnarray*}
Meanwhile,
\begin{eqnarray*}
  \nabla F\big( {{v^*}} \big) = \frac{1}{\lambda }\big( {W \odot W \odot \big( {{v^*} - Y} \big)} \big) = \frac{1}{\lambda }\big( {W \odot W \odot \big( {\mathcal P{u^*} - Y} \big)} \big).  
\end{eqnarray*}
We further derive that 
\begin{eqnarray}
    \label{split original model condition}
    \begin{aligned}
 0=&{\rm{div}}\big( {{p^*} + \alpha {q^*}} \big) + \nabla ^T{p^*} + \alpha\nabla ^T{q^*} \\
     \in& \frac{1}{\lambda }{{\mathcal P}^T}\big( {W \odot W \odot \big( {\mathcal P{u^*} - Y} \big)} \big) + \partial h\big( {{u^*}} \big)\\
    &+ {\nabla ^T}\big[ {\partial f\big( {\nabla {u^*}} \big)} \big] - \alpha {\nabla ^T}\big[ {\partial g\big( {\nabla {u^*}} \big)} \big].      \end{aligned}
\end{eqnarray}
Consequently, we get the desired lemma by (\ref{split original model condition}).
\end{proof}
\section{Proof of Lemma \ref{pdhg sufficient decrease condition}}
\label{proof pdhg energy}
Due to the subproblems of (\ref{Pdhg iteration}), the following relations hold:
\begin{eqnarray}
    \label{uuuuu}
{\rm{div}}\big( {{p^k} + \alpha {q^k}} \big) + {{\mathcal P}^T}{\Lambda ^{k + 1}} + \frac{{{u^k} - {u^{k + 1}}}}{{{\sigma _1}}} \in \partial {h}\big( {{u^{k + 1}}} \big), \\
    \label{vvvvv}
 - {\Lambda ^{k + 1}} + \frac{{{v^k} - {v^{k + 1}}}}{{{\sigma _2}}} = \nabla F\big( {{v^{k + 1}}} \big),\\
    \label{llll}
 {\Lambda ^{k + 1}}- {\Lambda ^k}= \rho \big( {{v^k} - \mathcal P{u^k}} \big).
 \end{eqnarray}
Note that the subproblems about $q$ and $p$ are same as Algorithm \ref{Pre Pdhg}, and we omit them.
\begin{proof}
Due to the subproblems with strongly convex property, we get the following estimate:
\begin{eqnarray}
\label{u ge}
  {h}\big( {{u^k}} \big) - {h }\big( {{u^{k + 1}}} \big) \ge& \big\langle {{\rm{div}}\big( {{p^k} + \alpha {q^k}} \big) + {{\mathcal P}^T}{\Lambda ^{k + 1}} ,{u^k} - {u^{k + 1}}} \big\rangle \\
  &+ \frac{1}{{{\sigma _1}}} {\big\| { {E_u^{k + 1}}} \big\|^2},\nonumber \\
\label{v ge}
 F\big( {{v^k}} \big) - F\big( {{v^{k + 1}}} \big) \ge&- \big\langle {{\Lambda ^{k + 1}},{v^k} - {v^{k + 1}}} \big\rangle  + \frac{{{{\big\| {{v^k} - {v^{k + 1}}} \big\|}^2}}}{{{\sigma _2}}}.
  \end{eqnarray}

First, we consider the subproblem in Step 1 of (\ref{Pdhg iteration}) as
\begin{eqnarray}
\label{energy split}
    \begin{aligned}
{L_\sigma }&\big( {{Z^k}} \big) - {L_\sigma }\big( {{\Lambda ^{k + 1}},{u^k},{v^k},{q^k},{p^k},{u^{k - 1}},{v^{k - 1}}} \big)\\
 =& \big\langle {{\Lambda ^k} - {\Lambda ^{k + 1}},{v^k} - P{u^k}} \big\rangle 
 =  - \frac{1}{\rho }{\big\| {E_\Lambda^{k+1}} \big\|^2}
    \end{aligned}
\end{eqnarray}
By (\ref{vvvvv}), we have 
\begin{eqnarray}
    \begin{aligned}
&{\big\| {E_\Lambda^{k+1}} \big\|^2}\\ = &{\big\| {\nabla F\big( {{v^k}} \big) - \nabla F\big( {{v^{k + 1}}} \big) + \frac{{{v^k} - {v^{k + 1}}}}{{{\sigma _2}}} - \frac{{{v^{k - 1}} - {v^k}}}{{{\sigma _2}}}} \big\|^2}\\
 \le &2{\Big\| {\nabla F\big( {{v^k}} \big) - \nabla F\big( {{v^{k + 1}}} \big) + \frac{{{v^k} - {v^{k + 1}}}}{{{\sigma _2}}}} \Big\|^2} +\frac{2}{\sigma_2^2} {\big\| {{{E_v^k}}}\big\|^2}\\
 \le &4{\big\| {\nabla F\big( {{v^k}} \big) - \nabla F\big( {{v^{k + 1}}} \big)} \big\|^2} + \frac{4}{{{\sigma_2^2}}}{\big\| {E_v^{k + 1}} \big\|^2} +\frac{2}{\sigma_2^2} {\big\| {{{E_v^k}}}\big\|^2}\\
\le& 4\big( {\frac{1}{{{\lambda ^2}}}\big\| {W} \big\|_{\rm{max}} ^4 + \frac{1}{{\sigma _2^2}}} \big){\big\| {{E_v^{k + 1}}} \big\|^2} +\frac{2}{\sigma_2^2} {\big\| {{{E_v^k}}}\big\|^2}.        
    \end{aligned}
\end{eqnarray}

Then we estimate the subproblem in Step 2 of (\ref{Pdhg iteration}) as
\begin{eqnarray}
    \begin{aligned}
{L_\sigma }&\big({{\Lambda ^{k + 1}},{u^k},{v^k},{q^k},{p^k},{u^{k - 1}},{v^{k - 1}}} \big)- {L_\sigma }\big( {{\Lambda ^{k + 1}},{u^{k + 1}},{v^k},{q^k},{p^k},{u^{k}},{v^{k - 1}}} \big)\\
 \ge&{\frac{1}{{{2\sigma _1}}} } {\big\| {E_u^{k+1}} \big\|^2} + \frac{1}{{2{\sigma _1}}}{\big\| {{E_u^{k}}} \big\|^2}. 
    \end{aligned}
\end{eqnarray}

For the subproblem in Step 4 of (\ref{Pdhg iteration}), similarly, we have
\begin{eqnarray}\small
    \begin{aligned}
{L_\sigma }&\big( {{\Lambda ^{k + 1}},{u^{k+1}},{v^k},{q^k},{p^k},{u^{k}},{v^{k - 1}}} \big)- {L_\sigma }\big( {{\Lambda ^{k + 1}},{u^{k + 1}},{v^{k+1}},{q^k},{p^k},{u^{k}},{v^{k}}} \big) \\
\ge& \frac{1}{{{2\sigma _2}}}{\big\| {{E_v^{k+1}}} \big\|^2} + \frac{1}{{2{\sigma _2}}}{\big\| {{E_v^{k}}} \big\|^2}.
    \end{aligned}
\end{eqnarray}

Regarding the subproblems of $q$ and $p$, (\ref{qqqqqqqq}) and (\ref{ppppppp}) give 
\begin{eqnarray}\small
\begin{aligned}
   {L_\sigma }&\big( {{\Lambda ^{k + 1}},{u^{k+1}},{v^{k+1}},{q^k},{p^k},{u^{k}},{v^{k}}} \big) 
- {L_\sigma }\big( {{\Lambda ^{k + 1}},{u^{k + 1}},{v^{k+1}},{q^{k+1}},{p^k},{u^{k}},{v^{k}}} \big)\\
\ge&\frac{1}{{{\tau}}}{\big\| {{E_q^{k+1}}} \big\|^2} +\alpha \big\langle {\nabla {u^{k + 1}} - \nabla {u^k},{q^{k + 1}} - {q^k}} \big\rangle,\\ 
\label{Energy split}
{L_\sigma }&\big( {{\Lambda ^{k + 1}},{u^{k+1}},{v^{k+1}},{q^{k+1}},{p^k},{u^{k}},{v^{k}}} \big) - {L_\sigma }\big( {Z^{k+1}} \big)\\
\ge&\frac{\eta}{{{2}}}{\big\| {{E_p^{k+1}}} \big\|^2}+ \big\langle {\nabla ({u^{k + 1}} - {u^k}+{u^k}-{u^{k+1}}),{p^{k + 1}} - {p^k}} \big\rangle. 
\end{aligned}
\end{eqnarray}

Summing up (\ref{energy split}) and (\ref{Energy split}), it follows
\begin{eqnarray}
\label{energy split estimation}
    \begin{aligned}
{L_\sigma }&\big( {{Z^k}} \big) - {L_\sigma }\big( {{Z^{k + 1}}} \big)\\
\ge &\big( {1-\frac{K\sigma_1}{\xi_2}}\big)\frac{{{\big\| {E_u^{k+1}} \big\|}^2}+{{\big\| {E_u^k} \big\|}^2}}{2\sigma_1}+\big( {\frac{2}{\tau}-K\alpha}\big)\frac{{\big\| {E_q^{k + 1}} \big\|^2}}{2} \\
&-{K\alpha}\frac{{\big\| {E_u^{k + 1}} \big\|^2}}{2}+ \big( {\eta-2K\xi_2}\big)\frac{{\big\| {E_p^{k + 1}} \big\|^2}}{2} +\big(\frac{1}{{2{\sigma _2}}} - \frac{2}{{\rho\sigma _2^2}}\big){\big\| {E_v^{k }} \big\|^2}\\
&+\Big(\frac{1}{{2{\sigma _2}}} - \frac{4}{\rho}\big( {\frac{1}{{{\lambda ^2}}}\big\| {W} \big\|_{\rm{max}} ^4 + \frac{1}{{\sigma _2^2}}} \big)\Big){\big\| {E_v^{k + 1}} \big\|^2}.
    \end{aligned}
\end{eqnarray}
Rearranging the above equation and setting 
\begin{eqnarray*}
    \begin{aligned}
     \tilde C := \min \left\{\frac{1}{{2\sigma_1 }} - \frac{K}{{2{\xi _2}}}-\frac{K\alpha}{2},\frac{1}{\tau } - \frac{{K\alpha }}{{2}},\frac{\eta }{2} - K{\xi _2},\right.\\
     \left.\frac{1}{{2{\sigma _2}}} - 4\Big( {\frac{1}{{{\lambda ^2}}}\big\| {W} \big\|_{\rm{max}} ^4 + \frac{1}{{\sigma _2^2}}} \Big) \right\},    \end{aligned}
\end{eqnarray*}
readily we obtain
\begin{eqnarray*}
    {L_\sigma }\big( {{Z^k}} \big) - {L_\sigma }\big( {{Z^{k + 1}}} \big) \ge \tilde{C} {\big\| {{Z^{k + 1}} - {Z^k}} \big\|^2}.
\end{eqnarray*}
\end{proof}
\end{appendices}

\section*{References}
\bibliographystyle{IEEEtran}
\bibliography{refs}
\end{document}